\documentclass[10pt]{elsarticle}
\usepackage{lineno,hyperref}
\usepackage{booktabs}
\usepackage{multirow}
\usepackage{amsthm}
\usepackage{hyperref}
\usepackage{crossreftools}
\usepackage{physics}
\usepackage{graphicx}
\usepackage{subfigure}
\usepackage{float}
\usepackage{bm}
\usepackage{lipsum}
\usepackage{listings} 
\usepackage{color}
\usepackage{amsfonts,enumerate,amsmath,amssymb}
\def\qed{\strut\hfill $\Box$}

\newtheorem{theorem}{Theorem}[section]
\newtheorem{lemma}[theorem]{Lemma}
\newtheorem{remark}[theorem]{Remark}
\newtheorem{proposition}[theorem]{Proposition}
\newtheorem{definition}[theorem]{Definition}
\newtheorem{example}{Example}

\def\para#1{\vskip .4\baselineskip\noindent{\bf #1}}
\newcommand{\vertiii}[1]{{\left\vert\kern-0.25ex\left\vert\kern-0.25ex\left\vert #1 
\right\vert\kern-0.25ex\right\vert\kern-0.25ex\right\vert}}
\def\para#1{\vskip .4\baselineskip\noindent{\bf #1}}
\numberwithin{equation}{section}

\def\f{\frac}

\begin{document}
\begin{frontmatter} 
\title{Central limit theorem for slow-fast system under mixed fractional Brownian motion}
\author[mymainaddress,mysecondaddress]{Xiaoyu Yang}
\ead{yangxiaoyu@yahoo.com}

\author[mymainaddress,mythirdaddress]{Yong Xu\corref{mycorrespondingauthor}}\cortext[mycorrespondingauthor]{Corresponding author}\ead{hsux3@nwpu.edu.cn}

\address[mymainaddress]{School of Mathematics and Statistics, Northwestern Polytechnical University, Xi'an, 710072, China}
\address[mysecondaddress]{Faculty of Mathematics, Kyushu University, Fukuoka, 8190395, Japan}
\address[mythirdaddress]{MOE Key Laboratory of Complexity Science in Aerospace, Northwestern Polytechnical University, Xi'an, 710072, China}

\begin{abstract}
This work considers a type of slow-fast system, where the slow component is driven by fractional Brownian motion (FBM) with \(H > 1/2\) and the fast component is a Markovian stationary process. Our solution mapping is defined based on the Young-Wiener sense, which is constructed via the stochastic sewing lemma. Then, we aim to show the fluctuation from the averaging limit by applying the Poisson PDE method. Unlike the case of standard Brownian motion, the Poisson PDE method must be developed to a non-Markovian fractional setting. The first part addresses the central limit theorem problem for a slow-fast system under small FBM, which is fully coupled with the fast varying one. Firstly, a Wiener-Young-It\^o formula is constructed for the Wiener-Young-It\^o integral. The behavior of the deviation component is related to solutions of Poisson PDEs associated with the Laplace operator of the fast process. The fast one is assumed to satisfy ``stricter” H\"older conditions related to the time scale parameter. The tightness is then derived through the H\"older semi-norm of the fast one, the properties of the Poisson PDE solution, and a Gronwall-type result for a linear Young differential equation. The weak limit shown includes an extra Gaussian process. The second part is dedicated to the problem under a general FBM. In contrast to the former one, here the FBM integral term will be more difficult to bound due to the coupling between the regularity of the Poisson PDE solution and dependence on the unbounded H\"older seminorm of the fast one. 

\vskip 0.08in 
\noindent{\bf AMS Math Classification.}
{60F05, 60H10, 60G22, 35J05}.
\vskip 0.08in 
\end{abstract} 
\end{frontmatter}

\tableofcontents

\section{Introduction}
In this article we study the deviation of the solutions to the slow one of the following system
\begin{align}\label{System}
\begin{cases}
dx^{\epsilon}_t = g(x^{\epsilon}_t, y^{\epsilon}_t)dt +a(\epsilon) f(x^{\epsilon}_t, y^{\epsilon}_t)dB^H_{t},\\[4pt]
dy^{\epsilon}_t = \frac{1}{\epsilon} b(x^{\epsilon}_t, y^{\epsilon}_t)dt + \frac{1}{\sqrt{\epsilon}}\sigma(x^{\epsilon}_t, y^{\epsilon}_t)dw_{t}.
\end{cases}
\end{align}
from its effective limit. Here $B_t^H$ denotes a fractional Brownian motion (FBM) in \(\mathbb{R}^d\) with similarity exponent $H>\f 12$ and $w$ denotes a standard  Brownian motion (BM) in \(\mathbb{R}^e\).
Moreover, \(y^{\epsilon}\) is stationary It\^o process which will be assumed owing the unique invariant measure \(\mu\). 
The stochastic process \(B^H\) and \(y\) is set on a standard probability space \((\Omega, \mathcal{F}, \mathbb{P})\).
Assume that the nonlinear coefficients \(f\)  is of \(BC^2\), \(g\in L(\mathbb{R}^d,\mathbb{R}^n)\) is \(BC^2\) and globally bounded.
In contrast, when $f\equiv 0$, a central limit theorem is available that identifies the fluctuation. However, with the addition of a simple additive non-Markovian noise, a completely different phenomenon emerges.

As \(\epsilon\) tends to zero, the above system converges to the effective averaged system under the \(\alpha\)-H\"older topology \((\alpha<H)\) \cite{Hairer-Li}, it is given by: 
\begin{equation}\label{averaged}
\dd \bar x_t = \bar g(x_t) \dd t+a \bar f(\bar x_t) \dd B^H_t,\ \bar x_0 = x.
\end{equation}
Here, \(\bar g = \int_{\mathbb{R}^{m}} g(x,y)\mu^x(\dd y)\)\(y\) where is \(\mu(\dd y)\) is the invariant probability measure of the stationary fast It\^o process with fixed \(x\). The averaging principle is a well-known topic, especially in the context of Markovian ordinary differential equations or stochastic differential equations \cite{Li_2008}.
When the driving noise is fractional Brownian motion with $H>1/2$, the averaging principle for stochastic ordinary/partial systems has been studied in \cite{Pei2023,Pei2024,bourguin-Gailus,LI-Gao}. For the case $H<1/2$, however, there are relatively few results on the averaging principle; we refer to \cite{PEI2021202, Yuzuru2025, inahama2025averaging, LI2025104683,li2025averaging}, where the problem is treated within the framework of rough path theory \cite{lyons-Qian,Friz-Victoir,Friz-Hairer}.

We now aim to study the deviation of the solutions to \eqref{System} from their effective limit \eqref{averaged}. To this end, we set
$$z_t^\epsilon =\f 1{\sqrt \epsilon } (x_t^\epsilon -\bar x_t),$$
then directly
$$dz_t^\epsilon =
\f 1{\sqrt \epsilon }(g(x_t^\epsilon,y_t^\epsilon)-\bar g(\bar x_t))\dd t+\f {a(\epsilon)}{\sqrt \epsilon }(f(x_t^\epsilon,y_t^\epsilon)-\bar f(\bar x_t))\dd B_t^H.$$
The asymptotic dynamics of \(z\) is also called the central limit theorem.
While by the functional limit theorem for Markovian ordinary differential equation, that is \(f=0\),
$$ \f 1{\sqrt \epsilon }\int_0^t (g(\bar x_s,y_{\frac{s}{\epsilon}})-\bar g(\bar x_s))ds$$
converges weakly to a Wiener process, so this classical setting which falls within the theory of diffusion creation, which has been studied a lot, see eg \cite{kipnis1986central} and the book \cite{komorowski2012fluctuations}.
When the fast one is Gaussian, the CLT and homogenization for a wide variety of situations are also solved via fourth moments, as discussed in references \cite{nourdin2012normal,nourdin2016multivariate} that based on the Malliavin calculus \cite{Nualart} and also see \cite{hairer2004periodic,BAILLEUL20174894}.
With the development of the rough path theory \cite{Friz-Victoir,Friz-Hairer}, it is worthy to state that the functional limit theorem is obtained by the rough path theory, where the fast one is fractional Ornstein-Uhlenbeck process driven by short and long-range dependent fractional noise \cite{gehringer2020diffusive,gehringer2022functional}.
Other types of fast processes were also considered in the functional limit theorem \cite{Gehring-Li,Djurdjevac-Helena,KELLY20174063,engel2021homogenization, GU20121069}. 
The homogenization are also given \cite{BAILLEUL20174894} by the rough flows \cite{Bailleul2019} which encompasses the theory of rough differential equations, and unifies them with the theory of stochastic flows. 
The method of martingale approximation also works for the central limit theorem \cite{jacod2003limit,Kelly-Melbourne}, then, it can be extended to the Markovian It\^0 process (\(f \neq 0\) and \(H=1/2\)) \cite{WANG2013822} and the Poisson jump process \cite{YANG2022107897}.
The CLT results via the Poisson equation method for recent work could see \cite{li2023functional,Rockner-Xie,ROCKNER-Xie23,Bourguin-Spiliopoulos}, also see diffusive approximation \cite{Bezemek-Spiliopoulos,Li-Wu-Xie,HONG2025405,sun2025diffusion}, and more about it could refer to \cite{Pardoux2001,Pardoux2003,Pardoux2005}. 
Just a quick note for some application of Poisson equation could be found at the averaging principle \cite{rockner2021strong}.

However, there has only been a little work on the CLT for fractional slow systems (see, for example \cite{Hairer:22,li2025fluctuations}). Firstly, we address the CLT problem for a stochastic system governed by a small FBM \eqref{slowfast}, that is $a(\epsilon)=\sqrt{\epsilon}$. Specifically, we assume that the slow and fast variable fully depends on each other. The integral with respect to the FBM is therefore in the Wiener-Young sense \cite{Hairer-Li} avoiding the difficulty caused by the unbounded H\"older semi-norm of the fast variable with respect to the time-scale parameter. Here, other methods are inefficient here, so we will utilise the Poisson equation method.
Note that to apply the Poisson equation method, the Young-It\^o formula for the stochastic dynamical system with mixed FBM in the mixed Young-It\^o sense is required. 
Moreover, the H\"older seminorm of the fast variable is assumed to satisfy ``stricter" H\"older conditions related to the time scale parameter.
Our second work considers a slow-fast system, where the slow component is driven by fractional Brownian motion and only the drift term depends on the fast one, but the fast one is a Markovian stationary process that fully depends on the slow variable. When $a(\epsilon)=1$, according to the \cite{Hairer:22,li2025fluctuations}, which states that \(\int_{0}^{t}\epsilon^{\frac{1}{2}-H}f(x_s^\epsilon, y^\epsilon_s) \dd B^H_t\) converges to an extra Gaussian process, so if \(f:=f(x,y)\) in \eqref{System}, the deviation \(z^\epsilon\) will blow up as \(\epsilon\) tends to zero. 
Then, the CLT problem is converted into a weak convergence problem, which is defined by the Laplace operator of the fast process. Initially, the tightness is derived through Wiener-Young integral computation, the H\"older seminorm of the fast variable, the properties of solutions to Poisson equations, and a Gronwall-type result for Young differential equations.
The limit process involves an additional Wiener process that is related to solving the Poisson equation.
The second part of this work is dedicated to addressing the CLT for the system derived above, that is the qualitative analysis of the deviation \(z^\epsilon\) when the small parameter \(\epsilon\) tends to zero.
In contrast to the former case, the FBM  integral here is more difficult to bound, owing to the coupling between the regularity of the solution to the Poisson equation and its dependence on the fast variable. This difficulty is circumvented by a delicate combination of the Young--Wiener--It\^o integral, the stronger H\"older regularity of the slow variable, and the ergodicity of the fast variable.
Finally, it is demonstrated that the limit satisfies with the linear Young-It\^o differential equation.

Throughout this paper, the symbols $c$, $C$, $c_1$, $C_1$, $\cdots$ will denote certain positive constants that may vary from line to line. The dependence of the constant on parameters will be explicit if necessary. 
\subsection{Notations}
\begin{itemize}
\item Assume \(\mathcal{V}\) and \(\mathcal{W}\) are both Euclidean spaces.
\item Denote the probability space by \((\Omega, \mathcal{F}, \mathbb{P})\) and \(\|\cdot\|_p\) 
is the norm in \(L^p(\Omega;\mathbb{P})\).
\item \(|\cdot|_{\infty}\) and \(|\cdot|_{Lip}\) denote the supremum norm and the minimal Lipschitz constant.
\item For a one-parameter process \(x_\cdot : [0,T] \to \mathcal{V}\), we denote \(x_{s,t} := x_t-x_s\) for any \(0\leq s\leq t\leq T\).
\item For \(0\leq s\leq u\leq t\leq T\) and two-parameter process \(A\), we denote \[\delta A_{s,u,t} := A_{s,t} -A_{s,u} -A_{u,t}. \]
For \(0<\eta< 1\), set 
$$\|A\|_{\eta,p}:=\sup_{s<t}\frac{\|A\|_{p}}{|t-s|^\eta}, \quad \vertiii{A}_{\eta,p}:=\sup_{s<t}\frac{\|\mathbf{E}[\delta A_{s,u,t}|\mathcal{F}_s]\|_{p}}{|t-s|^\alpha}.$$
\item \(L(\mathcal{W},\mathcal{V})\) denotes the set of bounded linear maps from \(\mathcal{W}\) to \(\mathcal{V}\).
\item For \(k\in \mathbb{N}\), denote \(BC^k\) the set of \(C^k\)-bounded functions whose derivatives up to order \(k-\) are also bounded.
\item \(H_\eta^p=\{A_{s,t}\in L^p(\Omega,\mathcal{F}_t,\mathbb{P}): \|A\|_{\eta,p}<\infty\}\) and \(\bar H_\eta^p=\{A_{s,t}: \vertiii{A}_{\eta,p}<\infty\}\).
\item Set $\eta \in(0,1]$. For a continuous path \(h:[0,T]\mapsto\mathcal{V}\), define the \(\eta-\)H\"older semi-norm by
\(\|h\|_{\eta}:=\sup_{0 \leq s<t \leq T} \frac{|h_t-h_s|}{(t-s)^\eta}<\infty\).
The set of all \(\eta\)-H\"older continuous path is denoted by \(\mathcal{C}^{\eta}([0,T],\mathcal{V})\).
\item \(\mathcal{B}_{\alpha,p}=\{x_t\in \mathcal{F}_t:\delta x_{s,t}\in H_\alpha^p\}\), \(\mathcal{B}_{\alpha,p}\subset L_p(\Omega,\mathcal{C}^\gamma)\) (up to modification) for \(\gamma<\alpha-\frac{1}{p}\). 
\item For a two-parameter process \(u(x,y): (\mathcal{V}\times \mathcal{W})\), define \(C_b^{k_1,k_2+\gamma}\) the space of functions that \(u(\cdot,y)\) is \(BC^{k_1}\) for any fixed \(y\in\mathcal{W}\) and \(u(x,\cdot)\) is \(BC^{k_2}\) for any fixed \(x\in\mathcal{V}\), meanwhile, \(k_2\) derivates of \(u(x,\cdot)\) are \(\gamma\)-H\"older continuous in \(x\).
\item For \(0<\eta< 1\) and \(h\in \mathcal{C}([0,T], \mathcal{V})\), set \(|h|_{-\kappa}:=\sup_{0\leq s<t\leq T} |t-s|^{\kappa-1}|\int_{s}^{t}h_rdr|\).
\end{itemize}
\section{Fractional Brownian motion}
In this section, we introduce the FBM and standard BM. A $\mathbb{R}^{d}$-valued continuous stochastic process $B^H_t=(B_t^{H,1},B_t^{H,2},\cdots,B_t^{H,d})_{t\in[0,T]}$ starting from 0 is called an FBM if it is a centred Gaussian process, satisfying that
$$\mathbb{E}\big[B_{t}^{H} B_{s}^{H}\big]=\frac{1}{2}\left[t^{2 H}+s^{2 H}-|t-s|^{2 H}\right]\times I_{d}, \quad(0\le s\le t \le T),$$
where $I_{d}$ denotes the identity matrix in $\mathbb{R}^{d\times d}$. 
Then, it is straightforward to see that 
$$\mathbb{E}\big[(B_{t}^{H}-B_{s}^{H})^{2}\big]=|t-s|^{2 H} \times I_{d},\quad (0\le s\le t \le T).$$
According to the
Kolmogorov's continuity criterion, 
the trajectories of $B^H$ are $H'$-H\"older continuous ($H'\in(0,H)$) and have $\lfloor 1 / H \rfloor<{p}<\lfloor 1 / H \rfloor+1$-variation almost surely. 

The FBM $B^H=(B_t^H)_{t\in [0,T]}$ can typically be constructed using a BM $B$,
\begin{equation}\label{vol-rep}
B_t^H :=\int_0^T K_H(t, s) d B_s, \qquad t \in[0,T].
\end{equation}
Here, we set for all $0\le s \le t \le T$,
$
K_H(t, s) :=k_H(t, s) 1_{[0, t]}(s) 
$
with
$$
k_H(t, s) :=\frac{c_H}{\Gamma\left(H+\frac{1}{2}\right)}(t-s)^{H-\frac{1}{2}} F\left(H-\frac{1}{2}, \frac{1}{2}-H, H+\frac{1}{2} ; 1-\frac{t}{s}\right),
$$
where
$c_H=\left[\frac{2 H \Gamma\left(\frac{3}{2}-H\right) \Gamma\left(H+\frac{1}{2}\right)}{\Gamma(2-2 H)}\right]^{1 / 2}$, $\Gamma$ is a gamma function and $F$ is the Gauss hypergeometric function (see e.g. \cite{1999Decreusefond})

We also consider a standard ${e}$-dimensional BM $(W_t)_{t\in [0,T]}$. Throughout, $(W_t)$ and $(B^H_t)$ are assumed to be independent. We set the filtration $\mathcal{F}_t=\mathcal{G}_t\vee \mathcal{\hat G}_t$ where
$$\mathcal{G}_t=\sigma\{B^H_u-B^H_r, r\le u\le t\}\quad \mathcal{\hat G}_t=\sigma\{W_u-W_r, r\le u\le t\}.$$
Note that $\mathcal{G}_t=\sigma\{B_u-B_r, r\le u\le t\}$ due to the representation \eqref{vol-rep}.

In particular, we introduce the Mandel-Van Ness representation of FBM $B^H$ by a Wiener process $B$ \cite{mandelbrot1968fractional}. Precisely, when Hurst parameter \(H\in(\frac{1}{2},1)\),
\begin{eqnarray}\label{FBM_MV}
B^H_r-B^H_u&=&\int_{-\infty}^u((r-v)^{H-\frac{1}{2}}-(u-v)^{H-\frac{1}{2}})\dd B_v+\int_u^r(r-v)^{\frac{1}{2}}\dd B_v\\
&=:&\bar B_r^{H,u}+\tilde B_r^{H,u}.
\end{eqnarray}
Recall that the filtration $\mathcal{G}_t$ generated by the increments of $B$ coincides with that generated by $B^H$ coincides with it, then $\bar B_t^{H,u}$ is $\mathcal{G}_u$-measurable and smooth on $t\in(u,\infty)$ but $\tilde B_t^{H,u}$ is independent of $\mathcal{G}_u$.

Then we introduce the Young-Wiener integral with respect to FBM $B^H$ with $H\in(\frac{1}{2},1)$. For a continuous measurable function $F:\mathbb{R}\times \mathbb{R}^m\mapsto \mathbb{R}^{m\times d}$ and an $\mathcal{G}_t$-measurable process $x_t$ for all $t\in[0,T]$, 
\begin{eqnarray}\label{integral_FBM}
A_{s,t}=\int_0^t F(r,x_r)\dd B^H_r=\int_0^t F(r,x_r)\dd \bar B_r^{H,s}+\int_0^t F(r,x_r)\dd \tilde B_r^{H,s}.
\end{eqnarray}
The first integral is considered as Riemann-Stieltjes integral with respect to the smooth process $\bar B^{H,s}$ and the second term is interpreted as a Wiener integral with respect to the Gaussian process $\tilde B^{H,s}$ separately. Note that $x_s$ is $\mathcal{G}_s$-measurable and independent of $\tilde B^{H,s}$. More details could be referred to \cite{Hairer-Li}.

Then, we introduce the stochatic sewing lemma whose details of the proof could refer to  \cite[Theorem 2.1]{Khoa2020}.
\begin{lemma}\label{sewing}
(Stochastic sewing Lemma)Let \(p\ge 2\) and \(\alpha>\frac{1}{2}\). One has \(A\in H_\eta^p\cap \bar H_{\bar \eta}^p\) with $\eta>1/2$ and $\bar \eta>1$. Then, for every $t\in[0,T]$, the limit
\begin{eqnarray}\label{Ist}
I_{s,t}:=\lim_{| \mathcal{P}|\to 0} \sum_{[u,v]\in \mathcal{P}} A_{u,v}
\end{eqnarray}
exists in $L_p$-sense. Here $ \mathcal{P}$ is a partition of time interval $[s,t]$. Moreover, there exists a constant $C>0$ such that 
$$\|I_{s,t}(A)\|_{p}\le C \big(\vertiii{A}_{\eta,p}|t-s|^{\bar \eta}+\|A\|_{\eta,p}|t-s|^{ \eta}\big)$$
$$\|\mathbb{E}\big(I_{s,t}(A)-A_{s,t}|\mathcal{F}_t\big)\|_p\le C\vertiii{A}_{\eta,p}|t-s|^{\bar \eta}.$$
Furthermore, if $\|\mathbb{E}\big(A_{s,t}|\mathcal{F}_t\big)\|_p \le C\vertiii{A}_{\eta,p}|t-s|^{\bar \eta}$, then $I(A)\equiv 0$. 
\end{lemma}

\begin{lemma}\label{Integral-dB}
(\cite[Lemma 3.10]{Hairer-Li})Let \(p\ge 2\) and \(\alpha>\frac{1}{2}\). Assume that \(x_{\cdot}\in \mathcal{B}_{\alpha,p}\), let \(f\in \mathcal{C}^{-\kappa,\gamma}\)(deterministic) for some \(\kappa, \gamma\) such that \(\eta=H-\kappa>\frac{1}{2}\) and \(\bar \eta=H-\kappa+\gamma\alpha\) and define
$A_{s,t}:=\int_{s}^{t}f(r, x_s)d B_{r}^H$ where the integral is defined as a conditional Wiener integral. Then \(A\in H_\eta^p\cap \bar H_{\bar \eta}^p\) and 
\begin{equation}\label{f(r,x_r)}
\|\int_{s}^{t}f(r, x_r)d B_{r}^H\|_{p} \lesssim |f|_{-\kappa,p}(|t-s|^{H-\kappa}+\|x\|_{\alpha}(t-s)^{(\bar \eta -\alpha)})
\end{equation}
and 
\begin{equation}\label{f(r,x_r)-f(r,x_s)}
\|\mathbb{E}\big(\int_{s}^{t}(f(r, x_r)-f(r, x_s))d B_{r}^H|\mathcal{F}_s\big)\|_{p} \lesssim |f|_{-\kappa,p}\|x\|_{\alpha}(t-s)^{(\bar \eta -\alpha)}.
\end{equation}
\end{lemma}

\begin{remark}\label{Young}
(\cite[Lemma 3.12]{Hairer-Li})
Further assume there exists \(\delta>0\) such that \(\delta+H>1\) and \(\sup_x\sup_{|t-s|\le 1}|t-s|^{-\delta}|f(t,x)-f(s,x)|<\infty\), the integral \(\int_{s}^{t}f(r, x_r)d B_{r}^H\) coincides with the usual Young integral.
\end{remark}

\section{Poisson equation}
Assume that the  process $Y^x$ is a stationary Markovian  and takes values in \(\mathbb{R}^{m}\) with fixed parameter $x$, that is, for any time \(t\in\mathbb{R}^{+}\),  that is
\begin{eqnarray}\label{Xi}
\dd Y^x_{t}= b(x, Y^x_{t})\dd t + \sigma(x, Y^x_{t})\dd W_t, Y^x_{0}=y.
\end{eqnarray}
Now the Markov transition semigroup denoted by \(\mathcal{P}_t^x\) is associated with the generator \(\mathcal{L}^x\),
\begin{eqnarray}\label{y_generator}
\mathcal{L}^x(x,y)= \frac{1}{2} a_{i,j}(x,y)\frac{\partial^2}{\partial y_i \partial y_j}+b_{i}(x,y)\frac{\partial}{\partial y_i}
\end{eqnarray}
where \(a=\sigma\sigma*\). Define that \(p_t(y,y';x)\) be the transition density function of \(Y_t^x\) which is also the fundamental solution for the operator \(\mathcal{L}^x\) with given \(x\in \mathbb{R}^n\),
then for any measurable function \(F\), set \(\mathcal{P}_t^x F(x,y):=\tilde{\mathbb{E}}[F(x, Y_t^x)]=\int_{\mathbb{R}^m}{p}_t^x(y,y';x)F(x,y')dy'\).

Furthermore, we assume the following holds.
\begin{itemize}
\item[H1.] Assume that \(b, \sigma \in C_b^{2,2+\gamma}\) for \(\gamma\in(0,1]\).
\item[H2.] Assume that \(a:=\sigma\sigma*\) is uniformly elliptic, that is, there exists a constant \(\lambda>1\) such that for any \(x\in \mathbb{R}^n\),
$$ \lambda^{-1}|\xi|^2 \le |a (x,y)\xi|\le \lambda |\xi|^2, \quad \forall \xi\in \mathbb{R}^{m}.$$
\item[H3.] There exist constants \(C>0\) and \(\beta_1>0\) such that for all \((x,y)\in \mathbb{R}^n\times\mathbb{R}^m\), \(\langle y, b(x,y)\rangle=-\beta_1|y|^2+C|x|^2+C\).
\end{itemize}

\begin{remark}\label{ergodicity}
Assume that \(b,\sigma\) are Lipschitz continuous and \textbf{H2}-\textbf{H3} hold. Then the transition semigroup \(\mathcal{P}_t^x\) related to \eqref{Xi} admits a unique invariant probability measure \(\mu^x\) on $\mathbb{R}^m$.
We have that for any fixed \(x_1,x_2\in \mathbb{R}^n\) and \(y\in \mathbb{R}^{m}\), any bounded measurable function \(u(y): \mathbb{R}^m \mapsto \mathbb{R}\),
$$\|\mathcal{P}_t^{x_1} u(y)-\mathcal{P}_t^{x_2} u(y)\|_{\infty}\le C |x_1-x_2|\|u\|_{Lip}.$$
We also have that for any fixed \(x\in \mathbb{R}^n\) and \(y\in \mathbb{R}^{m}\), any bounded measurable function \(u(y): \mathbb{R}^m \mapsto \mathbb{R}\),
$$\|\mathcal{P}_t^x u(y)-\int_{\mathbb{R}^{m}} u(y)\mu^x(dy)\|_{\infty}\le C e^{-ct}\|u\|_{Lip}(1+|x|+|y|)$$
and 
$$\|\mathcal{P}_t^x u(y)\|_{Lip}\le C e^{-ct}\|u\|_{Lip}(1+|x|+|y|).$$
\end{remark}

\begin{remark}\label{derivateofy}
Assume that \textbf{H1}-\textbf{H3} hold. Then, for any \(t\in\mathbb{R}^{+}\) and \(\xi\in \mathbb{R}^m\), \(\|D_xY_{t}^x\cdot \xi\|_{L_2}\leq C|\xi|^2\). Furthermore, for any \(t\in\mathbb{R}^{+}\) and \(\xi,\xi'\in \mathbb{R}^m\), 
\begin{eqnarray}\label{partial_y}
\|D_yY_{t}^x\cdot \xi\|_{L_2}\leq Ce^{-ct}|\xi|^2,\quad \|D^2_{yy}Y_{s,t}^x\cdot (\xi,\xi')\|_{L_2}\leq Ce^{-ct}(|\xi|^2+|\xi'|^2).
\end{eqnarray}
The above results could be obtained by \cite[Lemma 3.6]{rockner2021strong} or \cite[Proposition 3.1]{hong2023central}.
\end{remark}

Note that \(\mu^x(dy)\) is the unique invariant measure of an ergodic Markov process \(Y_t^x\) with the generator \(\mathcal{L}^x\) with given \(x\in \mathbb{R}^n\). It is direct to see that there exists a solution to the Poisson equation, and more properties could be found in \cite{Pardoux2001,Pardoux2003,Pardoux2005}.
\begin{theorem}\label{Poisson}
Under assumptions \textbf{H1}-\textbf{H3}. Then for any function \(h\in C^{2,2+\gamma}\) with \(\gamma\in(0,1]\) which is ``centering" that is \(\bar h = \int_{\mathbb{R}^{m}} h(x,y)\mu^x(dy)=0\), there exists a unique solution $u$ satisfying the following Poisson equation in \(\mathbb{R}^m\) which depends on the parameter \(x\), 
\begin{equation}\label{Poisson_eq}
-\mathcal{L}^x(y)u(x,y) = h(x,y).
\end{equation}
Here, \(\mathcal{L}^x\) is the infinitesimal generator of the fast process \(Y\) where \(x\) is considered as a parameter. 
Moreover, \(u(x,y)\in C^2\). 
\end{theorem}

Then we give some regularity properties of solution to the Poisson PDE \eqref{Poisson_eq}.
\begin{proposition}\label{Poisson_g}
Under assumptions \textbf{H1}-\textbf{H3} and \(\bar h (x)= \int_{\mathbb{R}^{m}} h(x,y)\mu(\dd y)=0\). Define that 
\begin{equation}\label{Phi}
u(x,y) =\int_{0}^{\infty}\mathcal{P}_t^x h(x, Y_t^x) \dd t.
\end{equation} 
Then \(u(x,y)\) is the solution to the Poisson equation \eqref{Poisson_eq}.
Moreover, \(u(x,y)\in C^2\). For all \((x,y)\in \mathbb{R}^n\times \mathbb{R}^m\), \(D_xu(x,y)\) and \(D^2_{xx}u(x,y)\) is linear growth, that is for all \((x,y)\in\mathbb{R}^m\times \mathbb{R}^n\), there exists a constant $L_u>0$ such that 
$$|D_xu(x,y)|\vee|D^2_{xx}u(x,y)|\le L_u(1+|x|+|y|).$$
And \(D_yu(x,y), D^2_{xy}u(x,y), D^2_{yy}u(x,y)\) is globally bounded, that is there exists a constant $L_u>0$ such that 
$$|D_yu(x,y)|\vee|D^2_{xy}u(x,y)|\vee|D^2_{yy}u(x,y)|\le L_u.$$
Moreover, when \(Y\) is independent of \(x\), \(\sup_{x\in \mathbb{R}^m}|u(x,y)|\le C(1+|y|)\).
\end{proposition}
\para{Proof}. Its proof could be found at \cite[Proposition 3.1]{hong2023central}, so details of the proof is omitted but we give more explanation for the case that there is no feedback from slow one into $y$. Recall that the transition semigroup \(\mathcal{P}_t\) related to \eqref{Xi} admits a unique invariant probability measure \(\mu\) on $\mathbb{R}^m$ which is irrelevant of $x\in \mathbb{R}^n$.
So for all \(y\in \mathbb{R}^{m}\), any bounded measurable function \(u(y): \mathbb{R}^m \mapsto \mathbb{R}\),
$$\|\mathcal{P}_t u(y)-\int_{\mathbb{R}^{m}} u(y)\mu(dy)\|_{\infty}\le C e^{-ct}\|u\|_{Lip}(1+|Y_0|)$$
and 
$$\|\mathcal{P}_t u(y)\|_{Lip}\le C e^{-ct}\|u\|_{Lip}(1+|Y_0|).$$ 
The proof is completed. \qed

We next give a remark regarding the regularity of $\bar h(x)$.
\begin{remark}\label{bar_h}
Let $\eta=1,2$. Assume \textbf{H1} and \(h \in BC^2\). Then,
\begin{equation}\label{partial_eta_f}
D_x^\eta\bar h (x)=\int_{\mathbb{R}^n}\big[D_x^\eta h(x,y)-\sum_{\iota=1}^{\eta}C_\eta^\iota\frac{D^\iota\mathcal{L}^x}{D x^\iota}(x,y)\cdot D_x^{\eta-\iota} u(x,y)\big] \mu^x(\dd y),
\end{equation}
where \(u\) is defined by \eqref{Phi} and \(C_\eta^\iota>0\) are different constants only depending on \(\eta,\iota\).
\end{remark}
\para{Proof}. First we focus on the situation of \(\eta=1\). By the chain rule, we have
\begin{equation}\label{partial_1_barf}
D_x\bar h(x) =\int_{\mathbb{R}^m} D_x h(x,y) \mu^x(\dd y)+\int_{\mathbb{R}^m} h(x,y) D_x p_\infty(x,y)\dd y.
\end{equation}
Note that 
\begin{equation}\label{partial_p}
D_xp_\infty (x,y) =-\int_{0}^{\infty}ds \int_{\mathbb{R}^m}dy' p_\infty(y',x) D_x\mathcal{L}^x(x,y')p_s(y',y;,x)\dd y',
\end{equation}
where 
\begin{equation}\label{D_L}
D_x\mathcal{L}^x(x,y)=\frac{1}{2} D_x a_{i,j}(x,y)\frac{\partial^2}{\partial y_i \partial y_j}+D_xb_{i}(x,y)\frac{\partial}{\partial y_i}.
\end{equation}
Therefore, by leveraging the \eqref{partial_p}, Fubini's theorem and property of transition density function, we have
\begin{equation}\label{partial_barf}
\begin{aligned}
&\int_{\mathbb{R}^m} h(x,y) D_x p_\infty(x,y)dy \\
&=- \int_{\mathbb{R}^m} dy'D_x\mathcal{L}^x(x,y')\big( \int_{0}^{\infty}\dd s \int_{\mathbb{R}^m}\dd y p_s(y',y;,x)h(x,y)\big)p_\infty(y',x)\\
&= - \int_{\mathbb{R}^m} dy'D_x\mathcal{L}^x(x,y')\big( \int_{0}^{\infty}\dd s \mathcal{P}_t^x h(x,y') \big)p_\infty(y',x)\\
&= - \int_{\mathbb{R}^m} D_x\mathcal{L}^x(x,y')u(x,y)\mu^x(\dd y').
\end{aligned}
\end{equation} 
So \eqref{partial_eta_f} holds for \(\eta=1\). Next, we turn to the case of \(\eta=2\).
\begin{equation}\label{partial_2_barf}
\begin{aligned}
&D^2_{xx}\bar h(x) \\
&=\int_{\mathbb{R}^m} D^2_{xx} h(x,y) \mu^x(dy)- \int_{\mathbb{R}^m} D^2_{xx}\mathcal{L}^x(x,y)u(x,y)\mu^x(\dd y)\\
&\quad-\int_{\mathbb{R}^m} D_{x}\mathcal{L}^x(x,y')D_x\Psi(x,y)\mu^x(\dd y)+\int_{\mathbb{R}^m} D_{x} h(x,y) D_xp_\infty(x,y)dy\\
&\quad- \int_{\mathbb{R}^m} D_{x}\mathcal{L}^x(x,y)u(x,y)D_xp_\infty(x,y)\dd y.
\end{aligned}
\end{equation}
Set
\[
H_1:=\int_{\mathbb{R}^m} D_x h(x,y) D_x p_\infty(x,y)\,dy
\quad\text{and}\quad
H_2:=-\int_{\mathbb{R}^m} D_x \mathcal{L}^x(x,y)\Psi(x,y) D_x p_\infty(x,y)\,dy.
\]
Then, analogously to the second equality in \eqref{partial_barf}, we have
\begin{eqnarray}\label{H1}
H_1&=& - \int_{\mathbb{R}^m} dy'D_x\mathcal{L}^x(x,y')\big( \int_{0}^{\infty}\dd s \int_{\mathbb{R}^m}\dd y p_s(y',y;,x)D_x\hat h(x,y)\big)\cr
&&p_\infty(y',x).
\end{eqnarray}
And 
\begin{eqnarray}\label{H2}
H_2&=& - \int_{\mathbb{R}^m} dy'D_x\mathcal{L}^x(x,y')\big( \int_{0}^{\infty}\dd s \int_{\mathbb{R}^m}\dd y p_s(y',y;,x) D_x\mathcal{L}^x(x,y)u(x,y)\big)\cr
&&p_\infty(y',x).
\end{eqnarray}
By combining \eqref{H1} and \eqref{H2}, we arrive at
\begin{eqnarray}\label{H}
H_1+H_2
&= &- \int_{\mathbb{R}^m} dy'D_x\mathcal{L}^x(x,y')\big( \int_{0}^{\infty}ds \int_{\mathbb{R}^m}\dd y D_x(p_s(y',y;,x) h(x,y))\big)p_\infty(y',x)\cr
&&+ \int_{\mathbb{R}^m} dy'D_x\mathcal{L}^x(x,y')\big( \int_{0}^{\infty}ds \int_{\mathbb{R}^m}\dd y \big[D_xp_s(y',y;,x) h(x,y)\big.\big.\cr
&&\big.\big.+p_s(y',y;,x)D_x\mathcal{L}^x(x,y)u(x,y)\big]\big)p_\infty(y',x)\cr
&= &- \int_{\mathbb{R}^m} dy'D_x\mathcal{L}^x(x,y')D_x\big( \int_{0}^{\infty}ds \int_{\mathbb{R}^m}\dd y p_s(y',y;,x) h(x,y)\big)p_\infty(y',x)\cr
&= & - \int_{\mathbb{R}^m} D_x\mathcal{L}^x(x,y')D_xu(x,y')\mu^x(\dd y')
\end{eqnarray}
Observe that the second equality is justified by the vanishing of the second term in the first line of \eqref{H}, which follows from the estimates in \cite[(34)]{Pardoux2003}. Consequently, it is not difficult to verify that \eqref{partial_eta_f} also holds for $\eta=2$. The proof is therefore complete.

\section{Central Limit Theorem: Case 1} 

In this section we consider the following slow-fast system
\begin{equation}\label{simiplified1_System}
\dd x_t^\epsilon = g(x_t^\epsilon,y^\epsilon_t) \dd t+\sqrt \epsilon f(x_t^\epsilon,y^\epsilon_t) \dd B^H_t ,\\ x^\epsilon_0 = x.
\end{equation}
here, \(f\) is assumed of \(BC^2\), \(g\) is both Lipschitz continuous and globally bounded. 
While the effective average dynamics, obtained as $\epsilon\to 0$, is quite simple and deterministic, it is given by: \\
\begin{equation}\label{averaged1-system}
\dd \bar x_t = \bar g(\bar x_t) \dd t,\\ \bar x_0 = x.
\end{equation}
where \(\bar g = \int_{\mathbb{R}^{m}} g(x,y)\mu^x(\dd y)=0\) with stationary measure of \(y\) which is denoted by \(\mu^x\).

Setting
$$z_t^\epsilon =\f 1{\sqrt \epsilon } (x_t^\epsilon -\bar x_t),$$
then
$$\dd z_t^\epsilon =
\f 1{\sqrt \epsilon }(g(x_t^\epsilon,y^\epsilon_t)-\bar g(\bar x_t))\dd t+f(x_t^\epsilon,y^\epsilon_t)\dd B_t^H.$$

\begin{proposition}\label{Poisson_g0}
Under assumptions \textbf{A1}-\textbf{A3} and \(\bar g (x)= \int_{\mathbb{R}^{n}} g(x,y)\mu(dy)=0\). Define that 
\begin{equation}\label{phi}
\Psi(x,y) =\int_{0}^{\infty}\mathcal{P}_t g(x, Y_t) \dd t.
\end{equation} 
The \(\Psi(x,y)\) is the solution to the Poisson equation \eqref{Poisson_eq} where \(h\) is replaced by \(g\).
Moreover, \(\Psi(x,y)\in C^2\). For all \((x,y)\in \mathbb{R}^n\times \mathbb{R}^m\), \(D_x\Psi(x,y)\) and \(D^2_{xx}\Psi(x,y)\) is linear growth. \(D_y\Psi(x,y), D^2_{xy}\Psi(x,y), D^2_{yy}\Psi(x,y)\) is globally bounded. 
Furthermore, since \(Y\) is independent of \(x\), \(\Psi\in C^2\) is Lipshitz continuous with respect to \(y\) uniformly over all \(x\) .
\end{proposition}
\para{Proof}. Since $\mathcal{L}$ is the infinitesimal generator of the Markov process $\{y_t\}_{t\in \mathbb{R}^{+}}$, and under assumptions \textbf{A1}--\textbf{A3}, it is straightforward to verify that $\Psi$ satisfies the Poisson equation \eqref{Poisson_eq}. Moreover, by \cite[Proposition 4.1]{rockner2021strong}, $D_x\Psi(x,y)$ and $D^2_{xx}\Psi(x,y)$ are of linear growth.

\begin{remark}
If \(f(x,y)=f(x)\), the assumption \textbf{A1} could be weaken to the assumption that \(b, \sigma \in C^{1,2+\gamma}\) for \(\gamma\in(0,1]\).
\end{remark} 

Moreover, a further assumption is required in order to bound the family $\{z^\epsilon\}_{\epsilon\in(0,1]}$.
\begin{itemize}
\item[A4.] Assume that for \(p\ge 1\) and \(\theta\in (0,\frac{1}{2})\), there exists a constant \(C>0\) such that for all \(0<\epsilon\le 1\),
$$ \mathbb{E}[\|y^\epsilon\|^p_{\infty}] \le C,\quad \mathbb{E}[\|y^\epsilon\|^p_{\theta-hld}] \le C \epsilon^{-\frac{p}{2}}.$$ 
\end{itemize}
Assumption \textbf{A4} may appear somewhat stringent. In fact, we provide an illustrative example in Appendix B to clarify this point.

Denote 
\begin{equation}\label{Dyphisigma}
D_y\Psi \sigma(x,y) u:=D_y\Phi (x,y)\cdot\sigma(x,y)u, \quad u \in \mathbb{R}^e,
\end{equation}
then, directly we have $D_y\Psi \sigma(x,y)\in \mathbb{R}^m\otimes \mathbb{R}^{e}$ and \\
\begin{equation}\label{Dyphisigma2}
(D_y \Psi \sigma)(D_y \Psi \sigma)^{*}(x,y):=\int_{\mathbb{R}^{m}} (D_y\Psi \sigma(x,y))(D_y\Psi \sigma(x,y))^{*}\mu^{x}(dy).
\end{equation}
\begin{theorem}\label{clt1}
Under \textbf{A1}--\textbf{A4}. As \(\epsilon\to 0\), \(z^\epsilon\) weakly converges to \(\bar z\), that is 
\begin{equation}\label{clt-limit0}
\dd \bar z_t = D\bar g(\bar x_t) \bar z_t \dd t+\bar f(\bar x_t) \dd B^H_t + \bar V^{\frac{1}{2}}(\bar x_t)\dd \hat W_t,
\end{equation}
where \(V(x,y):=(D_y \Psi \sigma)(D_y \Psi \sigma)^{*}(x,y)\) with $\bar V^{\frac{1}{2}}(x)=\int_{\mathbb{R}^{m}} V^{\frac{1}{2}}(x,y)\mu^{x}(dy)$.
\end{theorem}

We now outline our strategy for proving the CLT result in Theorem \ref{clt1}. To this end, we first recall the definition of a weak solution to the Young--It\^o SDE driven by mixed FBM.
\begin{definition}\label{weaksolution}
A weak solution to the following Young-It\^o SDE
\begin{equation}\label{weak_solution}
\dd X_t = F( X_t) \dd B^H_t + G( X_t) \dd t+ V( X_t)\dd W_t,\ X_0 = x \in \mathbb{R}^n,
\end{equation}
is a triple \((X,B^H, W)\) on the probability space \((\Omega, \mathcal{F}, \mathbb{P})\) with filtration \(\{\mathcal{F}_t, t\in[0,T]\}\), where
\begin{itemize}
\item[1.] \((\Omega, \mathcal{F}, \mathbb{P})\) is a complete probability space, \(\{\mathcal{F}_t, t\in[0,T]\}\) is a right-continuous filtration such that \(\mathcal{F}_0\) contains the \(\mathbb{P}\)-null sets;
\item[2.] \(W\) is a \(\mathcal{F}_t\)-\(r\)-dimensional BM;
\item[3.] \(B^H\) is a FBM of Hurst parameter \(H\) that is \(\mathcal{F}_0\)-measurable;
\item[4.] The process \(X\) is \(\mathcal{F}_t\)-adapted and has trajectories in \(C^{\alpha}\) almost surely;
\item[5.] \((X,B^H, W)\) satisfies Young-It\^o SDE \eqref{weak_solution}.
\end{itemize}
\end{definition}

\begin{theorem}\label{Solution_Limit}
Assume \(\bar f\in BC^1\) and \(g\) is Lipschitz continuous and globally bounded. Then, there exists a unique solution to the Young-It\^o SDE \eqref{clt-limit0}. 
\end{theorem}
\para{Proof}. Note that \(\Psi(\cdot,y)\in BC^1\) for every \(y\) and \(\sigma\) is assumed globally bounded and Lipschitz continuous, so we could show that \(V\) is globally bounded and Lipschitz continuous. 
The existence of the weak solution and uniqueness in law implies the existence of strong solution, which comes from the Yamada-Watanabe theorem. So the existence and uniqueness of the solution could be shown at \cite[Theorem 2.2]{Nualart03092008}.

To show Theorem \ref{clt1}, the strategy in our proof is as below.
\begin{itemize}
\item We first aim to show that the family $\{z^\epsilon\}_{\epsilon \in (0,1]}$ is relatively compact as a set of stochastic processes. By virtue of the embedding between different H\"older continuous spaces, we obtain the tightness of the laws of $z^\epsilon$ in the H\"older topology.
\item Then, by Prokhorov's theorem, every sequence of the laws of $z^\epsilon$ admits a subsequence that converges weakly to some limiting probability measure.
\item To apply Skorokhod's theorem, we introduce a new probability space (keeping the same notation for simplicity). Then there exists a subsequence of $\{z^\epsilon\}_{0<\epsilon\le 1}$ that converges almost surely to some limit in the H\"older topology. (Here, the little H\"older space is needed because it is separable and complete, whereas the usual H\"older space is not separable.)
\item Finally, we show that any such limit satisfies \eqref{Solution_Limit} in the uniform topology. By pathwise uniqueness in the uniform topology, the convergence in $C([0,T], \mathbb{R}^n)$ follows.
\end{itemize}

\subsection{A priori estimates}
In this subsection, we give some a priori estimates.

\begin{proposition}\label{boundness_x}
Let \((\mathcal{V}, \mathcal{W})\) be Banach spaces and suppose \(K \in L_p(C^\alpha([0, T], \mathcal{V}))\), \(B \in C^\beta([0, T], U)\) for some 
\((\beta \in (\frac{1}{2}, 1])\), \((0 < \alpha \le \gamma < 1)\) and \((\alpha + \beta > 1)\). Also suppose that $f$ is of $BC^2$ and $G$ are Lipschitz continuous with respect to \(x\) and globally bounded over \(r\in[0,t]\).
Then, if \(X : [0, T] \to \mathcal{V}\)
satisfies the controlled Young differential equation
\begin{equation}\label{eq:gronwall-YDE}
X_t = X_0 + \int_0^t f(s,X_s) \dd B_s + \int_0^t G(s) X_s\dd s + K_t.
\end{equation}
Then, 
\begin{equation}\label{infty-norm}
\|X\|_{\infty,p,[0,T]} \le C_2 e^{(2|f|^{p}_{-\kappa,p})^{\frac{1}{(\bar \eta-\alpha)p }}}(\|X_0\|^p+\mathbb{E}[\|K\|_\alpha^p]+\|f\|_{\infty}+1)
\end{equation}
and 
\begin{equation}\label{alpha-norm}
\|X\|_{\alpha,p} \le C_1(\mathbb{E}[\|K\|_\alpha^p]+\mathbb{E}[\|X\|^p_{\infty}]+1).
\end{equation}
for some constants \(C_1, C_2 > 0\) depending on \(\alpha, \beta, H, T\) and \(\|f\|_{BC^2}, \|G\|_{Lip}\). 
\end{proposition}
\para{Proof}. First, by leveraging the stochastic sewing lemma, \cite[Lemma 3.10]{Hairer-Li}, that is if $f\in \mathcal{C}^{-\kappa,\gamma}$ for some $\kappa,\gamma>0$ such that $\eta>1/2$ and $\bar \eta>1$. 
Then it has $A_{u,v}:=\int_{u}^{v}f(r, x_u)d B_{r}^H$ belongs to $H_\eta^p \cap \bar H_{\bar \eta}^p$ with $\eta>1/2$ and $\bar \eta>1$. 
Set \(I_{u,v}:=\int_{u}^{v}f(r, x_r)d B_{r}^H\), then it has $\|I_{u,v}(A)\|_p \lesssim |u-v|^\eta$ and $\|\mathbb{E}(I_{u,v}(A)-A_{u,v}|\mathcal{F}_u)\|_p \lesssim |u-v|^{\bar\eta}$. Then, we have 
\begin{equation}\label{f(r,x_r)}
\mathbb{E}[\|\int_{s}^{t}f(r, x_r)d B_{r}^H\|_{p}] \le |f|^p_{-\kappa,p}(|t-s|^{H-\kappa}+\|x\|_{\alpha}^p(t-s)^{(\bar \eta -\alpha)p}).
\end{equation}
By combining \eqref{f(r,x_r)} and some straight computation, for any $0\le s<t\le T$ such that $|t-s|<1$,
\begin{eqnarray}\label{X_p}
\mathbb{E}[\|X_{s,t}\|^p]&\le& 3^{p-1} \mathbb{E}[\|K_{s,t}\|]+3^{p-1}\mathbb{E}[\|\int_s^t f(r, X_r) \dd B_r\|^p]+3^{p-1}\mathbb{E}[\|\int_s^t G(r, X_r) \dd r\|^p]\cr
&\le&3^{p-1}\mathbb{E}[\|K\|^p_\alpha|t-s|^{\alpha p}]+3^{p-1}\mathbb{E}[|f|^p_{-\kappa,p}(|t-s|^{H-\kappa}+\|x\|_{\alpha}^p(t-s)^{(\bar \eta -\alpha)p})]\cr
&&+3^{p-1}\mathbb{E}[(1+\|x\|_{\infty})(t-s)^p].
\end{eqnarray}
Note that \(f\) is assumed globally bounded, so there exists a constant \(C_1\) such that \(|f|^p_{-\kappa,p}\le C_1\). 
Then, by choosing $s,t$ small enough such that $|f|^p_{-\kappa,p}(t-s)^{(\bar \eta -\alpha)p}<1/2$, we have that 
\begin{equation}\label{X_alpha_p}
\begin{aligned}
\mathbb{E}[\|X\|_{\alpha,[s,t]}^p] &\le 2\times 3^{p-1}\mathbb{E}[\|K\|_\alpha^p]+2\times 3^{p-1}|f|^p_{-\kappa,p}|t-s|^{(H-\kappa-\alpha)p}\cr
&\quad+2\times 3^{p-1}|t-s|^{(1-\alpha)p}+2\times 3^{p-1}\mathbb{E}[\|X\|^p_{\infty,[s,t]}]|t-s|^{(1-\alpha)p}.
\end{aligned}
\end{equation}
By leveraging the sub-additivity of the H\"older norm in all time interval, the estimate \eqref{alpha-norm} has arrived.

On the other hand, we have 
\begin{eqnarray}\label{X_infty}
\mathbb{E}[\|X\|_{\infty,[s,t]}^p] &\le& \mathbb{E}[\|X_s\|_{\infty}^p] +\mathbb{E}[\|X\|_{\alpha,[s,t]}^p]|t-s|^{\alpha p}\cr
&\le& \mathbb{E}[\|X_s\|_{\infty}^p] +2\times 3^{p-1}\mathbb{E}[\|K\|_\alpha^p]|t-s|^{\alpha p}+2\times 3^{p-1}|f|^p_{-\kappa,p}|t-s|^{(H-\kappa)p}\cr
&&+2\times 3^{p-1}|t-s|^{(1-\alpha)p}+2\times 3^{p-1}\mathbb{E}[\|X\|^p_{\infty,[s,t]}]|t-s|^{p}.
\end{eqnarray}
Then, by taking $s,t$ sufficiently small such that $(3|t-s|)^p<1/4$, so we have
\begin{eqnarray}\label{X_infty_2}
\mathbb{E}[\|X\|_{\infty,[s,t]}^p] & \le& 2\times 3^{p-1}\mathbb{E}[\|X_s\|_{\infty}^p] \cr
&&+4\times 3^{p-1}\big[\mathbb{E}[\|K\|_\alpha^p]+|f|^p_{-\kappa,p}|t-s|^{(H-\kappa)p}+|t-s|^{(1-\alpha)p}\big].
\end{eqnarray}
Note that \(s,t\) are assumed \(t-s\le 1\). Take the uniform partition $\mathcal{P}:=\{0=t_0<t_1<\cdots<t_n=T\}$ on the time interval $[0,T]$ with mesh size $|t_{i+1}-t_i|=(2\times 3^{p-1}|f|^{p}_{-\kappa,p})^{-\frac{1}{(\bar \eta-\alpha)p }}\wedge 3\times 4^{-1/p}.$ 
So $|\mathcal{P}|=(2\times 3^{p-1}|f|^{p}_{-\kappa,p})^{\frac{1}{(\bar \eta-\alpha)p }}$. Consequently, by iterating the estimates in all subintervals, we have \\
\begin{eqnarray}\label{X_infty_1}
\mathbb{E}[\|X\|_{\infty}^p] &\le& (2\times 3^{p-1})^{|\mathcal{P}|}[\|X_0\|^p] \cr
&&+4\times 3^{p-1}\big[\mathbb{E}[\|K\|_\alpha^p]+|f|^p_{-\kappa,p}+1\big](1+2+2^2+\cdots+2^{|\mathcal{P}|}).
\end{eqnarray}
So the \eqref{infty-norm} is arrived. \qed

Let $x\in B_{\alpha,p}$ with $p\ge 2$ and $\alpha\in(0,1/2)$ be such that $1+\frac{1}{p}<\alpha+H$. For a given $\mathcal{F}_t$-adapted process $x$ (not necessarily the solution to \eqref{simiplified_System}), we introduce a fast-varying process $\tilde y^{\epsilon}$ as follows:
\begin{equation}\label{fast2}
d\tilde y^{\epsilon}_t = \frac{1}{\epsilon} b(x_t, \tilde y^{\epsilon}_t)\,dt + \frac{1}{\sqrt{\epsilon}}\sigma(x_t, \tilde y^{\epsilon}_t)\,dw_t, \qquad \tilde y^\epsilon_0=y.
\end{equation}
Then the following estimates hold whose proof is a small extension of \cite[Theorem 4.16]{Hairer-Li}, and hence the details are omitted for brevity.
\begin{theorem}\label{f_x_y}
Let \(p\ge 2\) and \(\alpha\in(0,1/2)\) be such that \(1+\frac{1}{p}<\alpha+H\). There are exponents \(\eta>\frac{1}{2}\) and \(\bar \eta>1\) such that for \(h:\mathbb{R}^n\times \mathbb{R}^m\mapsto \mathbb{R}\) a bounded uniformly Lipschitz continuous function, one has the following:
\begin{itemize}
\item[1]. Suppose that \(\bar h(x)=\int_{\mathbb{R}^m}h(x,y)\mu(\dd y)=0\) for all \(x\in \mathbb{R}^n\). Then for any \(\beta<H\), there exists a constant \(\kappa>0\) such that 
\begin{equation}
\|\int_{s}^{t}h(x_r,y_r^\epsilon)\dd B_{r}^H\|_{p} \le C \epsilon^\kappa |h|_{BC^1}(|t-s|^\eta+\|x\|_{\alpha}^p(t-s)^{\bar \eta }).
\end{equation}
\item[2]. For general \(h\), one has the bound 
\begin{equation}
\|\int_{s}^{t}h(x_r,y_r^\epsilon)\dd B_{r}^H\|_{p} \le C |h|_{BC^1}(|t-s|^\eta+\|x\|_{\alpha}^p(t-s)^{\bar \eta }).
\end{equation}
where $C>0$ is a constant.
\end{itemize}
\end{theorem}

\subsection{Tightness criterion}
Note that \( \Psi(x,y) =\int_{0}^{\infty}\mathcal{P}_t g(x, Y_t) \dd t\). 
According to the Proposition \ref{Poisson_g}, it is the solution to the following Poisson equation 
$$ -\mathcal{L}^x(y)\Psi(x,y) = g(x,y)-\bar g(x)$$ where \(\mathcal{L}\) is the infinitesimal generator of the fast process \(y^\epsilon\) with fixed slow one. Note that \(\Psi\in C^2\) as a function of \((x,y)\).

By leveraging the Young-It\^o formula to the function \(\Psi\), then
\begin{eqnarray}\label{Phi}
&&\Psi(x^\epsilon_t,y_t^\epsilon)-\Psi(x_0,y_0)\cr
&=& \Psi(x_0,y_0)-\Psi(x^\epsilon_t,y_t^\epsilon)
+D_x\Psi(x^\epsilon_r,y_r^\epsilon)g(x^\epsilon_r, y_r^\epsilon)\dd r\cr
&&\big.+\sqrt{\epsilon}\int_{0}^{t}D_x\Psi(x^\epsilon_r,y_r^\epsilon)f(x^\epsilon_r,y_r^\epsilon)\dd B_r^H+\frac{1}{{\epsilon}}\int_{0}^{t}\mathcal{L}^x\Psi(x^\epsilon_r,y_r^\epsilon)\dd r
+\frac{1}{\sqrt{\epsilon}}M_t^\epsilon
\end{eqnarray}
where $M_t^\epsilon$ is a local martingale related to $y^\epsilon$ and its coefficient is denoted by $V^{\frac{1}{2}}(x,y)$ where $V(x,y):=(D_y \Psi \sigma)(D_y \Psi \sigma)^{*}(x,y)$.

Then, $z^\epsilon$ could be rewritten as below
\begin{eqnarray}\label{homogenization-Poisson}
z_t^\epsilon -z_0
&= &\int_0^t f(x_r^\epsilon,y_r^\epsilon) \dd B^H_r + \sqrt{\epsilon}\big[ \Psi(x_0,y_0)-\Psi(x^\epsilon_t,y_t^\epsilon)
+\int_{0}^{t}D_x\Psi(x^\epsilon_r,y_r^\epsilon)g(x^\epsilon_r, y_r^\epsilon)\dd r\big.\cr
&&\big.+\sqrt{\epsilon}\int_{0}^{t}D_x\Psi(x^\epsilon_r,y_r^\epsilon)f(x^\epsilon_r,y_r^\epsilon)\dd B_r^H
\big]+\int_{0}^{t}\big(\int_{0}^{1}D\bar g(\theta x^\epsilon_r+(1-\theta)\bar x_r)\dd \theta\big)z^\epsilon_r \dd r
\cr
&&+M_t^\epsilon\cr
&=& \int_0^t f(x_r^\epsilon,y_r^\epsilon) \dd B^H_r +\sqrt{\epsilon}\big[ \Psi(x_0,y_0)-\Psi(x^\epsilon_t,y_t^\epsilon)+ \int_{0}^{t}D_x\Psi(x^\epsilon_r,y_r^\epsilon)g(x^\epsilon_r, y_r^\epsilon)\dd r\big.\cr
&&\big.+\sqrt{\epsilon}\int_{0}^{t}D_x\Psi(x^\epsilon_r,y_r^\epsilon)f(x^\epsilon_r,y_r^\epsilon)\dd B_r^H
\big]+\int_{0}^{t}\big(\int_{0}^{1}D\bar g(\theta x^\epsilon_r+(1-\theta)\bar x_r)\dd \theta\big)z^\epsilon_r \dd r
\cr
&& +M_t^\epsilon.
\end{eqnarray}
We next show that $z^\epsilon$ is bounded under the H\"older norm. By virtue of the compact embedding between different finite dimensional H\"older continuous spaces, the tightness of the laws of $z^\epsilon$ in the H\"older topology follows. Before showing   this boundedness result, we first establish some necessary a priori estimates.

\begin{lemma}\label{clt1_x_bounded}
Let \(p\ge 2\) and \(\alpha\in(0,H)\). Under \textbf{H1}--\textbf{H3}. Then for all \(\epsilon\in(0,1]\),
\begin{eqnarray}\label{clt_x_bdd}
\mathbb{E}[\|x^\epsilon\|^p_{\alpha}] \le C_{x_0,T}, \quad \|\bar x\|_{\alpha}\le C_{x_0,T}
\end{eqnarray} 
where the constant \(C_{x_0,T}>0\) only depends on \(x_0\) and \(T\).
\end{lemma}
\para{Proof}.  Along the approach in \cite[Corollary 4.20]{Hairer-Li}, the above result is easy to obtain, so details of the proof are omitted. \qed

The following result can be directly obtained by extending \cite{Hairer-Li}. so proof details are omitted. (Due to the existence of small FBM not general one, so a general Gronwall lemma is enough and \cite[Lemma 2.2]{Hairer-Li} is not needed here).
\begin{lemma}\label{x1_convergence}
Let \(p\ge 2\) and \(\gamma\in(\frac{1}{2},H)\). Under \textbf{H1}--\textbf{H3}. As \(\epsilon\) tends to zero, there exists a \(\kappa>0\) such that
\begin{eqnarray}\label{clt_x_bdd0}
\|x^\epsilon-\bar x\|_{\gamma,p} \le C_{x_0,T} \epsilon^\kappa 
\end{eqnarray} 
where \(C_{x_0,T}>0\) only depends on \(x_0\) and \(T\) and \(\bar x\) is the solution to the deterministic ODE \eqref{averaged1-system}.
Moreover, \(\bar g\) is globally Lipschitz continuous.
\end{lemma}

\begin{proposition}\label{Poisson_CLT}
Under assumptions \textbf{H1}-\textbf{H3}. Define \(\hat g(x,y):=g(x,y)-\bar g(x)\) and 
\begin{equation}\label{psi}
\Psi(x,y) =\int_{0}^{\infty}\mathcal{P}_t^x \hat g(x, Y_t^x) \dd t.
\end{equation} 
The \(\Psi(x,y)\) is the solution to the Poisson equation \eqref{Poisson_eq} where \(h\) is replaced by \(\hat g\).
Moreover, \(\Psi(x,y)\in C^2\). For all \((x,y)\in \mathbb{R}^n\times \mathbb{R}^m\), \(D_x\Psi(x,y)\) and \(D^2_{xx}\Psi(x,y)\) are linear growth. Specially,
\(D_y\Psi(x,y), D^2_{xy}\Psi(x,y), D^2_{yy}\Psi(x,y)\) are globally bounded. 
\end{proposition}

\begin{remark}\label{Poisson_g0}
If Assumptions \textbf{H1}--\textbf{H3} hold, according to Proposition \ref{Poisson_CLT},
for every \(\hat g(x,y):=g(x,y)-\bar g(x)\), there exists a solution \(\Psi(x,\cdot)\) satisfying that \(\Psi(x,\cdot) \in BC^3\) for every \(x\in \mathbb{R}^n\) and \(u(\cdot,y)\in BC^1\) for every \(y\in \mathbb{R}^m\), to the Poisson equation \eqref{psi}.
Then, according to \cite[Lemma 3.2]{rockner2021strong}, it deduces that \(\bar g\in BC^2\). 
When there is no feedback from the slow dynamics into the fast one, it ensures that \(\bar g\in BC^2\) by only assuming \(D_{x} g( x,y)\) is Lipschitz with respect to \(x\in\mathbb{R}^n\) for each \(y\in\mathbb{R}^m\).
\end{remark}

\begin{lemma}\label{partial_x_phi}
Let \(p\ge 2\) and \(\alpha\in(0,H)\). There are exponents \(\eta>\frac{1}{2}\) such that for \(f:\mathbb{R}^n\times \mathbb{R}^m\mapsto \mathbb{R}\) a bounded uniformly Lipschitz continuous function, we have
$$ \|\int_{s}^{t}\sqrt{\epsilon}D_x\Phi(x_r^\epsilon,y_r^\epsilon)f(x_r^\epsilon,y_r^\epsilon)\dd B_{r}^H\|_{p} \le C_{x_0,T}|t-s|^\eta
$$
where $ C_{x_0,T}>0$ is a constant only depending on $x_0$ and $T$.
\end{lemma}
\para{Proof}. 
Set \(I_{s,t}(A^\epsilon):=\int_{s}^{t}\sqrt{\epsilon}D_x\Psi(x_r^\epsilon,y_r^\epsilon)f(x_r^\epsilon,y_r^\epsilon)d B_{r}^H\) and \(A^\epsilon_{s,t}:= \sqrt{\epsilon}\int_{s}^{t}D_x\Psi(x_s^\epsilon,y_r^\epsilon)f(x_s^\epsilon,y_r^\epsilon)d B_{r}^H\). 
The integral \(I_{s,t}(A^\epsilon)\) is constructed via the stochastic sewing lemma \cite[Section 3.3]{Hairer-Li}. In fact, since we assume the condition \textbf{A4}, it has for all \(p\ge 2\) and \(\theta\in (0,\frac{1}{2})\), 
$$ \mathbb{E}[\|y^\epsilon\|^p_{\infty}] \le C,\quad \mathbb{E}[\|y^\epsilon\|^p_{\theta-hld}] \le C \epsilon^{-\frac{p}{2}}.$$ 
Set \(\tilde A^\epsilon_{s,t}:= \sqrt{\epsilon}D_x\Psi(x_s^\epsilon,y_s^\epsilon)f(x_s^\epsilon,y_s^\epsilon)( B_{t}^H- B_{s}^H)\). Recall that from Proposition \ref{Poisson_g}, $D_x\Phi$ is Lipschitz continuous with respect to all \((x,y)\in\mathbb{R}^m\times \mathbb{R}^n\), so it is straightforward to see that, 
\begin{eqnarray}\label{AtildeA}
&&\mathbb{E}[|A^\epsilon_{s,t}-\tilde A^\epsilon_{s,t}|^p]\cr
&\le& \mathbb{E}[|\int_{s}^{t} \sqrt{\epsilon}(D_x\Psi(x_s^\epsilon,y_r^\epsilon)f(x_s^\epsilon,y_r^\epsilon)-D_x\Psi(x_s^\epsilon,y_s^\epsilon)f(x_s^\epsilon,y_s^\epsilon))d B_{r}^H|^p]\cr
&\le&\|y^\epsilon\|^p_{\theta,p}\|B^H\|_\alpha^p|t-s|^{(\alpha+\beta)p}
\end{eqnarray}
Since $\theta<\frac{1}{2}$ and $\beta\in(\frac{1}{2},H)$ such that $\theta+\beta>1$,
so \(I_{s,t}(A^\epsilon)\) equals the Young integral for all \(\epsilon\in(0,1]\). 

With some direct computation,
\begin{equation}\label{D_xpsi_f}
\begin{aligned}
\mathbb{E}\big[\big|\sqrt{\epsilon}\int_{s}^{t}D_x\Psi(x_r^\epsilon,y_r^\epsilon)f(x_r^\epsilon)\dd B_{r}^H\big|^p\big]
&\le C\epsilon^{\frac{p}{2}}\mathbb{E}\big[\big|\int_{s}^{t}(D_x\Psi(x_r^\epsilon,y_r^\epsilon)-D_x\Psi(x_s^\epsilon,y_s^\epsilon))f(x_r^\epsilon,y_r^\epsilon)\dd B_{r}^H\big|^p\big]\cr
&\quad+C\epsilon^{\frac{p}{2}}\mathbb{E}\big[\big|D_x\Psi(x_s^\epsilon,y_s^\epsilon)\int_{s}^{t}(f(x_r^\epsilon,y_r^\epsilon)-f(x_s^\epsilon,y_s^\epsilon))\dd B_{r}^H\big|^p\big]\cr
&\quad+C\epsilon^{\frac{p}{2}}\mathbb{E}\big[\big|D_x\Psi(x_s^\epsilon,y_s^\epsilon)f(x_s^\epsilon,y_s^\epsilon)\int_{s}^{t}\dd B_{r}^H\big|^p\big]\cr
&=:K_1+K_2+K_3.
\end{aligned}
\end{equation}
So, according to assumption \textbf{A4}, Lemma \ref{clt1_x_bounded} and assumption that $f\in BC^2$, it deduces
\begin{equation}\label{K1}
K_1 \le C\epsilon^{\frac{p}{2}}\mathbb{E}\big[|(\|x^\epsilon\|_{\alpha}+\|y^\epsilon\|_{\theta})\|B^H\|_{\beta}|^p\big](t-s)^{(\theta+\beta)p}
\le C_{x_0,T}(t-s)^{(\theta+\beta)p}.
\end{equation}
where $ C_{x_0,T}>0$ is a constant only depending on $x_0$ and $T$.

Recall that according to the Proposition \ref{Poisson_g}, for all \((x,y)\in\mathbb{R}^m\times \mathbb{R}^n\), there exists a constant $L_\Psi>0$ such that 
$$|D_x\Psi(x,y)|\vee|D^2_{xx}\Psi(x,y)|\le L_\Psi(1+|x|+|y|).$$
So based on above result, assumption \textbf{A4}, Lemma \ref{clt1_x_bounded} and assumption that $f\in BC^2$, it is also straightforward to arrive at \(K_2\le C_{x_0,T}(t-s)^{(\alpha+\beta)p}\) and \(K_3\le C_{x_0,T}(t-s)^{\beta p}\).

Thus the proof is completed. \qed

\begin{lemma}\label{x_equicontinutiy0}
Let \(p\ge 2\) and \(\alpha\in(0,H)\) be such that \(1+\frac{1}{p}<\alpha+H\). Then for all \(\epsilon\in(0,1]\), $z^\epsilon\in \mathcal{B}_{\alpha, p}$.
\end{lemma}
\para{Proof}. It requires to show that for all \(\epsilon\in(0,1]\) and \(s,t\in[0,T]\)
\begin{eqnarray}\label{z1_equicontinutiy}
\mathbb{E}[\|z_t^\epsilon-z_s^\epsilon\|^p_{\infty}] \le C_{T}(t-s)^{\frac{p}{2}}
\end{eqnarray} 
where the constant \(C_{T}>0\) only depends on \(T\).

According to the Jensen inequality and straightforward computation \eqref{homogenization-Poisson}, 
\begin{eqnarray}\label{x_st}
|z_t^\epsilon -z^\epsilon_s|^p 
&\le& 6^{p-1} |\int_s^t f(x_r^\epsilon,y_r^\epsilon) \dd B^H_r|^p  + 6^{p-1}|\sqrt{\epsilon}\int_{s}^{t}D_x\Psi(x^\epsilon_r,y_r^\epsilon)g(x^\epsilon_r, y_r^\epsilon)dr|^p\cr
&& +6^{p-1}|\sqrt{\epsilon}\Psi(x^\epsilon_s,y^\epsilon_s)-\Psi(x^\epsilon_t,y_t^\epsilon)|^p +6^{p-1}|{\epsilon}\int_{0}^{t}D_x\Psi(x^\epsilon_r,y_r^\epsilon)f(x^\epsilon_r)\dd B_r^H|^p\cr
&& +6^{p-1}|\int_{s}^{t}\big(\int_{0}^{1}D\bar g(\theta x^\epsilon_r+(1-\theta)\bar x_r)\dd \theta\big)z^\epsilon_r \dd r|^p\cr
&& +6^{p-1}|\int_{s}^{t}D_y\Psi(x^\epsilon_r,y_r^\epsilon)\sigma(y^\epsilon_r)\dd W_r|^p\cr
&=:& I_1+I_2+I_3+I_4+I_5+|M_t^\epsilon|^p.
\end{eqnarray}
Next we estimate the terms on the right-hand side of \eqref{x_st}. From now on, we denote by $\mathbb{E}^B$ and $\mathbb{E}^W$ the expectations with respect to $B$ and $W$, respectively, so that $\mathbb{E}= \mathbb{E}^B \mathbb{E}^W$.

According to Lemma \ref{partial_x_phi}, it is easy to see that there are exponents \(\eta>\frac{1}{2}\) and \(\bar \eta>1\) such that for 
\( \mathbb{E}[I_1]\le C_1 6^{p-1} |f(x_r^\epsilon,y_r^\epsilon)|^p_{-\kappa,p}(|t-s|^{\eta p}+\|x\|_{\alpha}^p(t-s)^{\bar \eta p})
\) with some constant $C_1>0$.
Then, we have\\ 
\begin{equation}\label{I1}
\begin{aligned}
\mathbb{E}[I_1] 
& \le C_1 6^{p-1} \mathbb{E}[|f(x_r^\epsilon,y_r^\epsilon)|^p_{-\kappa,p}(|t-s|^{\eta p}+\|x\|_{\alpha}^p(t-s)^{\bar \eta p})]\\
& \le C_1 6^{p-1} L^p_f(|t-s|^{\eta p}+\|x\|_{\alpha}^p(t-s)^{\bar \eta p}),
\end{aligned}
\end{equation}
here, \(f\) is controlled by the positive constant \(L_f\) globally. 

Next, we estimate the second term \(I_2\). With aid of the Proposition \ref{Poisson_g}, for all \((x,y)\in \mathbb{R}^n\times \mathbb{R}^m\), \(D_x \Phi (x,y)\) is Lipschitz continuous. Then, by leveraging the H\"older inequality, we have
\begin{equation}\label{I2}
\begin{aligned}
\mathbb{E}[I_2] 
& \le 6^{p-1}\int_{s}^{t}\mathbb{E}\big[D_x\Psi(x^\epsilon_r,y_r^\epsilon)g(x^\epsilon_r, y_r^\epsilon)\big]^p \dd r(t-s)^{p-1}\\
& \le 6^{p-1}L^p_{D\Psi} \int_{s}^{t}\mathbb{E}\big[1+|x_r^\epsilon|+|y_r^\epsilon| dr\big]^p \dd r(t-s)^{p-1}\\
& \le 6^{p-1}L^p_{D\Psi} \sup_{r\in[0,t]}\mathbb{E}\big[1+|x_r^\epsilon|+|y_r^\epsilon|\big]^p(t-s)^p\\
& \le 6^{p-1}L^p_{D\Psi} C (1+C_y)(t-s)^p+6^{p-1}L^p_{D\Phi} C \mathbb{E}[\|x_r^\epsilon\|_{\infty,[s,t]}]^p(t-s)^p.
\end{aligned}
\end{equation}
Note that the second inequality comes from Proposition \ref{Poisson_g} with Lipschitz constant of \(D_{\Psi}\) which is denote by \(L_{D\Psi}>0\).
The fourth inequality is based on fundamental inequality with constant \(C>0\) independent of \(p\) and global boundedness of the fast component \(y\).

Then, we estimate the term \(I_3\). By using the assumption \textbf{A4} and Proposition \ref{Poisson_g}, then
\begin{equation}\label{I3_1}
\begin{aligned}
\mathbb{E}[I_3] 
& \le 6^{p-1} \epsilon^{\frac{p}{2}}\mathbb{E}[|x^\epsilon_t-x^\epsilon_s|+|y^\epsilon_t-y^\epsilon_s|]^p\\
& \le 6^{p-1} \epsilon^{\frac{p}{2}}(\mathbb{E}[\|x^\epsilon\|^{p}_\alpha(t-s)^{\alpha p}]+\mathbb{E}[|y^\epsilon_t-y^\epsilon_s|]^p)\\
& \le 6^{p-1}\epsilon^{\frac{p}{2}}\mathbb{E}[\|x^\epsilon\|^{p}_\alpha(t-s)^{\alpha p}]+C^26^{p-1}|t-s|^{p/2}
\end{aligned}
\end{equation}
where \(L_{\Psi}\) is the Lipschitz constant of \(\Psi\).

The estimate for the term \(I_4\) is given as below. According to Lemma \ref{partial_x_phi}, there are exponents \(\eta>\frac{1}{2}\) such that for \(f:\mathbb{R}^n\times \mathbb{R}^m\mapsto \mathbb{R}\) a bounded uniformly Lipschitz continuous function, we have
$$ \mathbb{E}[I_4] \le \|\int_{s}^{t}\sqrt{\epsilon}D_x\Psi(x_r^\epsilon,y_r^\epsilon)f(x_r^\epsilon,y_r^\epsilon)\dd B_{r}^H\|_{p} \le C_{x_0,T}6^{p-1}|t-s|^\eta
$$
where $ C_{x_0,T}>0$ is a constant only depending on $x_0$ and $T$.
With aid of the Proposition \ref{Poisson_g}, for all \((x,y)\in \mathbb{R}^n\times \mathbb{R}^m\), \(D_x\bar g\) is controlled by \(L_{g}\) globally. 
Then, we have
\begin{equation}\label{I5}
\mathbb{E}[I_5] 
\le 6^{p-1}C_1 \mathbb{E}[|D_x\bar g|^p_{sup}\|z\|_{sup,[s,t]}^p(t-s)^{p})]
\end{equation}
Finally, we will estimate the last term in \eqref{x_st}. According to the Burkholder-Davis-Gundy inequality and Proposition \ref{Poisson_g}, we immediately obtain the bound\\
\begin{equation}\label{M_st}
\begin{aligned}
\mathbb{E}[M_{s,t}]^p 
= 6^{p-1}\mathbb{E} \big[\int_{s}^{t}D_y\Psi(x^\epsilon_r,y_r^\epsilon)\sigma(y^\epsilon_r)\dd W_r\big]^p
& \le C6^{p-1}\mathbb{E}\big[\int_{s}^{t} (D_y\Psi(x^\epsilon_r,y_r^\epsilon)\sigma(y^\epsilon_r))^2dr \big]^{\frac{p}{2}}\\
& \le C6^{p-1}(t-s)^{\frac{p}{2}}.
\end{aligned}
\end{equation}
By combining \eqref{x_st}--\eqref{M_st}, it is trivial that with a constant \(C'\) controlling all constants here, for \(\eta>\frac{1}{2}\) and \(\bar \eta>1\)
\begin{equation}\label{xp_st}
\begin{aligned}
\mathbb{E}[z_t^\epsilon -z^\epsilon_s]^p 
&\le C'6^{p-1} (|t-s|^{\eta p}+\|x\|_{\alpha}^p(t-s)^{\bar \eta p})+C'6^{p-1}(t-s)^p\\
&\quad+C'6^{p-1}\mathbb{E}[\|x_r^\epsilon\|_{\infty,[s,t]}]^p|t-s|^p+ C'6^{p-1}|t-s|^{p/2}+ C_{x_0,T}6^{p-1}|t-s|^\eta\\
&\quad+C'6^{p-1}\|z\|_{sup,[s,t]}^p(t-s)^{p}
\end{aligned}
\end{equation}
Subsequently, we choose suitable \([s,t] \subset[0,T]\) such that \(|t-s|<1/2\). Then,
\begin{equation}\label{x_p_st}
\begin{aligned}
\|z^\epsilon\|_{\alpha, p} 
&\le C'\|z_s\|_p+C'\|x\|_{\alpha}(t-s)^{\bar \eta-\alpha,[s,t]}+\frac{1}{2}\|x^\epsilon\|_{\alpha,p}|t-s|^{\eta-\alpha}\\
&\quad+\frac{C'}{2}\|z^\epsilon\|_{\alpha,p}(t-s)^{1-\alpha}+\|x_r^\epsilon\|_{\infty,p,[s,t]}(t-s)^{1-\alpha}+C'.
\end{aligned}
\end{equation}
Then
\begin{eqnarray}\label{x_p_st}
&& \|z^\epsilon\|_{\alpha, p,[s,t]} \cr
&\le& C'\|z_s\|_p+2C'\|x\|_{\alpha}(t-s)^{\bar \eta-\alpha}+C'\|x^\epsilon\|_{\alpha,p}(t-s)^{\eta-\alpha}\cr
&&+\|x_r^\epsilon\|_{\infty,p,[s,t]}(t-s)^{1-\alpha}+C'.
\end{eqnarray}
By repeating \eqref{X_infty}--\eqref{X_infty_2} in Proposition \ref{boundness_x} and with the help of Lemma \ref{clt_x_bdd0},  the estimate \eqref{x_equicontinutiy0} is arrived.
The proof is completed.\qed

Then we introduce ``Little H\"older space" to avoid the separability issue afterwards. Set \(0<\alpha<1/2\), the space $\mathcal{C}^{\alpha}$ is not separable, but we could find a ``Little H\"older space" $ H^\alpha$
which is the space that for all $g\in \mathcal{C}^\alpha$, with the norm
\begin{equation}\label{little_Holder}
\lim_{\delta\to0+}\sup_{\substack{|t-s|\le \delta\\0 \leq s<t \le T}} \frac{|g(t)-g(s)|}{(t-s)^\alpha}=0. 
\end{equation}
Note $ H^\alpha$ is separable and \( H^\alpha=\overline{\bigcup_{\kappa>0}\mathcal{C}^{\alpha+\kappa}}\).
For any given \(0<\alpha<1/2\), we could find a slight large exponent \(\alpha+\kappa\) such that \(\alpha<\alpha+\kappa<1/2\), then our process takes values in \(\mathcal{C}^{\alpha+\kappa} \), directly, it also belongs to the space $ H^\alpha$.

Based on Lemma \ref{x_equicontinutiy0}, the following result directly arrived.\\
\begin{proposition}\label{homogenization_tight}
Set \(0<\alpha<1/2\) and let \(\mathbb{P}^\epsilon=\mathbb{P}\circ x^\epsilon\) with \(\epsilon\in(0,1]\) be the sequence of probability measures induced by \(z^\epsilon\) on \(H^\alpha\). 
By the tightness criterion established in \cite{lamperti1962convergence}, the sequence of probabilities \(\mathbb{P}^\epsilon, \epsilon\in(0,1]\) is tight in \(H^\alpha\).
\end{proposition}

\subsection{Weak limit}

\para{Proof of Theorem \ref{clt1}}.

With Prokhorov's theorem, there is a subsequence \(\{\epsilon_{n_k}\}_{k\ge 1}\) of any sequence \(\{\epsilon_n\}_{n\ge 1}\) that tends to \(0\) such that sequence of the law of \(z^\epsilon\) that is denoted by \(\mathbb{P}^{\epsilon_{n_k}}\) 
is relatively compact, that is it weakly converges to a some limit probability \(\mathbb{P}^0\) in \(H^\alpha\) with \(0<\alpha<1/2\). 
Then we introduce a new probability space which is still denoted by \((\Omega, \mathcal{F}, \mathbb{P})\)(for the sake of simplicity). 
By Skorokhod's theorem, there exists a subsequence of \((z^{\epsilon_{n_k}}, M^{\epsilon_{n_k}}, B^{H,{\epsilon_{n_k}}}, W^{\epsilon_{n_k}})_{k\ge 1}\) that takes values in \(H^\alpha\), such that 
for every \(k\ge 1\), \((z^{\epsilon_{n_k}}, B^{H,{\epsilon_{n_k}}}, W^{\epsilon_{n_k}})\) owns same law \(\mathbb{P}^{\epsilon_{n_k}}\) 
converges to some limit \((z, M, B^{H}, W)\) under uniform topology almost surely. 
Moreover, \((B^{H,{\epsilon_{n_k}}}(\omega), W^{\epsilon_{n_k}}(\omega))=(B^{H}(\omega), W(\omega))\) for every \(\omega\in \Omega\).

The remain task is only to show that the weak limit satisfies \eqref{clt-limit0} with uniform topology, as it implies the convergence with H\"older topology. To do this, we give some straight computation,
\begin{equation}\label{x_varepsilon}
\begin{aligned}
z_t^{\epsilon_{n_k}} -z_t
&= \int_0^t (f(x_r^{\epsilon_{n_k}},y_r^{\epsilon_{n_k}}) -\bar f(x_r^{\epsilon_{n_k}}) )\dd B^H_r +\int_0^t \big(\bar f(x_r^{\epsilon_{n_k}}) -\bar f(x_r)\big)\dd B^H_r\\
&\quad+\int_{0}^{t}\big[\big(\int_{0}^{1}D\bar g(\theta x^{\epsilon_{n_k}}_r+(1-\theta)\bar x_r)\dd \theta\big)z^{\epsilon_{n_k}}_r -D\bar g(\bar x_r)\dd \theta z_r\big]\dd r\\
&\quad+\sqrt{{\epsilon_{n_k}}}\big[ \Psi(x_0,y_0)-\Psi(x^{\epsilon_{n_k}}_t,y_t^{\epsilon_{n_k}})+\int_{0}^{t}D_x\Psi(x^{\epsilon_{n_k}}_r,y_r^{\epsilon_{n_k}})g(x^{\epsilon_{n_k}}_r, y_r^{\epsilon_{n_k}})\dd r\big.\\
&\quad\big.+\sqrt{{\epsilon_{n_k}}}\int_{0}^{t}D_x\Psi(x^{\epsilon_{n_k}}_r,y_r^{\epsilon_{n_k}})f(x^{\epsilon_{n_k}}_r,y_r^{\epsilon_{n_k}})\dd B_r^H \big]+M_t^{\epsilon_{n_k}}\\
&=:I_1^{\epsilon_{n_k}}+I_2^{\epsilon_{n_k}}+I_3^{\epsilon_{n_k}}+I_4^{\epsilon_{n_k}}+M_t^{\epsilon_{n_k}}.
\end{aligned}
\end{equation}
We apply with Theorem \ref{f_x_y} with \(h=f-\bar f\), which yields the following bound for \(\alpha<H\) 
\begin{equation}\label{f-barf}
\|\int_{0}^{\cdot}(f(x_r,y_r^{\epsilon_{n_k}})-\bar f(x_r))dB^H_r\|_{\alpha,p} \lesssim (\epsilon_{n_k})^\kappa(1+\|x\|_{\alpha,p})
\end{equation}
uniformly over \(x\in\mathbb{R}^n\) and \(\epsilon_{n_k}\) for all \(k\ge 1\). Here, \(\kappa>0\) is small enough. \(\alpha<1/2\) and \(p>1\) such that \(\alpha+H>1+\frac{1}{p}\).\\
According to Lemma \ref{clt1_x_bounded}, we have that \(\sup_{\epsilon_{n_k}}\|x^\epsilon\|_{\alpha,p}<\infty\). 
Then by combining with Theorem \ref{f_x_y}, we conclude that 
\begin{equation}\label{I_1}
\|\int_{0}^{\cdot}(f(x^\epsilon_r,y_r^{\epsilon_{n_k}})-\bar f(x^\epsilon_r))\dd B^H_r\|_{\alpha,p} \lesssim (\epsilon_{n_k})^\kappa.
\end{equation}
It is straight to see that \(\|I_1^{\epsilon_{n_k}}\|_{\infty} \to 0\) as \(k\) tends to infinity. 

Next, we estimate \(I_2^{\epsilon_{n_k}}\). First we recall that \(\bar f\in BC^2\) which comes from Remark \ref{bar_h}. Then, with direct computation, we have for \(\alpha<1/2\) and \(\beta<H\) such that \(\alpha+\beta>1\), 
\begin{equation}\label{I_2}
\begin{aligned}
\big|\int_0^t [\bar f(x_r^{\epsilon_{n_k}}) -\bar f(\bar x_r)]\dd B^H_r\big|&=\big|\int_0^t \int_{0}^{1}[ D\bar f(\theta x_r^{\epsilon_{n_k}}+(1-\theta)\bar x_r)(x_r^{\epsilon_{n_k}}-\bar x_r)]\dd \theta \dd B^H_r\big|\\
&\le C \big\|\int_{0}^{1}[ D\bar f(\theta x^{\epsilon_{n_k}}+(1-\theta)\bar x)(x^{\epsilon_{n_k}}-\bar x)] \big\|_{\alpha}\|B^H\|_{\beta}.
\end{aligned}
\end{equation}
Then, we will show the convergence under \(\alpha\) H\"older topology in \eqref{I_2}. We set \(J_r:=\int_{0}^{1}[ D\bar f(\theta x_r^{\epsilon_{n_k}}+(1-\theta)\bar x_r)(x_r^{\epsilon_{n_k}}-\bar x_r)] \), with Taylor's expansion,
\begin{equation}\label{J}
\begin{aligned}
J_r-J_s &= \int_{0}^{1}[D\bar f(\theta x_r^{\epsilon_{n_k}}+(1-\theta)\bar x_r)(x_r^{\epsilon_{n_k}}-\bar x_r-x_s^{\epsilon_{n_k}}+\bar x_s)] \dd \theta\\
&\quad+\big[\int_{0}^{1}[D\bar f(\theta x_r^{\epsilon_{n_k}}+(1-\theta)\bar x_r)-D\bar f(\theta x_s^{\epsilon_{n_k}}+(1-\theta)\bar x_s)](x_s^{\epsilon_{n_k}}-\bar x_s)\big]\dd \theta\\
&= \int_{0}^{1}[D\bar f(\theta x_r^{\epsilon_{n_k}}+(1-\theta)\bar x_r)(x_r^{\epsilon_{n_k}}-\bar x_r-x_s^{\epsilon_{n_k}}+\bar x_s)] \dd \theta\\
&\quad+\int_{0}^{1}\int_{0}^{1}[ D^2\bar f(\theta x_s^{\epsilon_{n_k}}+(1-\theta)\bar x_s+\theta' v)]v\dd \theta \dd \theta'(x_s^{\epsilon_{n_k}}-\bar x_s)
\end{aligned}
\end{equation}
where \(v:=\theta (x_r^{\epsilon_{n_k}}-x_s^{\epsilon_{n_k}})+(1-\theta)((x_r^{\epsilon_{n_k}}-x_s^{\epsilon_{n_k}}))\). Then, by combining with the result that \(x^{\epsilon_{n_k}}\) converges to some limit \(x\) under H\"older topology almost surely as \(k\) tends to infinity, so \eqref{I_2} converges to zero uniformly all over \(t\in[0,T]\).

Firstly, we divide the third term into two parts,
\begin{equation}\label{I_3}
\begin{aligned}
I_{3}^{\epsilon_{n_k}}
&= \int_{0}^{t}\big[\big(\int_{0}^{1}D\bar g(\theta x^{\epsilon_{n_k}}_r+(1-\theta)\bar x_r)\dd \theta\big)-D\bar g(\bar x_r) \big]z^{\epsilon_{n_k}}_r\dd r+\int_0^t \big[D\bar g(\bar x_r) \big( z^{\epsilon_{n_k}}_r-z_r\big)\big]\dd r\\
&=:I_{31}^{\epsilon_{n_k}}+I_{32}^{\epsilon_{n_k}}.
\end{aligned}
\end{equation}
Next, we estimate \(I_{31}^{\epsilon_{n_k}}\). First, we recall that \(\bar g\in BC^2\) which comes from Remark \ref{Poisson_g}. 
With Taylor's expansion, by combining with the result \(\bar g\in C_b^2\) that shown in Remark \ref{Poisson_g0}and Lemma \ref{x_equicontinutiy0}, it arrives at 
\begin{equation}\label{I_31}
\begin{aligned}
\mathbb{E}[|I_{31}^{\epsilon_{n_k}} |]
&= \mathbb{E} \big[\big|\int_{0}^{t}\big[\int_{0}^{1}\int_{0}^{1}D^2\bar g(\theta \theta' x_r^{\epsilon_{n_k}}+(1-\theta \theta')\bar x_r)(x_r^{\epsilon_{n_k}}-\bar x_r)\dd \theta \dd \theta'\big]z^{\epsilon_{n_k}}_r\dd r\big|\big]\\
& \le L_{D\bar g} \mathbb{E} \big[|\int_{0}^{t}(x_r^{\epsilon_{n_k}}-\bar x_r)z^{\epsilon_{n_k}}_r\dd r|\big]
\le L_{D\bar g} \sqrt{\epsilon}.
\end{aligned}
\end{equation}
Next, we estimate \(I_{32}^{\epsilon_{n_k}} \). First we recall that \(\bar g\in BC^2\) which comes from Remark \ref{Poisson_g0}.
Then, with direct computation, we have 
\begin{equation}\label{I_32}
\big\|I_{32}^{\epsilon_{n_k}}-\int_{0}^{t}D\bar g(x_r)z_r\dd r\big\|_\infty \to 0
\end{equation}
Take direct computation, we also conclude that \(\|I_4^{\epsilon_{n_k}}\|_\infty \to 0\) as \(k\) tends to infinity according to Lemma \ref{partial_x_phi} and the result that \(\Phi\in C^2\) shown in Remark \ref{Poisson_g0}.
\begin{equation}\label{I_4}
\begin{aligned}
I_4^{\epsilon_{n_k}} &= \sqrt{{\epsilon_{n_k}}}\big[ \Psi(x_0,y_0)-\Psi(x^{\epsilon_{n_k}}_t,y_t^{\epsilon_{n_k}})\big]+ \sqrt{{\epsilon_{n_k}}}\big[\int_{0}^{t}(D_x\Psi(x^{\epsilon_{n_k}}_r,y_r^{\epsilon_{n_k}})g(x^{\epsilon_{n_k}}_r, y_r^{\epsilon_{n_k}})\dd r\big]\\
&\quad +{{\epsilon_{n_k}}}\big[\int_{0}^{t}D_x\Psi(x^{\epsilon_{n_k}}_r,y_r^{\epsilon_{n_k}})f(x^{\epsilon_{n_k}}_r,y_r^{\epsilon_{n_k}})\dd B_r^H\big]\\
&=: I_{41}^{\epsilon_{n_k}}+I_{42}^{\epsilon_{n_k}}+I_{43}^{\epsilon_{n_k}}.
\end{aligned}
\end{equation}
Then, it deduces that \(\|I_{41}^{\epsilon_{n_k}}\|_\infty \to 0\) by Lemma \ref{x_equicontinutiy0}, assumption \textbf{A4} and the result that \(\Psi\in C^2\) that shown in Lemma \ref{Poisson_g0}.

Recall that according to the Remark \ref{Poisson_g0}, for all \((x,y)\in\mathbb{R}^m\times \mathbb{R}^n\), there exists a constant $L_\Psi>0$ such that 
$|D_x\Psi(x,y)|\vee|D^2_{xx}\Psi(x,y)|\le L_\Psi(1+|x|+|y|)$.
So based on above result and assumption \textbf{A4}, we have
\begin{equation}\label{I_42}
\begin{aligned}
\mathbb{E}[|I_{41}^{\epsilon_{n_k}}| ]&\le \sqrt{{\epsilon_{n_k}}}\mathbb{E}\big[\int_0^t(1+|x_r^{{\epsilon_{n_k}}}|+|y_r^{{\epsilon_{n_k}}}|)\dd r\big]\\
&\le \sqrt{{\epsilon_{n_k}}}\big[1+\mathbb{E}[\|x^{{\epsilon_{n_k}}}\|_{\sup, [0,t]}]+\mathbb{E}[\|y^{{\epsilon_{n_k}}}\|_{\sup,[0,t]}]\big]t\\
&\le \sqrt{{\epsilon_{n_k}}}L_{\Phi}(1+x+\mathbb{E}[\|x^{{\epsilon_{n_k}}}\|_{\alpha}]t+C_y)t.
\end{aligned}
\end{equation}
Thus $\|I_{42}^{\epsilon_{n_k}}\|_\infty \to 0$ as $k\to\infty$. Meanwhile, by Lemma \ref{partial_x_phi} and Lemma \ref{Poisson_g}, we also have $\|I_{43}^{\epsilon_{n_k}}\|_\infty \to 0$ as $k\to\infty$. Consequently, we conclude that $\|I_4^{\epsilon_{n_k}}\|_\infty \to 0$ as $k\to\infty$.

It remains to clarify the limit process of the final term \(M^{\epsilon_{n_k}}\).
Note that it is a continuous local martingale on the probability space \((\Omega, \mathcal{F}, \mathbb{P})\) with quadratic variational process as below
\begin{equation}\label{Mepsilon-quadratic}
\langle M^{\epsilon_{n_k}}\rangle_t=\int_{0}^{t} (D_y \Psi g)(D_y\Psi g)^*(x_r^{\epsilon_{n_k}},y_r^{\epsilon_{n_k}})\dd s, \quad t\in[0,T].
\end{equation}
Since \(M_t^{\epsilon_{n_k}}\) converges to \(M_t\) under \(\mathbb{P}\) almost surely, so \(M\) is also a continuous local martingale on the same probability space \((\Omega, \mathcal{F}, \mathbb{P})\). If we could show that the quadratic variational process of \(M\) is\\
\begin{equation}\label{M-quadratic}
\langle M\rangle_t=\int_{0}^{t} \overline{(D_y \Psi g)(D_y\Psi g)^*}(x_r)\dd s, \quad t\in[0,T],
\end{equation}
so according to the martingale representation theorem \cite[Theorem 4.5.1]{stroock2007multidimensional}, 
there exists another probability space \((\tilde \Omega, \mathcal{\tilde F}, \mathbb{\tilde P})\), filtration \(\{\mathcal{\tilde F}_t\}_{t\in[0,T]}\) and standard Brownian motion \(\tilde W\) defined on \((\Omega\times \tilde \Omega, \mathcal{ F} \times\mathcal{\tilde F}, \mathbb{ P}\times\mathbb{\tilde P})\) adapted to \(\mathcal{ F}_t\times \mathcal{\tilde F}_t\) such that
\begin{equation}\label{tildeM-quadratic}
\tilde M_t=\int_{0}^{t} (\overline{(D_y \Psi g)(D_y\Psi g)^*})^{\frac{1}{2}}(\tilde x_r)d\tilde W_s, \quad t\in[0,T],
\end{equation}
with \(\tilde M_t(\omega,\tilde \omega)=\tilde M_t(\omega)\) and \(\tilde x_r(\omega,\tilde \omega)=\tilde x_r(\omega)\) for all \((\omega,\tilde \omega)\in (\Omega\times \tilde \Omega)\).
So it remains to show the quadratic variational process of \(M\) is \eqref{M-quadratic}. To do this, it is sufficient to show that 
\begin{equation}\label{M-tildeM}
\lim_{k \to \infty} \mathbb{E}\big[\int_{0}^{t} (D_y \Psi g)(D_y\Psi g)^*(\tilde x_r^{\epsilon_{n_k}},\tilde y_r^{\epsilon_{n_k}})\dd s-\int_{0}^{t} \overline{(D_y \Psi g)(D_y\Psi g)^*}(\tilde x_r)\dd r\big]^2=0.
\end{equation}
Then it implies that there exists a sequence which still be denoted by \(\epsilon_{n_k}\) such that \(\langle M^{\epsilon_{n_k}}\rangle_t \to \langle M\rangle_t\) as \(k\) tends to infinity.
Note that \(M^{\epsilon_{n_k}}\otimes M^{\epsilon_{n_k}}-\langle M^{\epsilon_{n_k}}\rangle_t\) is a matrix-valued continuous martingale, then with aid of the Vitali convergence theorem, 
we conclude that for \(t\in[0,T]\),
\begin{equation}\label{M-martingale}
\begin{aligned}
& \mathbb{E}\big[\big(\tilde M_t\otimes \tilde M_t-\int_{0}^{t} \overline{(D_y \Psi g)(D_y\Psi g)^*}(\tilde x_r)dr\big)-\big(\tilde M_s\otimes \tilde M_s-\int_{0}^{s} \overline{(D_y \Psi g)(D_y\Psi g)^*}(\tilde x_r)\dd r\big)|\mathcal{\hat F}_s\big]\\
&=\lim_{k \to \infty}\mathbb{E}\big[\big(\tilde M^{\epsilon_{n_k}}_t\otimes \tilde M^{\epsilon_{n_k}}_t-\int_{0}^{t} (D_y \Psi g)(D_y\Psi g)^*(\tilde x^{\epsilon_{n_k}}_r,\tilde y^{\epsilon_{n_k}}_r)\dd r\big)\big.\\
&\quad\quad\quad\quad\big.-\big(\tilde M^{\epsilon_{n_k}}_s\otimes \tilde M^{\epsilon_{n_k}}_s-\int_{0}^{s} (D_y \Psi g)(D_y\Psi g)^*(\tilde x^{\epsilon_{n_k}}_r,\tilde y^{\epsilon_{n_k}}_r)\dd r\big)|\mathcal{\hat F}_s\big]\\
&=0
\end{aligned}
\end{equation}
under \(\mathbb{P}\) almost surely.

Hence, we could conclude that the quadratic variational process of \(\langle M\rangle_t\) is in the sense of \eqref{M-quadratic}. The remaining task is to show that \eqref{M-tildeM} holds. First, we divide it into two terms as below
\begin{equation}\label{M1M2}
\begin{aligned}
& \mathbb{E}\big[\int_{0}^{t} (D_y \Psi g)(D_y\Psi g)^*(\tilde x_r^{\epsilon_{n_k}},\tilde y_r^{\epsilon_{n_k}})dr-\int_{0}^{t} \overline{(D_y \Psi g)(D_y\Psi g)^*}(\tilde x_r)\dd r\big]^2\\
&=\mathbb{E}\big[\int_{0}^{t} (D_y \Psi g)(D_y\Psi g)^*(\tilde x_r^{\epsilon_{n_k}},\tilde y_r^{\epsilon_{n_k}})dr-\int_{0}^{t} \overline{(D_y \Psi g)(D_y\Psi g)^*}(\tilde x_r^{\epsilon_{n_k}})\dd r\big]^2\\
&\quad+\mathbb{E}\big[\int_{0}^{t} \overline{(D_y \Psi g)(D_y\Psi g)^*}(\tilde x_r^{\epsilon_{n_k}})dr-\int_{0}^{t} \overline{(D_y \Psi g)(D_y\Psi g)^*}(\tilde x_r)\dd r\big]^2\\
&=:M^{\epsilon_{n_k}}_1+M^{\epsilon_{n_k}}_2
\end{aligned}
\end{equation}
Then with Theorem \ref{f_x_y} by the integral against time not FBM, we conclude that \(M^{\epsilon_{n_k}}_1\) converges to $0$ as \(k\) tends to infinity. 
By combining with the assumption that \(g\in BC^2\) and the result that $D_y\Psi(x,y)$ is globally bounded which comes from Proposition \ref{Poisson_g}, it is directly to verify that \((D_y \Psi g)(D_y\Psi g)^*(\cdot)\) is globally Lipschitz. 
So it is straightforward to see that \(M^{\epsilon_{n_k}}_2\) converges to $0$ as \(k\) tends to infinity. The proof is completed.\qed

\section{Central Limit Theorem: Case 2} 

In this section we study the deviation of the slow one  in the following slow-fast system
\begin{equation}\label{simiplified_System}
\begin{cases}
dx^{\epsilon}_t = g(x^{\epsilon}_t, y^{\epsilon}_t)dt + f(x^{\epsilon}_t)dB^H_{t},\\[4pt]
dy^{\epsilon}_t = \frac{1}{\epsilon} b(x^{\epsilon}_t, y^{\epsilon}_t)dt + \frac{1}{\sqrt{\epsilon}}\sigma(x^{\epsilon}_t, y^{\epsilon}_t)dw_{t}.
\end{cases}
\end{equation}
from its effective limit. Here $B_t^H$ denotes a FBM with similarity exponent $H>\f 12$. 

While the effective average dynamics, obtained as $\epsilon\to 0$, is quite simpler than original one, it is given by: 
\begin{equation}\label{averaged-system}
\dd \bar x_t = \bar g(\bar x_t) \dd t+f(\bar x_t) \dd B^H_t,\ \bar x_0 = x.
\end{equation}
Setting
$$z_t^\epsilon =\f 1{\sqrt \epsilon } (x_t^\epsilon -\bar x_t),$$
then
$$\dd z_t^\epsilon =
\f 1{\sqrt \epsilon }(g(x_t^\epsilon,y^\epsilon_t)-\bar g(\bar x_t))\dd t+\f 1{\sqrt \epsilon }(f(x_t^\epsilon)-f(\bar x_t))\dd B_t^H.$$

\begin{remark}\label{property_f _g}
Under assumption \textbf{H1}-\textbf{H3}. Since we assume \(g\in L(\mathbb{R}^d,\mathbb{R}^n)\) is \(BC^2\) and globally bounded, then \(\bar g \in BC^1\) which is straightforward.
\end{remark}

\begin{theorem}\label{clt_Result}
Under assumptions \textbf{H1}--\textbf{H3} and \textbf{A4}. As \(\epsilon\to 0\), \(x^\epsilon\) weakly converges to \(\bar x\), that is 
\begin{equation}\label{clt2-limit}
\dd \bar z_t = D\bar g(\bar x_t) \bar z_t \dd t+Df(\bar x_t) \bar z_t \dd B^H_t+ \bar V^{\frac{1}{2}}(\bar x_t)\dd \hat W_t,\ \bar x_0 = x,
\end{equation}
where \(\bar x\) satisfies the ODE \eqref{averaged-system}, and \(V(x,y):=(D_y \Phi \sigma)(D_y \Phi \sigma)^{*}(x,y)\) with $\bar V^{\frac{1}{2}}(x)=\int_{\mathbb{R}^{m}} V^{\frac{1}{2}}(x,y)\mu^{x}(dy)$.
\end{theorem}

To show Theorem \ref{limit_solution}, the strategy in our proof is as below.
\begin{itemize}
\item Firstly, the deviation component could be rewritten as following, 
\begin{equation}\label{z-epsilon}
\begin{aligned}
z^\epsilon&=G^\epsilon +\int_{0}^{t}\big(\int_{0}^{1}D\bar g(\theta x^\epsilon+(1-\theta)\bar x_r)\dd \theta\big)z^\epsilon \dd r\\
&\quad+\int_{0}^{t}\big(\int_{0}^{1}Df(\theta x^\epsilon + (1-\theta)\bar x_r)\dd \theta\big)z^\epsilon_r \dd B_r^H
\end{aligned}
\end{equation}
where 
\begin{equation*}
G^\epsilon:=\frac{1}{\sqrt{\epsilon}}\int_{0}^{t}[g(x^\epsilon_r,y_r^\epsilon)-\bar g(x^\epsilon_r,y_r^\epsilon)]\dd r.
\end{equation*}
\item According to the Proposition \ref{Poisson_g}, \( \Psi(x,y) =\int_{0}^{\infty}\mathcal{P}_t (g(x, Y_t)-\bar g(x)) \dd t\)
is the solution to the following Poisson equation 
$$ -\mathcal{L}(y)\Psi(x,y) = g(x,y)-\bar g(x)$$ where \(\mathcal{L}^x\) is the infinitesimal generator of the fast process \(y^\epsilon\) with fixed-$x$. Moreover, \(\Psi\in C^2\).
\item Then, by the Young-It\^o formula \ref{Young-Ito-SDE}, 
\begin{equation}\label{psi_epsilon}
\begin{aligned}
\Psi(x^\epsilon_t,y_t^\epsilon)&=\Psi(x_0,y_0)+\int_{0}^{t}D_x\Psi(x^\epsilon_r,y_r^\epsilon)g(x^\epsilon_r, y_r^\epsilon)\dd r\\
&\quad+\int_{0}^{t}D_x\Psi(x^\epsilon_r,y_r^\epsilon)f(x^\epsilon_r)\dd B_r^H+\frac{1}{\epsilon}\int_{0}^{t}\mathcal{L}^x(y_r^\epsilon)\Psi(x^\epsilon_r, y_r^\epsilon)\dd r +\frac{1}{\sqrt{\epsilon}}M^\epsilon_t
\end{aligned}
\end{equation}
with \(M_t^\epsilon\) is a local martingale related to \(y^\epsilon\). 
\item Moreover, we rewrite \(G^\epsilon\) as below
\begin{equation}\label{G_epsilon}
\begin{aligned}
G^\epsilon_t &=\sqrt{\epsilon}\big[ \Psi(x_0,y_0)-\Psi(x^\epsilon_t,y_t^\epsilon)
+\int_{0}^{t}D_x\Psi(x^\epsilon_r,y_r^\epsilon)g(x^\epsilon_r,y_r^\epsilon)\dd r\big.\\
&\quad\big.+\int_{0}^{t}D_x\Psi(x^\epsilon_r,y_r^\epsilon)f(x^\epsilon_r)\dd B_r^H
+\int_{0}^{t}\mathcal{L}(y_r^\epsilon)\Psi(x^\epsilon_r, y_r^\epsilon)\dd r\big]
+M_t^\epsilon\\
&=:N_t^\epsilon+M_t^\epsilon.
\end{aligned}
\end{equation}
Note that \(\|x^\epsilon\|_{\alpha,p}\) is bounded with \(\frac{1}{2}<\alpha<H\) and \(y^\epsilon\) is globally bounded in \(L_p\). 
\item \(N^\epsilon\) is bounded with H\"older topology in \(L_p\) sense; (For \(N^\epsilon\), some estimation for \(\Phi\) are needed.)
\item \(M^\epsilon\) is bounded with H\"older topology in \(L_p\) sense.
\item With above results, we further rewrite \(z^\epsilon\) as below,
\begin{equation}\label{simiplified_System2}
\begin{aligned}
z^\epsilon_t&=N^\epsilon_t + M^\epsilon_t +\big[\int_{0}^{t}\big(\int_{0}^{1}D\bar g(\theta x^\epsilon_r+(1-\theta)\bar x_r)\dd \theta\big)z^\epsilon_r \dd r
-\int_{0}^{t}D\bar g(\bar x_r)z^\epsilon_r \dd r\big]\\
&\quad+\big[\int_{0}^{t}\big(\int_{0}^{1}Df(\theta x^\epsilon_r + (1-\theta)\bar x_r)d\theta\big)z^\epsilon_r \dd B_r^H
-\int_{0}^{t}Df(\bar x_r)z^\epsilon_r \dd B_r^H\big]\\
&\quad+\int_{0}^{t}D\bar g(\bar x_r)z^\epsilon_r \dd r+\int_{0}^{t}Df(\bar x_r)z^\epsilon_r \dd B_r^H\\
&=: N^\epsilon_t + M^\epsilon_t+A^\epsilon_t+B^\epsilon_t+\int_{0}^{t}Dg^\epsilon(\bar x_r)\dd z^\epsilon_r \dd r+\int_{0}^{t}Df(\bar x_r)z^\epsilon_r \dd B_r^H
\end{aligned}
\end{equation}
\item \(A^\epsilon\), \(B^\epsilon\) converges to \(0\) in H\"older topology in the \(L_p\) sense. 
\item It then remains to show that $z^\epsilon$ is bounded in the $L_p$-sense under the H\"older norm. By virtue of the embedding between different H\"older continuous spaces, the tightness of the laws of $z^\epsilon$ in the H\"older topology follows. (Here, the little H\"older space is required because it is separable and complete, whereas the usual H\"older space is not separable.)
\item Then, by Prokhorov's theorem, every sequence of the laws of $z^\epsilon$ admits a subsequence that converges weakly to some limiting probability measure. By Skorokhod's theorem, after passing to a new probability space (for simplicity), there exists a subsequence of $\{z^\epsilon\}_{0<\epsilon\le 1}$ that converges to some limit almost surely under the uniform topology.
\end{itemize}

\subsection{A priori estimates}
The following deterministic Gronwall-type lemma is a straightforward extension of \cite[Lemma 3.3]{li2025fluctuations} to include a drift term; we therefore omit the details of the proof for brevity.
\begin{proposition}\label{boundness_z_epsilon}
Let $(U, V)$ be Banach spaces and suppose $f \in C^\alpha([0, T], V)$, $B \in C^\beta([0, T], U)$ for some 
$(\beta \in (\frac{1}{2}, 1])$, $(0 < \alpha \le \gamma < 1)$ and $(\alpha + \beta > 1)$. Also suppose that \(A \in C^\gamma([0,T],L(V,(U,V)))\) and \(\hat A\in C^\gamma([0,T],L(V,V))\) for \(\alpha<\gamma<1\).
Then, if \(X : [0, T] \to V\)
satisfies the controlled Young differential equation
\begin{equation}\label{eq:gronwall-YDE}
Z_t = Z_0 + \int_0^t \hat A_sZ_s \dd s + \int_0^t A_s Z_s\dd B_s+ f_t.
\end{equation}
Then, 
\begin{equation}\label{zinfty-norm}
\|Z\|_{\infty,[0,T]} \le C_1 e^{C_2 (\|A\|_{\gamma}^{\frac{1}{\gamma}}+\|\hat A\|_{\gamma}^{\frac{1}{\gamma}}+\|B\|_{\beta}^{\frac{1}{\beta}})} (\|Z_0\|+\|f\|_{\alpha})
\end{equation}
and 
\begin{equation}\label{zalpha-norm}
\|Z\|_{\alpha,[0,T]} \le C_3 \|B\|_{\beta}^{\frac{1}{\beta}}((\|B\|_{\beta}+\|A\|_{\gamma}+\|\hat A\|_{\gamma})\|Z\|_\alpha+\|f\|_{\alpha})
\end{equation}
for some constants \(C_1, C_2, C_3 > 0\) depending on \(\alpha, \beta, \gamma, T\) and \(\|f\|_{BC^2}\). 
\end{proposition}

\begin{remark}\label{clt2_weaksolution}
With small extension of \cite[Lemma 3.6]{li2025fluctuations} with drift term, it is straight to see that a Picard iteration will converge to the solution of \eqref{eq:gronwall-YDE}. So in the case where $f$,$A$ $\hat A$ and $B$ in \eqref{boundness_z_epsilon} are random, the weak uniqueness of the solution by realizing that at each step.
\end{remark}

The existence of the weak solution and uniqueness in law implies the existence of strong solution, which comes from the Yamada-Watanabe theorem, like the approach that in \cite{Nualart03092008}. The weak existence comes from the \cite[Theorem 5.3]{Nualart03092008} directly. The strong uniqueness could be derived straightly from residue lemma \cite[Lemma]{li2025fluctuations} with adding an drift term and the fact that the solution to the average dynamics \eqref{averaged-system} is strongly unique. So, we just give the result below with omitting the proof and just assume the coefficient only depends on $t$.
\begin{theorem}\label{limit_solution}
Fix \(H\in(\frac{1}{2},1)\) and \(\beta<H\). Assume that \(G,H\) are Lipschitz continuous and globally bounded, \(F \) is $\gamma$-H\"older continuous \(\gamma<1\) such that \(\gamma+\beta>1\). Then, there exists a unique weak solution to the following Young-It\^o SDE
\begin{equation}\label{clt-limit}
\dd Z_t = G_tZ_t \dd t+F_t Z_t \dd B^H_t + H_t\dd \hat W_t,\ \bar Z_0 = z,
\end{equation}
where the first integral is defined in the Young sense and the last one is in the It\^o sense.
\end{theorem}

\subsection{Tightness criterion}
Before showing the tightness criterion, some a-prior estimates are necessary.
\begin{lemma}\label{clt_x_bounded}
Let \(p\ge 2\) and \(\gamma\in(\frac{1}{2},H)\). Under \textbf{H1}--\textbf{H3}. Then for all \(\epsilon\in(0,1]\),
\begin{eqnarray}\label{clt_x_bdd}
\mathbb{E}[\|x^\epsilon\|^p_{\gamma}] \le C_{x_0,T}, \quad \mathbb{E}[\|\bar x\|^p_{\gamma}] \le C_{x_0,T}
\end{eqnarray} 
where \(C_{x_0,T}>0\) only depends on \(x_0\) and \(T\).
\end{lemma}
\para{Proof}. Along the approach in \cite[Corollary 4.20]{Hairer-Li}, with the special condition that $f(x,y):=f(x)$ here, the above result is straightforward to arrive at. \qed

The following result can be directly obtained by extending \cite{Hairer-Li}, so proof details are omitted.
\begin{lemma}\label{x_convergence}
Let \(p\ge 2\) and \(\gamma\in(\frac{1}{2},H)\). Under \textbf{H1}--\textbf{H3}. As \(\epsilon\) tends to zero, there exists a \(\kappa>0\) such that
\begin{eqnarray}\label{clt_x_bdd}
\|x^\epsilon-\bar x\|_{\gamma,p} \le C_{x_0,T} \epsilon^\kappa 
\end{eqnarray} 
where \(C_{x_0,T}>0\) only depends on \(x_0\) and \(T\) and \(\bar x\) is the solution to the Young SDE \eqref{averaged-system}.
Moreover, \(\bar g\) is globally Lipschitz continuous.
\end{lemma}

\begin{remark}\label{Dg_property}
If Assumptions \textbf{H1}--\textbf{H3} hold, according to Proposition \ref{Poisson_CLT},
for every \(\hat g(x,y):=g(x,y)-\bar g(x)\), there exists a solution \(\Psi(x,\cdot)\) satisfying that \(\Psi(x,\cdot) \in BC^3\) for every \(x\in \mathbb{R}^n\) and \(u(\cdot,y)\in BC^1\) for every \(y\in \mathbb{R}^m\), to the Poisson equation \eqref{psi}.
Then, according to \cite[Lemma 3.2]{rockner2021strong}, it deduces that \(\bar g\in BC^2\). 
When there is no feedback from the slow dynamics into the fast one, it ensures that \(\bar g\in BC^2\) by only assuming \(D_{x} g( x,y)\) is Lipschitz with respect to \(x\in\mathbb{R}^n\) for each \(y\in\mathbb{R}^m\).
\end{remark}

Let $x\in B_{\alpha,p}$ with $p\ge 2$ and $\alpha\in(0,1/2)$ be such that $1+\frac{1}{p}<\alpha+H$. For a given $\mathcal{F}_t$-adapted process $x$ (not necessarily the solution to \eqref{simiplified_System}), we introduce a fast-varying process $\tilde y^{\epsilon}$ by
\begin{equation}\label{fast1}
d\tilde y^{\epsilon}_t = \frac{1}{\epsilon} b(x_t, \tilde y^{\epsilon}_t)\,dt + \frac{1}{\sqrt{\epsilon}}\sigma(x_t, \tilde y^{\epsilon}_t)\,dw_t, \qquad \tilde y^\epsilon_0=y
\end{equation}
where $t\in[0,T]$.
Denote its solution flow by $\Phi_{s,t}^x$ for $t>s$, with $\Phi_{s,s}^x=y$.
The following estimates then hold. These estimates differ substantially from those in Theorem \ref{f_x_y}, and a more elaborate proof is required.
\begin{lemma}\label{zpartial_x_phi}
Let \(p\ge 2\) and \(0<\alpha<\beta<H\) be such that \(1+\frac{1}{p}<\alpha+H\). There are exponents \(\eta>\frac{1}{2}\) and \(\bar \eta>1\) such that for \(f:\mathbb{R}^n\mapsto \mathbb{R}^m\times \mathbb{R}^d\) a bounded uniformly Lipschitz continuous function, we have
$$ \|\int_{s}^{t}D_x\Psi(x_r,\tilde y_r^\epsilon)f(x_r)\dd B_{r}^H\|_{\alpha,p} \le C_{x_0,T}|t-s|^\eta
$$
where \(C_{x_0,T}>0\) is a constant only depending on \(x_0\) and \(T\).
\end{lemma}
\para{Proof}. Firstly, we introduce the other fast varying process where $x$ is a given $\mathcal{F}_t$-adapted process (not necessarily the solution to \eqref{simiplified_System}),
\begin{equation}\label{fast2}
\dd Y_{s,t} = \frac{1}{\epsilon} b(x, Y_{s,t})\dd t + \frac{1}{\sqrt{\epsilon}}\sigma(x, Y_{s,t})\dd w_{t}, Y_{s,s}=y.
\end{equation}
Respectively, denote its solution flow by $\bar \Phi_{s,t}^x$ for any $t>s$ with $\Phi_{s,s}^x=y$. Hence, $Y_{s,t}:=\bar \Phi_{s,t}^x(\Phi_{s,t}^x(y_0))$.

According to the Proposition \ref{Poisson_CLT}, \(\Psi \in C^2\).
Moreover, for all \((x,y)\in \mathbb{R}^n\times \mathbb{R}^m\), \(D_x\Psi(x,y)\) is linear growth with respect to \((x,y)\in \mathbb{R}^n\times \mathbb{R}^m\). Along the approach as in Lemma \ref{partial_x_phi}, it is straightforward to see that 
$$ \|\sqrt{\epsilon}\int_{s}^{t}D_x\Psi(x_r,y_r^\epsilon)f(x_r^\epsilon)\dd B_{r}^H\|_{\alpha,p} \le C_{x_0,T}|t-s|^\eta
$$
where \(C_{x_0,T}>0\) is a constant only depending on \(x_0\) and \(T\).

However, sharper estimates are needed here, which necessitates the use of the Wiener--Young integral. Recall from Proposition \ref{Poisson_g} that for all $(x,y)\in\mathbb{R}^m\times \mathbb{R}^n$, there exists a constant $L_\Phi>0$ such that
$|D_x\Psi(x,y)|\le L_\Psi(1+|x|+|y|)$. Directly, we have for all $x_1,x_2\in \mathbb{R}^n$ and $y\in \mathbb{R}^m$, there exists a constant $C>0$ such that 
\begin{eqnarray}\label{psifx}
&&|D_{x}\Psi({x_1},y)f({x_1})-D_{x}\Psi({x_2},y)f({x_2})| \cr
&&\le |D_{x}\Psi({x_1},y)-D_{x}\Psi({x_2},y)||f({x_1})|+|D_{x}\Psi({x_2},y)||f({x_1})-f({x_2})|\cr
&&\le C(1+|x_1|+|x_2|+|y|)|x_1-x_2|.
\end{eqnarray}
and 
\begin{eqnarray}\label{psifxlinear}
|D_{x}\Psi({x_1},y)f({x_1})|\le C(1+|x_1|+|x_2|+|y|).
\end{eqnarray}
Meanwhile, we have for all $y_1,y_2\in \mathbb{R}^m$ and $x\in \mathbb{R}^n$, there exists a constant $C>0$ such that 
\begin{eqnarray}\label{psify}
&&|D_{x}\Psi({x},y_1)f({x})-D_{x}\Psi({x},y_2)f({x})|\cr
&& \le |D_{x}\Psi({x},y_1)-D_{x}\Psi({x},y_2)||f({x})|_{\sup}
\le L_{Dx\Phi f}|y_1-y_2|.
\end{eqnarray}
According to the Lemma \ref{sewing}, the integral is defined by $\lim_{| \mathcal{P}|\to 0} \sum_{[u,v]\in \mathcal{P}} A_{u,v}$ where \\
$$A_{s,t}=\int_s^t D_{x}\Psi({x_s},Y_{s,r})f({x}) \dd B_r^H.$$
Directly,
\begin{eqnarray}\label{deltaAst}
\delta A_{sut}=\int_u^t (D_{x}\Psi({x_s},Y_{s,u})f({x})-D_{x}\Psi({x_u},Y_{u,r})f({x})) \dd B_r^H.
\end{eqnarray}
Then, it requires to estimate $\|A\|_{\eta,p}$ and $\vertiii{A}_{\bar \eta,p}$ for $\eta>1/2$, $\bar \eta>1$ and $p\ge 1$.

First we aim to  bound the $\|A\|_{\eta,p}$. Let $F_s(r):=D_{x}\Psi({x_s},Y_{s,r})f({x})$. According to the \cite[Lemma 3.4]{Hairer-Li}, for $q>p$ and a small exponent $0<\kappa<H-\frac{1}{2}$,
\begin{equation}\label{Astp}
\|A_{s,t}\|_p \le C \||F_s|_{-\kappa}\|_q|t-s|^{H-\kappa}.
\end{equation}
Directly, with aid of the property \eqref{psifxlinear}, Lemma \ref{clt_x_bdd} and Assumption \textbf{A4}, we have 
\begin{equation}\label{Fs}
\begin{aligned}
\||F_s|_{-\kappa}\|_q &= \| \sup_{u<v}|u-v|^{\kappa-1}\int_{u}^vF_x(r)\dd r\|_q\\
&\le C\sup_{u<v}|u-v|^{\kappa-1}\int_{u}^v(1+\|x\|_q+\|Y_{s,r}\|_q)\dd r\\
&\le C_{x_0,y_0}
\end{aligned}
\end{equation}
where the constant $ C_{x_0,y_0,T}>0$ only depends on $x_0,y_0,T$. Then, we have that for every $p\ge 2$ and small exponent $\kappa>0$, $\{A_{s,t}\}\in H_\eta^p$ with $\eta=H-\kappa$.

Then it remains to show for every $p\ge 2$, $A_{s,t}\in \bar H_{\bar \eta}^p$ for $\bar \eta>1$. According to the Young-Wiener integral \eqref{integral_FBM}, it is straightforward to see 
\begin{equation}\label{deltaAst1}
\begin{aligned}
\delta A_{sut}&=\int_u^t (D_{x}\Psi({x_s},Y_{s,u})f({x})-D_{x}\Psi({x_u},Y_{u,r})f({x})) \dd \tilde B_r^{H,u}\\
&\quad+\int_u^t (D_{x}\Psi({x_s},Y_{s,u})f({x})-D_{x}\Psi({x_u},Y_{u,r})f({x})) \dd \bar B_r^{H,u}\\
&=\int_u^t (D_{x}\Psi({x_s},Y_{s,u})f({x})-D_{x}\Psi({x_u},Y_{u,r})f({x})) \dd \tilde B_r^{H,u}\\
&\quad+\int_u^t (D_{x}\Psi({x_s},Y_{s,u})f({x})-D_{x}\Psi({x_u},Y_{u,r})f({x})) \dd \dot{\bar B}_r^{H,u}
\end{aligned}
\end{equation}
Recall that this integral is the sum of the Wiener integral with respect to $\tilde B_r^{H,u}$ and the Riemann-Stieltjes integral against smooth process $\bar B_r^{H,u}$. Notably, $\tilde B_r^{H,u}$ is independent of filtration $\mathcal{F}_u \vee \mathcal{G}_t$, and the integrand in \eqref{deltaAst1} is measurable with respect to it, so the conditional expectation of the first term in the second line against filtration $\mathcal{F}_u$ will vanish. Let $D_x\Psi f(x_s, \cdot):=D_{x}\Psi({x_s},\cdot)f({x})$, then
\begin{equation}\label{deltaAst_tilde}
\begin{aligned}
\mathbb{E}[\delta A_{sut}|\mathcal{F}_u]&=\int_u^t \mathbb{E}[D_{x}\Psi({x_s},Y_{s,u})f({x})-D_{x}\Psi({x_u},Y_{u,r})f({x})|\mathcal{F}_u] \dot{\bar B}_r^{H,u}\dd r\\
&=\int_u^t (\mathcal{P}_{r-u}^{x_s}D_x\Psi f(x_s, \cdot)(Y_{s,u})-\mathcal{P}_{r-u}^{x_u}D_x\Psi f(x_u, \cdot)(y_u))\dot{\bar B}_r^{H,u}\dd r
\end{aligned}
\end{equation}
where the last line comes from the property that $Y_{s,r}=\bar \Phi_{s,r}^{x_s}(y_s)=\bar \Phi_{u,r}^{x_s}\bar \Phi_{s,u}^{x_s}(y_s)=\bar \Phi_{u,r}^{x_s}(Y_{s,u})$ and $Y_{u,r}=\bar \Phi_{u,r}^{x_u}(y_u)$.

Since it is difficult to get this bound directly and keep it separate from $\dot{\bar B}_r^{H,u}$, so we bound the conditional expectation against the filtration $\mathcal{F}_s\vee \mathcal{G}_u$. So
\begin{equation}\label{tildeF_s}
\begin{aligned}
&\mathbb{E}[\delta A_{sut}|\mathcal{F}_s\vee \mathcal{G}_u]\\
&=\int_u^t \mathbb{E}[D_{x}\Psi({x_s},Y_{s,u})f({x}_s)-D_{x}\Psi({x_u},Y_{u,r})f({x}_s)|\mathcal{F}_s\vee \mathcal{G}_u] \dot{\bar B}_r^{H,u}\dd r\\
&=\int_u^t (\mathcal{P}_{u,-s}^{x_s}\mathcal{P}_{r-u}^{x_s}D_x\Psi f(x_s, \cdot)(y_s)-\mathbb{E}[\mathcal{P}_{r-u}^{x_u}D_x\Psi f(x_u, \cdot)(y_u)|\mathcal{F}_s\vee \mathcal{G}_u])\dot{\bar B}_r^{H,u}\dd r.
\end{aligned}
\end{equation}
Then we give the other set $\mathcal{U}_s^u=\{F:\Omega\times \mathbb{R}^m \mapsto \mathbb{R}: F \quad \text{is bounded in the}\quad  L_p \quad \text{sense and} (\mathcal{F}_s\vee \mathcal{G}_u \otimes \mathcal{B}(\mathbb{R}^m)) \}$. Before estimating \eqref{tildeF_s}, we introduce two operators which will be used below. 

First, we introduce a bounded linear operator $\mathcal{Q}_{r,v}^x:\mathcal{U}_v^u \mapsto \mathcal{U}_r^u$ by following sense
\begin{equation}\label{q_st}
\begin{aligned}
(\mathcal{Q}_{r,v}^x F)(\omega, y):=\mathbb{E}[F(\cdot,\Phi_{r,v}^{x}(y, \cdot) )|\mathcal{F}_s\vee \mathcal{G}_u](\omega).
\end{aligned}
\end{equation}
According to the \cite[Lemma 4.23]{Hairer-Li}, it has for $s<r<v$,
\begin{equation}\label{q_cocycle}
\begin{aligned}
\mathcal{Q}_{s,r}^x \circ\mathcal{Q}_{r,v}^xF=\mathcal{Q}_{s,v}^xF.
\end{aligned}
\end{equation}
Since for any fixed $\bar x$, $\bar \Phi_{r,v}^{\bar x}(y)$ is independent of the filtration $\mathcal{G}_u$ for all $r\le v \le u$, then we define an operator $\hat{\mathcal{P}}_{r,v}^{\bar x}:\mathcal{U}_v^u \mapsto \mathcal{U}_r^u$ by following
\begin{equation}\label{hatp_st}
\begin{aligned}
(\hat{\mathcal{P}}_{r,v}^{\bar x} F)(\omega, y):=\mathbb{E}[(\mathcal{P}_{v-r}^{\bar x}F)(\cdot, y)|\mathcal{F}_r\vee \mathcal{G}_u](\omega).
\end{aligned}
\end{equation}
Note that if $F$ is measurable with respect to the filtration $\mathcal{F}_r\vee \mathcal{G}_u$, the two operators $\hat{\mathcal{P}}_{r,v}^{\bar x} $ and $\mathcal{P}_{v-r}^{\bar x}$ will coincide.

Then we rewrite \eqref{tildeF_s} as below
\begin{equation}\label{tildeF_s1}
\begin{aligned}
\mathbb{E}[\delta A_{sut}|\mathcal{F}_s\vee \mathcal{G}_u]
&=\int_u^t (\hat{\mathcal{P}}_{s,u}^{x_s} \mathcal{P}_{r-u}^{x_s}D_x\Psi f(x_s, \cdot)(y_s)-\mathcal{Q}_{s,u}^x\mathcal{P}_{r-u}^{x_u}D_x\Psi f(x_u, \cdot)(y_u) )\dot{\bar B}_r^{H,u}\dd r\\
&=\int_u^t \hat{\mathcal{P}}_{s,u}^{x_s}( \mathcal{P}_{r-u}^{x_s}D_x\Psi f(x_s, \cdot)(y_s)-\mathcal{P}_{r-u}^{x_u}D_x\Psi f(x_u, \cdot)(y_u) )\dot{\bar B}_r^{H,u}\dd r\\
&\quad+\int_u^t ((\hat{\mathcal{P}}_{s,u}^{x_s}-\mathcal{Q}_{s,u}^x)\mathcal{P}_{r-u}^{x_s}D_x\Psi f(x_u, \cdot)(y_s))\dot{\bar B}_r^{H,u}\dd r\\
&=:R_1+R_2
\end{aligned}
\end{equation}
We begin by bounding the term $R_1$. In view of Remark \ref{ergodicity} and the property of $D_x\Phi$ given in \eqref{psifx}, it follows that
\begin{equation}\label{DxPhif_x}
\begin{aligned}
&|\mathcal{P}_t^x D_x\Phi f(x, \cdot)-\mathcal{P}_t^{\bar x} D_x\Phi f(\bar x, \cdot)|\\
&\le C(1+|x|+\|Y_{s,u}\|_{p})|x-\bar x|.
\end{aligned}
\end{equation}
Then,
\begin{equation}\label{R1}
\begin{aligned}
|R_1|
&\le \mathbb{E}[C(1+|x_s|+|Y_{s,u}|)|x_s-x_u||\mathcal{F}_s\vee \mathcal{G}_u]\int_u^t|\dot{\bar B}_r^{H,u}|\dd r 
\end{aligned}
\end{equation}
With aid of Lemma \ref{clt_x_bounded} and assumption \textbf{A4}, choose $p',q'>1$ such that $\frac{1}{p'}+\frac{1}{q'}$, it has for $p\ge 2$ \\
\begin{equation}\label{R1}
\begin{aligned}
\|R_1\|_p
&\le \|\mathbb{E}[(1+|x_s|+|Y_{s,u}|)|x_s-x_u||\mathcal{F}_s\vee \mathcal{G}_u]\|_{pp'}\|\int_u^t|\dot{\bar B}_r^{H,u}|\dd r \|_{pq'}\\
&\le C (1+\|x\|_{\alpha, pp'}+\|Y_{s,u}\|_{pp'})|t-s|^\alpha\cdot \|\int_u^t|\dot{\bar B}_r^{H,u}|\dd r \|_{pq'}\\
& \le C (1+\|x\|_{\alpha, pp'}+\|Y_{s,u}\|_{pp'})|t-s|^{H+\alpha}
\end{aligned}
\end{equation}
where the final line comes from the $\|\dot{\bar B}_r^{H,u}\|_{q} \lesssim |r-u|^H$ for every $q>1$.

We next to bound the term $R_2$. Recall that 
\begin{equation}\label{hatP-Q}
\begin{aligned}
|(\hat{\mathcal{P}}_{s,u}^{x_s}-\mathcal{Q}_{s,u}^x)F|\lesssim\sqrt{\mathbb{E}[|x|_\beta^2|\mathcal{F}_s\vee\mathcal{G}_u]}|u-s|^\beta|F|_{Lip}
\end{aligned}
\end{equation}
for any $\beta<H$ which comes from \cite[Lemma 4.24]{Hairer-Li}.
Then, 
\begin{equation}\label{R2}
\begin{aligned}
\|R_2\|_p
&\le \big|\int_u^t ((\hat{\mathcal{P}}_{s,u}^{x_s}-\mathcal{Q}_{s,u}^x)\mathcal{P}_{r-u}^{x_s}D_x\Psi f(x_u, \cdot)(y_s))\dot{\bar B}_r^{H,u}\dd r\big|\\
& \lesssim \int_u^t\sqrt{\mathbb{E}[|x|_\beta^2|\mathcal{F}_s\vee\mathcal{G}_u]}|u-s|^\beta
|\mathcal{P}_{r-u}^{x_s}D_x \Phi f(x_s,\cdot)|_{Lip}|\dot{\bar B}_r^{H,u}|\dd r.
\end{aligned}
\end{equation}
With the Remark \ref{ergodicity}, assumption \textbf{A4} and Lemma \ref{clt_x_bounded}, it is straightforward to derive that 
\begin{equation}\label{PDPhif}
\begin{aligned}
|\mathcal{P}_{r-u}^{x_s}D_x \Phi f(x_s,\cdot)|_{Lip}
\le Ce^{-c(r-u)/\epsilon} |D_x \Phi f(x_s,\cdot)|_{Lip}\le Ce^{-c(r-u)/\epsilon} L_{Dx\Phi f}
\end{aligned}
\end{equation}
Here, $L_{Dx\Phi f}>0$ is a constant shown in \eqref{psify}.
Combining \eqref{deltaAst_tilde}--\eqref{PDPhif}, for every $p\ge 2$, we could see that $A_{s,t}\in \bar H_{\bar \eta}^p$ for $\bar \eta>1$. 
The proof is completed. \qed

\begin{lemma}\label{z_equicontinutiy}
Let \(p\ge 2\) and \(\alpha\in(0,\frac{1}{2})\) be such that \(1+\frac{1}{p}<\alpha+H\). Then for all \(\epsilon\in(0,1]\), $z^\epsilon\in \mathcal{B}_{\alpha, p}$.
\end{lemma}
\para{Proof}. It requires to show that for all \(\epsilon\in(0,1]\) and \(s,t\in[0,T]\)
\begin{eqnarray}\label{z_equicontinutiy}
\mathbb{E}[\|z_t^\epsilon-z_s^\epsilon\|^p_{\infty}] \le C_{T}(t-s)^{\frac{p}{2}}
\end{eqnarray} 
where the constant \(C_{T}>0\) only depends on \(T\).

First, it is straight to obtain the below equality with Young-It\^o formula, 
\begin{equation}\label{z_st}
\begin{aligned}
z^\epsilon_t-z^\epsilon_s&=N^\epsilon_{s,t} + M^\epsilon_{s,t} +\big[\int_{s}^{t}\big(\int_{0}^{1}D\bar g(\theta x^\epsilon_r+(1-\theta)\bar x_r)\dd \theta\big)z^\epsilon_r \dd r\big]\\
&\quad+\big[\int_{s}^{t}\big(\int_{0}^{1}Df(\theta x^\epsilon_r + (1-\theta)\bar x_r)\dd \theta\big)z^\epsilon_r \dd B_r^H\big]\\
&=: N^\epsilon_{s,t} + M^\epsilon_{s,t}+A^\epsilon_{s,t}+B^\epsilon_{s,t}.
\end{aligned}
\end{equation}
We denote that $\mathbb{E}^B$, $\mathbb{E}^W$ are the expectation with respect to $B$ and $W$ respectively, so it has that $\mathbb{E}= \mathbb{E}^B\times \mathbb{E}^W$. 

Firstly, we intend to show that the term \(N^\epsilon\) is bounded with H\"older topology in $L_p$ sense. So,
\begin{equation}\label{N_epsilon}
\begin{aligned}
N_{s,t}^\epsilon &=\sqrt{\epsilon}\big[ \Psi(x^\epsilon_s,y^\epsilon_s)-\Psi(x^\epsilon_t,y_t^\epsilon)\big]
+\sqrt{\epsilon}\big[\int_{s}^{t}D_x\Psi(x^\epsilon_r,y_r^\epsilon)g(x^\epsilon_r,y_r^\epsilon)\dd r\big]\\
&\quad+\sqrt{\epsilon}\big[\int_{s}^{t}D_x\Psi(x^\epsilon_r,y_r^\epsilon)f(x^\epsilon_r)\dd B_r^H\big]
\\
&=: N_{s,t}^{1,\epsilon}+ N_{s,t}^{2,\epsilon}+ N_{s,t}^{3,\epsilon}.
\end{aligned}
\end{equation}
For the first term \(N_{s,t}^{1,\epsilon}\), similar to \eqref{I3_1}, by using the assumption \textbf{A4} and Proposition \ref{Poisson_CLT}, then we have
\begin{equation}\label{N_st^1}
\begin{aligned}
\mathbb{E}[|N_{s,t}^{1,\epsilon}|^p] 
\le \epsilon^{\frac{p}{2}}L_{\Psi}\mathbb{E}[|x^\epsilon_t-x^\epsilon_s|+|y^\epsilon_t-y^\epsilon_s|]^p
&\le \epsilon^{\frac{p}{2}}L_{\Psi}(\mathbb{E}[\|x^\epsilon\|^{p}_\alpha(t-s)^{\alpha p}]+\mathbb{E}[|y^\epsilon_t-y^\epsilon_s|]^p)\\
&\le \epsilon^{\frac{p}{2}}L_{\Psi}\mathbb{E}[\|x^\epsilon\|^{p}_\alpha(t-s)^{\alpha p}]+C^2_pL_{\Psi}(t-s)^{\frac{p}{2}}
\end{aligned}
\end{equation}
where \(L_{\Psi}\) is the Lipschitz constant of \(\Psi\).

Then we estimate the second term \(N_{s,t}^{2,\epsilon}\), due to the assumption \textbf{A4} and result that \(D_x\Psi\) is linear growth with respect to \((x,y)\in \mathbb{R}^n\times \mathbb{R}^m\) shown in Proposition \ref{Poisson_CLT}, it arrives at
\begin{equation}\label{N_st^2}
\begin{aligned}
\mathbb{E}[|N_{s,t}^{2,\epsilon}|^p] 
& \le L_{D\Psi} \int_{s}^{t}\mathbb{E}\big[1+|x_r^\epsilon|+|y_r^\epsilon| \dd r\big]^p \dd r(t-s)^{p-1}\\
& \le L_{D\Psi} \sup_{r\in[0,t]}\mathbb{E}\big[1+|x_r^\epsilon|+|y_r^\epsilon|\big]^p(t-s)^p\\
& \le L_{D\Psi} C_p (1+C_y)(t-s)^p+L_{D\Phi} C_p\mathbb{E}[\|x_r^\epsilon\|_{\infty,[s,t]}]^p(t-s)^p.
\end{aligned}
\end{equation}
where \(L_{D\Psi}\) is the Lipschitz constant of \(D\Psi\).

With the help of Lemma \ref{zpartial_x_phi}, we readily derive
\[
\mathbb{E}\left[|N_{s,t}^{3,\epsilon}|^p\right] \le C_{x_0,T}|t-s|^\eta.
\]
Combining the preceding estimates, we conclude that the term $N^\epsilon$ is bounded in the $L_p$ sense with H\"older regularity of order $\alpha<1/2$.

We now turn to bound  the $M^\epsilon$ in the H\"older topology in the $L_p$ sense. The argument is similar to that used in the derivation of \eqref{M_st} in Lemma \ref{z_equicontinutiy}. Indeed, by applying the Burkholder--Davis--Gundy inequality, Proposition \ref{Poisson_CLT}, Assumption \textbf{A4}, and Lemma \ref{clt_x_bounded}, we obtain the bound immediately.
\begin{eqnarray}\label{M_stepsilon}
\mathbb{E}[M^\epsilon_{s,t}]^p 
&=& \mathbb{E} \big[\int_{s}^{t}D_y\Psi(x^\epsilon_r,y_r^\epsilon)\sigma(x^\epsilon_r,y^\epsilon_r)\dd W_r\big]^p\cr
& \le &C\mathbb{E}\big[\int_{s}^{t} (D_y\Psi(x^\epsilon_r,y_r^\epsilon)\sigma(x^\epsilon_r,y^\epsilon_r))^2\dd r \big]^{\frac{p}{2}}\le C(t-s)^{\frac{p}{2}}.
\end{eqnarray}
To obtain the desired result \eqref{z_equicontinutiy}, we intend to apply the Gronwall-type lemma stated in Proposition \ref{eq:gronwall-YDE}. It therefore remains to verify that both $D\bar g$ and $Df$ are $\alpha$-H\"older continuous in the $L_p$ sense.
First, by Remark \ref{Dg_property} and Lemma \ref{clt_x_bounded}, it is straightforward to verify that $D\bar g$ is $\alpha$-H\"older continuous in the $L_p$ sense. Similarly, by employing the assumption $f\in BC^2$ together with Lemma \ref{clt_x_bounded}, we obtain that $Df$ is also $\alpha$-H\"older continuous in the $L_p$ sense.
The proof is completed.\qed

Based on Lemma \ref{z_equicontinutiy}, the following result is directly arrived.
\begin{proposition}\label{z_tight}
Set \(0<\alpha<1/2\) and let \(\mathbb{\hat P}^\epsilon=\mathbb{P}\circ z^\epsilon\) with \(\epsilon\in(0,1]\) be the sequence of probability measures induced by $z^\epsilon$ on \(H^\alpha\).
By the tightness criterion established in \cite{lamperti1962convergence}, the sequence of probabilities \(\mathbb{\hat P}^\epsilon, \epsilon\in(0,1]\) is tight in \(H^\alpha\).
\end{proposition}

\subsection{Weak limit}
\para{Proof of Theorem \ref{clt_Result}}.\\

With aid of Prokhorov's theorem, there is a subsequence \(\{\epsilon_{n_k}\}_{k\ge 1}\) of any sequence \(\{\epsilon_n\}_{n\ge 1}\) that tends to \(0\) such that sequence of the law of \(z^\epsilon\) that is denoted by \(\mathbb{\hat P}^{\epsilon_{n_k}}\) 
is relatively compact, that is it weakly converges to some limit probability \(\mathbb{\hat P}^0\) in \(H^\alpha\) with \(0<\alpha<1/2\). 
Then we introduce a new probability space which is still denoted by \((\Omega, \mathcal{F}, \mathbb{P})\) (for the sake of simplicity). 
By Skorokhod's theorem, there exists a subsequence of \((z^{\epsilon_{n_k}}, M^{\epsilon_{n_k}}, B^{H,{\epsilon_{n_k}}}, W^{\epsilon_{n_k}})_{k\ge 1}\) that takes values in \(H^\alpha\), such that 
for every \(k\ge 1\), \((z^{\epsilon_{n_k}}, B^{H,{\epsilon_{n_k}}}, W^{\epsilon_{n_k}})\) owns same law \(\mathbb{\hat P}^{\epsilon_{n_k}}\) 
converges to some limit \((z, M, B^{H}, W)\) under uniform topology almost surely. 
Moreover, \((B^{H,{\epsilon_{n_k}}}(\omega), W^{\epsilon_{n_k}}(\omega))=(B^{H}(\omega), W(\omega))\) for every \(\omega\in \Omega\).

The remaining task is only to verify that the weak limit satisfies \eqref{clt-limit} with uniform topology which implies the convergence with H\"older topology. To do this, we rewrite \(z^{\epsilon_{n_k}}\) below,
\begin{equation}\label{z_varepsilon}
\begin{aligned}
z_t^{\epsilon_{n_k}} 
&=\int_{0}^{t}\big[\int_{0}^{1}\big(D\bar g(\theta x^{\epsilon_{n_k}}_r + (1-\theta)\bar x_r)-D\bar g(\bar x_r)\big)\dd \theta\big]z^{\epsilon_{n_k}}_r \dd r+\int_{0}^{t}D\bar g(\bar x_r)z^{\epsilon_{n_k}}_r \dd r\\
&\quad+\int_{0}^{t}\big[\int_{0}^{1}\big(Df(\theta x^\epsilon_r + (1-\theta)\bar x_r)-Df(\bar x_r)\big)\dd \theta\big]z^{\epsilon_{n_k}}_r \dd B_r^H+\int_{0}^{t}Df(\bar x_r)z^{\epsilon_{n_k}}_r \dd B_r^H\\
&\quad+\sqrt{{\epsilon_{n_k}}}\big[ \Psi(x_0,y_0)-\Psi(x^{\epsilon_{n_k}}_t,y_t^{\epsilon_{n_k}})+\int_{0}^{t}D_x \Psi(x^{\epsilon_{n_k}}_r,y_r^{\epsilon_{n_k}}) g(x^{\epsilon_{n_k}}_r,y_r^{\epsilon_{n_k}})\dd r\big.\\
&\quad \big.+\int_{0}^{t}D_x\Psi(x^{\epsilon_{n_k}}_r,y_r^{\epsilon_{n_k}})f(x^{\epsilon_{n_k}}_r)\dd B_r^H\big]+M_t^{\epsilon_{n_k}}\\
&=:L_1^{\epsilon_{n_k}}+L_2^{\epsilon_{n_k}}+L_3^{\epsilon_{n_k}}+L_4^{\epsilon_{n_k}}+L_5^{\epsilon_{n_k}}+M_t^{\epsilon_{n_k}}.
\end{aligned}
\end{equation}
Next, we estimate $L_1^{\epsilon_{n_k}}$. We first recall from Remark \ref{Poisson_g} that $\bar g\in BC^2$. By Taylor's expansion, together with the fact that $\bar g\in BC^2$ (see Remark \ref{Poisson_g}) and Lemma \ref{z_equicontinutiy}, we obtain 
\begin{equation}\label{L_1}
\begin{aligned}
\mathbb{E}[|L_1^{\epsilon_{n_k}} |]
&= \mathbb{E} \big[\big|\int_{0}^{t}\big[\int_{0}^{1}\int_{0}^{1}D^2\bar g(\theta \theta' x_r^{\epsilon_{n_k}}+(1-\theta \theta')\bar x_r)(x_r^{\epsilon_{n_k}}-\bar x_r)\dd \theta \dd \theta'\big]z^{\epsilon_{n_k}}_r\dd r\big|\big]\\
& \le L_{D\bar g} \mathbb{E} \big[|\int_{0}^{t}(x_r^{\epsilon_{n_k}}-\bar x_r)z^{\epsilon_{n_k}}_r\dd r|\big]
\le L_{D\bar g} \sqrt{\epsilon}.
\end{aligned}
\end{equation}
Next, we estimate \(L_2^{\epsilon_{n_k}} \). First we recall that \(\bar g\in BC^2\) which comes from Remark \ref{Dg_property}.
Then, with direct computation, we have 
\begin{equation}\label{L_2}
\big\|L_2^{\epsilon_{n_k}}-\int_{0}^{t}D\bar g(x_r)z_r\dd r\big\|_\infty \to 0
\end{equation}
Then, it remains to  estimate \(L_3^{\epsilon_{n_k}}\). Due to the result that \(\bar g\in BC^2\) which comes from Remark \ref{Dg_property}, and  some straightforward computation, we have for \(\alpha<1/2\) and \(\beta<H\) such that \(\alpha+\beta>1\), 
\begin{equation}\label{L_3}
\begin{aligned}
\mathbb{E}[|L_3^{\epsilon_{n_k}} |]
&=\mathbb{E}\big[\big|\int_{0}^{t}\big[\int_{0}^{1}\big(Df(\theta x^\epsilon_r + (1-\theta)\bar x_r)-Df(\bar x_r)\big)\dd \theta\big]z^{\epsilon_{n_k}}_r\dd B^H_r\big|\big]\\
&\le C \mathbb{E}\big[\big\|\big[\int_{0}^{1}\big(Df(\theta x^\epsilon + (1-\theta)\bar x)-Df(\bar x)\big)\dd \theta\big]\big\|_{\alpha}\|z^{\epsilon_{n_k}} \|_{\alpha}\|B^H\|_{\beta}\big].
\end{aligned}
\end{equation}
Then, we will show the convergence under \(\alpha\) H\"older topology in \eqref{L_3}. We set \\
\begin{equation}\label{J'}
J_r:=\int_{0}^{1}\big(Df(\theta x_r^\epsilon + (1-\theta)\bar x_r)-Df(\bar x_r)\big)\dd \theta
=\int_{0}^{1}\int_{0}^{1}D^2f(\theta\theta' x_r^\epsilon + (1-\theta\theta')\bar x_r)(x_r^\epsilon-\bar x_r)\dd \theta \dd \theta'
\end{equation}
with Taylor's expansion
\begin{equation}\label{L_3'}
\begin{aligned}
J_r-J_s &= \int_{0}^{1}\int_{0}^{1}[D^2 f(\theta\theta' x_r^{\epsilon_{n_k}}+(1-\theta\theta')\bar x_r)\theta\theta'(x_r^{\epsilon_{n_k}}-\bar x_r-x_s^{\epsilon_{n_k}}+\bar x_s)] \dd \theta\dd \theta'\\
&\quad+\big[\int_{0}^{1}\int_{0}^{1}[D^2 f(\theta x_r^{\epsilon_{n_k}}+(1-\theta)\bar x_r)-D^2 f(\theta x_s^{\epsilon_{n_k}}+(1-\theta)\bar x_s)]\theta\theta'(x_s^{\epsilon_{n_k}}-\bar x_s)\big]\dd \theta\dd\theta'\\
&= \int_{0}^{1}\int_{0}^{1}[D^2 f(\theta\theta' x_r^{\epsilon_{n_k}}+(1-\theta\theta')\bar x_r)\theta\theta'(x_r^{\epsilon_{n_k}}-\bar x_r-x_s^{\epsilon_{n_k}}+\bar x_s)] \dd \theta\dd \theta'\\
&\quad+\int_{0}^{1}\int_{0}^{1}\int_{0}^{1}[ D^3 f(\theta\theta' x_s^{\epsilon_{n_k}}+(1-\theta\theta')\bar x_s+\theta'' v)]v\dd \theta \dd \theta'\dd \theta''(x_s^{\epsilon_{n_k}}-\bar x_s)
\end{aligned}
\end{equation}
where \(v:=\theta\theta' (x_r^{\epsilon_{n_k}}-x_s^{\epsilon_{n_k}})+(1-\theta\theta')(x_r^{\epsilon_{n_k}}-x_s^{\epsilon_{n_k}})\). Then, by combining with Lemma \ref{x_convergence}, so as \(k\) tends to infinity, \eqref{L_3} converges to zero uniformly all over \(t\in[0,T]\).

For the fourth term \(L_4^{\epsilon_{n_k}}\), according to the straight computation, 
\begin{equation}\label{L_4}
\begin{aligned}
\mathbb{E}[|L_4^{\epsilon_{n_k}} |]
&=\mathbb{E}\big[\big|\int_{0}^{t}Df(\bar x_r)(z^{\epsilon_{n_k}}_r-z_r)\dd B^H_r\big|\big]\le C \mathbb{E}\big[\|Df(\bar x)\|_{\alpha}\|z^{\epsilon_{n_k}}-z \|_{\alpha}\|B^H\|_{\beta}\big].
\end{aligned}
\end{equation}
It is then straightforward to see that \eqref{L_4} converges to zero as $k\to\infty$. Moreover, the fifth term in \eqref{z_varepsilon} vanishes in the limit $k\to\infty$ by an argument analogous to that used for $I_4^{\epsilon_{n_k}}$ in \eqref{x_varepsilon}. For the last term $M^{\epsilon_{n_k}}$, the same estimates as in \eqref{Mepsilon-quadratic}--\eqref{M1M2} yield the desired result.
The proof is completed. \qed

\begin{remark}
Consider the following slow-fast system,
\begin{align}\label{slowfastmore}
\begin{cases}
dx^{\epsilon,\delta}_t = g(x^{\epsilon,\delta}_t, y^{\epsilon,\delta}_t)dt + \sqrt{\epsilon} f(x^{\epsilon,\delta}_t)dB^H_{t},\\[4pt]
dy^{\epsilon,\delta}_t = \frac{1}{\delta} b(x^{\epsilon,\delta}_t, y^{\epsilon,\delta}_t)dt + \frac{1}{\sqrt{\delta}}\sigma(x^{\epsilon,\delta}_t, y^{\epsilon,\delta}_t)dw_{t}.
\end{cases}
\end{align}
Here we assume $\delta=o(\epsilon)$.  As \(\epsilon\) tends to zero, the above system converges to the effective averaged system under \(\alpha\)-H\"older topology \((\alpha<H)\) \cite{Hairer-Li}, it is given by: 
\begin{equation}\label{averaged-system}
\dd \bar x_t = \bar g(x_t) \dd t,\ \bar x_0 = x.
\end{equation}
Here, \(\bar g = \int_{\mathbb{R}^{m}} g(x,y)\mu^x(\dd y)\)\(y\) where is \(\mu(\dd y)\) is the invariant probability measure of the stationary fast It\^o process with fixed \(x\). 
Setting
$$z_t^{\epsilon} =\f 1{\sqrt \epsilon } (x_t^{\epsilon,\delta} -\bar x_t).$$
Under \textbf{A1}--\textbf{A4}. As \(\epsilon\to 0\), \(z^\epsilon\) weakly converges to \(\bar z\), that is 
\begin{equation}\label{clt-limit1}
\dd \bar z_t = D\bar g(\bar x_t) \bar z_t \dd t+ f(\bar x_t) \dd B^H_t.
\end{equation}
Note that there is no extra Gaussian process, that is different from Theorem \ref{clt_Result}. The reason is as below,
\begin{itemize}
\item The deviation component could be rewritten as following, 
\begin{equation}\label{z-epsilon1}
\begin{aligned}
z^\epsilon&=G^\epsilon +\int_{0}^{t}\big(\int_{0}^{1}D\bar g(\theta x^{\epsilon,\delta}+(1-\theta)\bar x_r)\dd \theta\big)z^\epsilon \dd r+\int_{0}^{t}\big(\int_{0}^{1}Df(\theta x^{\epsilon,\delta} + (1-\theta)\bar x_r)\dd \theta\big)z^\epsilon_r \dd B_r^H
\end{aligned}
\end{equation}
where
\begin{equation*}
G^\epsilon:=\frac{1}{\sqrt{\epsilon}}\int_{0}^{t}[g(x^{\epsilon,\delta}_r,y_r^{\epsilon,\delta})-\bar g(x^{\epsilon,\delta}_r,y_r^{\epsilon,\delta})]\dd r.
\end{equation*}
\item According to the Proposition \ref{Poisson_g}, \( \Psi(x,y) =\int_{0}^{\infty}\mathcal{P}_t (g(x, Y_t)-\bar g(x)) \dd t\) is the solution to the Poisson equation $ -\mathcal{L}(y)\Psi(x,y) = g(x,y)-\bar g(x)$ where \(\mathcal{L}^x\) is the infinitesimal generator of the fast process \(y^{\epsilon,\delta}\) with fixed-$x$. Moreover, \(\Psi\in C^2\).
\item Then, by the Young-It\^o formula \ref{Young-Ito-SDE}, 
\begin{equation}\label{psi_epsilon1}
\begin{aligned}
\Psi(x^{\epsilon,\delta}_t,y_t^{\epsilon,\delta})&=\Psi(x_0,y_0)+\int_{0}^{t}D_x\Psi(x^{\epsilon,\delta}_r,y_r^{\epsilon,\delta})g(x^{\epsilon,\delta}_r, y_r^{\epsilon,\delta})\dd r\\
&\quad+\int_{0}^{t}D_x\Psi(x^{\epsilon,\delta}_r,y_r^{\epsilon,\delta})f(x^{\epsilon,\delta}_r)\dd B_r^H+\frac{1}{\delta}\int_{0}^{t}\mathcal{L}^x(y_r^{\epsilon,\delta})\Psi(x^{\epsilon,\delta}_r, y_r^{\epsilon,\delta})\dd r \\
&\quad+\frac{1}{\sqrt{\delta}}M^\epsilon_t
\end{aligned}
\end{equation}
with \(M_t^\epsilon\) is a local martingale related to \(y^\epsilon\). Then the term \(G^\epsilon\) is rewritten as below
\begin{equation}\label{G_epsilon1}
\begin{aligned}
G^\epsilon_t &=\frac{\delta}{\sqrt{\epsilon}}\big[ \Psi(x_0,y_0)-\Psi(x^{\epsilon,\delta}_t,y_t^{\epsilon,\delta})
+\int_{0}^{t}D_x\Psi(x^{\epsilon,\delta}_r,y_r^{\epsilon,\delta})g(x^{\epsilon,\delta}_r,y_r^\epsilon)\dd r\big.\\
&\quad\big.+\int_{0}^{t}D_x\Psi(x^{\epsilon,\delta}_r,y_r^{\epsilon,\delta})f(x^{\epsilon,\delta}_r)\dd B_r^H
+\int_{0}^{t}\mathcal{L}(y_r^{\epsilon,\delta})\Psi(x^{\epsilon,\delta}_r, y_r^{\epsilon,\delta})\dd r\big]
+\frac{\sqrt{\delta}}{\sqrt{\epsilon}}M_t^\epsilon\\
&=:N_t^\epsilon+\frac{\sqrt{\delta}}{\sqrt{\epsilon}}M_t^\epsilon.
\end{aligned}
\end{equation}
By the same approach as in Theorem \ref{clt_Result} and under the assumption that $\delta=o(\epsilon)$, it is straightforward to see that no extra Gaussian term appears.
\end{itemize}
\end{remark}

\section{Acknowledgement} 
The first author was supported by JSPS Grant-in-Aid for JSPS Fellows (Grant No.
25KF0152). The second author was supported by the Key International (Regional) Cooperative Research Projects of the NSF of China (Grant 12120101002).

\section{Appendix A: Young-It\^o Formula}
\begin{proposition}
Consider the Young-It\^o SDE as below
\begin{equation}\label{Young-Ito-SDE}
\dd x_t =g(x_t) \dd t+ f(t, x_t) \dd B^H_t + \sigma(x_t) \dd W_t,\ x_0 = x.
\end{equation}
Here, the integral with respect to fractional Brownian motion \(B^H\) is in the pathwise Young integrals by Wiener integrals that introduced in Section 2. The integral with respect to BM \(W\) is in the It\^o sense.
Assume \(f\) is of \(BC^2\), \(\sigma\) and \(g\) are both Lipschitz continuous and linear growth with respect to \(x\). Then, for any function \(\Psi\in C_b^2\), 
\begin{equation}\label{phi_epsilon}
\begin{aligned}
\Psi(x_t)&=\Psi(x_0)+\int_{0}^{t}D_x\Psi(x_r)g(x_r)\dd r+\int_{0}^{t}D_x\Psi(x_r)f(r, x_r)\dd B_r^H\\
&+\int_{0}^{t}D_x\Psi(x_r)\sigma(x_r)\dd W_r+\frac{1}{2}\int_{0}^{t}D^2_{xx}\Psi(x_r)a(x_r)\dd r
\end{aligned}
\end{equation}
holds. Here, let \(a:=\sigma\sigma^*\).
\end{proposition}

{Proof}. By leveraging the Taylor's expansion, that is for any partition \(\mathcal{P}\) of time interval \([0,t]\), 
\begin{equation}\label{Taylor}
\begin{aligned}
\Psi(x_t)-\Psi(x_0)&=\sum_{[u,v]\in \mathcal{P}}\big(\langle D\Psi(x_u), \delta x_{u,v}\rangle+\frac{1}{2}\langle D^2\Psi(x_u), (\delta x_{u,v})^{\otimes 2}\rangle +o(|\delta x_{u,v}|^\gamma) \big)\cr
&=:\sum_{[u,v]\in \mathcal{P}} J_{u,v} +o(| \mathcal{P}|^{\alpha\gamma-1}).
\end{aligned}
\end{equation}
Note that \(\gamma>3\). Note that the sample path of FBM $B^H$ is of $\beta$-H\"older continuous almost surely via the Kolmogorov continuity theorem with $\beta<H$. Then recall the result that \(x\in \mathcal{B}_{\alpha, p}\) with \(1-\beta<\alpha<\frac{1}{2}\) which is an extension of \cite[Theorem 4.6]{Hairer-Li}, so \(\alpha\gamma-1>0\) (Specially, if $\sigma=0$, it has $1-\beta<\alpha<\beta$). Then we aim to establish the fact that 
\begin{equation}\label{Taylor-1}
\Psi(x_t)-\Psi(x_0)=L_2-\lim_{| \mathcal{P}|\to 0} \sum_{[u,v]\in \mathcal{P}} J_{u,v}
\end{equation}
Then, rewrite \(J_{u,v}\) as follows,
\begin{equation}\label{Taylor-2}
\begin{aligned}
J_{u,v}&=\langle D\Psi(x_u), g(x_u)(v-u)\rangle+\langle D\Psi(x_u), \int_{u}^{v}f(r, x_u)\dd B_{r}^H\rangle \\
&\quad+\langle D\Psi(x_u), \sigma(x_u)\delta W_{u,v}\rangle+\frac{1}{2}\langle D^2\Psi(x_u), (\int_{u}^{v}f(r, x_u)\dd B_{r}^H)^{\otimes 2}\rangle \\
&\quad+\frac{1}{2}\langle D^2\Psi(x_u), (\sigma(x_u)\delta W_{u,v})^{\otimes 2}\rangle +\frac{1}{2}\langle D^2\Psi(x_u), (g(x_u)(v-u))^{\otimes 2}\rangle \\
&\quad+\langle D^2\Psi(x_u), (\int_{u}^{v}f(r, x_u)\dd B_{r}^H)\otimes(\sigma(x_u)\delta W_{u,v})\rangle\\
&\quad+\langle D^2\Psi(x_u), (\int_{u}^{v}f(r, x_u)\dd B_{r}^H)\otimes(g(x_u)(v-u))\rangle \\
&\quad+\langle D^2\Psi(x_u), (\sigma(x_u)\delta W_{u,v})\otimes(g(x_u)(v-u))\rangle\\
&=:J^1_{u,v}+\cdots+J^9_{u,v}.
\end{aligned}
\end{equation}
We first focus on the term \(J^1_{u,v}\). First fix any \(\omega\in \Omega\), the convergence will directly follow from the definition of Riemann integral by \({| \mathcal{P}|\to 0}\), then convergence almost surely is obtained. With the dominated convergence theorem, the \(L_2\) convergence takes place.

Next, we shall to give the convergence of $\sum_{[u,v]\in \mathcal{P}} J^2_{u,v}$. Then we leverage the stochastic sewing lemma, \cite[Lemma 3.10]{Hairer-Li}. Precisely speaking, if \(f\in \mathcal{C}^{-\kappa,\gamma}\) for some \(\kappa,\gamma>0\) such that \(\eta>1/2\) and \(\bar \eta>1\). Then it has \(A_{u,v}:=\int_{u}^{v}f(r, x_u)\dd B_{r}^H\) belongs to \(H_\eta^p \cap \bar H_{\bar \eta}^p\) with \(\eta>1/2\) and \(\bar \eta>1\), and \(\|I_{u,v}(A)\|_p \lesssim |u-v|^\eta\) and \(\|\mathbb{E}(I_{u,v}(A)-A_{u,v}|\mathcal{F}_u)\|_p \lesssim |u-v|^{\bar\eta}\). So that we have \(\|J^2_{u,v}\|_p \lesssim |u-v|^\eta\) with \(\eta>1/2\), then by the H\"older inequality, assumption $\Psi\in C^2$ and fact that \(x\in \mathcal{B}_{\alpha, p}\) with \(1-\beta<\alpha<\beta\),
\begin{equation}\label{dB_t}
\begin{aligned}
&\|\sum_{[u,v]\in \mathcal{P}}(\int_{u}^v\langle D\Psi(x_r), f(r, x_r)\dd B_{r}^H\rangle-J^2_{u,v})\|^2_2\\
&\lesssim|\mathcal{P}|^{-1}\sum_{[u,v]\in \mathcal{P}}\|\mathbb{E}(\int_{u}^v\langle D\Psi(x_r), f(r, x_r)\dd B_{r}^H\rangle-J^2_{u,v}) |\mathcal{F}_u)\|^2_2\lesssim |\mathcal{P}|^{2\bar \eta-2}.
\end{aligned}
\end{equation}
Then \(\bar \eta-1>0 \) implies that 
\(\int_{0}^t\langle D\Psi(x_r), f(r, x_r)\dd B_{r}^H\rangle=L_p-\lim_{| \mathcal{P}|\to 0}\sum_{[u,v]\in \mathcal{P}} J^2_{u,v}\).

Next, we solve the convergence relate to the term \(J^3\). Consider a simple process \(\tilde \Phi\) which is defined by \(D\Phi(x_r)\) for any \(r\in[u,v]\). Then the It\^o integral for \(\tilde \Phi\) is defined in the sense that \(\sum_{[u,v]\in \mathcal{P}}J^3_{u,v}\). Then the convergence for \(L_2-\lim\sum_{[u,v]\in \mathcal{P}} J^3_{u,v}\) follows from the following convergence 
\begin{equation}\label{dW}
\lim_{| \mathcal{P}|\to 0} \mathbb{E}\big[\int_0^t(\tilde \Psi_s\sigma_s-D\Psi_s \sigma_s)^2\dd s\big] =0
\end{equation}
which is easy to verify. 

Next, we show the convergence of $\sum_{[u,v]\in \mathcal{P}} J^4_{u,v}$. Directly, it has
\begin{eqnarray}\label{J4}
\mathbb{E}[|J^4_{u,v}|^2]&=& \mathbb{E}[\mathbb{E}[|J^4_{u,v}|^2|\mathcal{F}_u]]\cr
&=& \mathbb{E}[|D_x\Psi(x_u)|^2\mathbb{E}[|(\int_{u}^{v}f(r, x_u)\dd B_{r}^H)^{\otimes 2}|^2|\mathcal{F}_u]]
\end{eqnarray}
Due to the \cite[Lemma 3.4]{Hairer-Li}, so we have \(\mathbb{E}[\|(\int_{u}^{v}f(r, x_u)\dd B_{r}^H)^{\otimes 2}\|^p|\mathcal{F}_u] \lesssim |u-v|^{2\eta p}\) with \(H>\eta>1/2\), then by result that \(x\in \mathcal{B}_{\alpha, p}\) with $p>1$ and \(1-\beta<\alpha<\frac{1}{2}\), the H\"older inequality and condition that $\Psi\in C_2$,\\
$$
\lim_{| \mathcal{P}|\to 0}\|\sum_{[u,v]\in \mathcal{P}}J^4_{u,v}\|^2_2\lesssim\lim_{| \mathcal{P}|\to 0}| \mathcal{P}|^{-1}\sum_{[u,v]\in \mathcal{P}}\|J^4_{u,v}\|^2_2 \lesssim \lim_{| \mathcal{P}|\to 0}|\mathcal{P}|^{4\eta -2}=0.$$
Then we estimate the term \(J^5_{u,v}\). We denote that \(\mathbb{E}^B\) and \(\mathbb{E}^W\) are the expectation with respect to \(B\) and \(W\) respectively, so it has that $\mathbb{E}= \mathbb{E}^B\times \mathbb{E}^W$. Then it has
\begin{equation}\label{dW_t^2}
\begin{aligned}
&\mathbb{E}^W\big[\sum_{[u,v]\in \mathcal{P}}(J^5_{u,v}-\frac{1}{2}\langle D^2\Psi(x_u), \sigma(x_u) (v-u)\rangle)\big]^2\cr
&=\mathbb{E}^W\big[\sum_{[u,v]\in \mathcal{P}}(\langle D^2\Psi(x_u), (\sigma(x_u)\delta W_{u,v})^{\otimes 2}\rangle-\langle D^2\Psi(x_u), \sigma(x_u) (v-u)\rangle) \big]^2\cr
&=\sum_{[u,v]\in \mathcal{P}}\mathbb{E}^W\big[\mathbb{E}^W_u\big(\langle D^2\Psi(x_u), a(x_u)( (\delta W_{u,v})^{\otimes 2}- (v-u) )\rangle \big)\big]^2
\end{aligned}
\end{equation}
Here the second equality comes from the independence between \(W_{s,t}\) and \(W_{s',t'}\) for different time interval \([s,t]\) and \([s't']\). Then, we have \(\sum_{[u,v]\in \mathcal{P}}J^5_{u,v}\) converges to \(\int_0^t \langle D^2\Phi(x_r), a(x_r)\rangle \dd r\) in the \(L_2\) sense.
Meanwhile, we have that \(\sum_{[u,v]\in \mathcal{P}}J^6_{u,v}\lesssim | \mathcal{P}|\to 0\).

Next, we consider the term \(J^7\), with some straight computation, 
\begin{equation}\label{dBdW_t}
\begin{aligned}
\big[\sum_{[u,v]\in \mathcal{P}}J^7_{u,v}\big]^2
&=\big[\sum_{[u,v]\in \mathcal{P}}\langle D^2\Psi(x_u), (\int_{u}^{v}f(r, x_u)\dd B_{r}^H)\otimes(\sigma(x_u)\delta W_{u,v})\rangle \big]^2\\
&=\sum_{[u,v]\in \mathcal{P}}\langle D^2\Psi(x_u), (\int_{u}^{v}f(r, x_u)\dd B_{r}^H)\otimes(\sigma(x_u)\delta W_{u,v})\rangle ^2\\
&\quad+\sum_{\substack{[u,v],[u',v']\in \mathcal{P}\\ [u,v]\cap[u',v']=\emptyset}}\big[\langle D^2\Psi(x_u), (\int_{u}^{v}f(r, x_u)\dd B_{r}^H)\otimes(\sigma(x_u)\delta W_{u,v})\rangle\big.\\
&\qquad\qquad\qquad\big.\times\langle D^2\Psi(x_u'), (\int_{u'}^{v'}f(r, x_u')\dd B_{r}^H)\otimes(\sigma(x_u')\delta W_{u',v'})\rangle \big]\\
&=:\sum_{[u,v]\in \mathcal{P}}J^{7,1}_{u,v}+\sum_{\substack{[u,v],[u',v']\in \mathcal{P}\\ [u,v]\cap[u',v']=\emptyset}}J^{7,2}_{u,v,u'v'}.
\end{aligned}
\end{equation}
Since \(\mathbb{E}^W[\sum_{\substack{[u,v],[u',v']\in \mathcal{P}\\ [u,v]\cap[u',v']=\emptyset}}J^{7,2}_{u,v,u'v'}]=0\), it remains to estimate the first convergence.
Then, due to the result that \(A_{u,v}:=\int_{u}^{v}f(r, x_u)\dd B_{r}^H\) belongs to \(H_\eta^p \cap \bar H_{\bar \eta}^p\) with \(\eta>1/2\) and \(\bar \eta>1\),
\begin{equation}\label{dBdW_t1}
\begin{aligned}
\mathbb{E}^B\big[\mathbb{E}^W\big[\sum_{[u,v]\in \mathcal{P}}J^{7,1}_{u,v}\big]\big]&=\sum_{[u,v]\in \mathcal{P}}\mathbb{E}^B\big[\mathbb{E}^W\big[\langle D^2\Psi(x_u), (\int_{u}^{v}f(r, x_u)\dd B_{r}^H)\otimes(\sigma(x_u)\delta W_{u,v})\rangle^2 \big]\big]\\
&\lesssim \sum_{[u,v]\in \mathcal{P}}|\mathcal{P}|^{ 2\eta+1}.
\end{aligned}
\end{equation}
Finally, the term related to \(J^8\) converges to 0 along the similar approach. Of course, the term \(J^9\) also converges to \(0\) and it is easy to verify. Just along the similar way and by applying $ \Psi\in C^2$ the result that \(x\in \mathcal{B}_{\alpha, p}\) with \(1-\beta<\alpha<\frac{1}{2}\), it is straightforward to see that $\lim_{| \mathcal{P}|\to 0} \sum_{[u,v]\in \mathcal{P}} J_{u,v}$ is uniformly bounded in $L^p$ sense for $p>2$, so the convergence in the $L^p$ sense is arrived.
The proof is completed. \qed

\section{Appendix B: Examples of fast one}
Assume that the fast process \(y^{\epsilon}\) satisfy the following It\^o SDE for \(t\in[0,T]\).
\begin{equation}\label{y_system}
\dd y^{\epsilon}_t =\frac{1}{\epsilon} b( x^{\epsilon}_t, y^{\epsilon}_t)\dd t + \frac{1}{{\sqrt \epsilon}}\sigma( x^{\epsilon}_t, y^{\epsilon}_t)\dd W_{t},\quad y^{\epsilon}_0=y\in \mathbb{R}^{m}.
\end{equation}
Here, \(x^\epsilon\) is the slow one in the system \eqref{simiplified_System}. This section aims to give an example for the Assumption \textbf{A4}.
\begin{example}\label{YHolder}
Let \(p\ge 2\). Under assumptions \textbf{H2}-\textbf{H3} and assume that \(b,\sigma\) are globally Lipschitz continuous and linear growth, we have 
\begin{itemize}
\item[(i).] There exists a constant $C>0$ such that for all $t\in [0,T]$,
\begin{equation}\label{y_Lp}
\sup_{0<\epsilon\le 1}\sup_{t\in[0,T]}\mathbb{E} [|y_t^\epsilon|^p ] \le C.
\end{equation}
\item[(ii).] Further assume that \(b(x,y):=-\Gamma y+\zeta(x,y)\) where \(\Gamma\) is a matrix. Let \(\theta<1/2\). There exists a constant \(C>0\) such that for all \(0<\epsilon\le 1\),
\begin{equation}\label{y_holder}
\mathbb{E}[\|y^\epsilon\|^p_{\theta-hld}] \le C \epsilon^{-\frac{p}{2}}.
\end{equation}
\end{itemize}
\end{example}
{Proof}. Firstly, we will show that (i).
By the It\^o's formula and differentiating on both sides, we have
\begin{equation}\label{y-3}
\begin{aligned}
\frac{\dd}{\dd t}\mathbb{E}[{| {y}_{t}^{\epsilon} |^p}] &\le \frac{p}{2\epsilon }\mathbb{E}\big[|y_t^\epsilon|^{p-2}\{2\langle {y}_{t}^{\epsilon},b( {x}_{t}^{\epsilon}, {y}_{t}^{\epsilon} )\rangle +(p-1)|\sigma( {x}_{t}^{\epsilon}, {y}_{t}^{\epsilon})|^2 \}\big]\cr
&\le \frac{p}{2\epsilon }\mathbb{E}\big[|y_t^\epsilon|^{p-2}\{-\beta_1|{y}_{t}^{\epsilon} |^2 +|{x}_{t}^{\epsilon} |^2 +C\}\big]\cr
&\le - \frac{p\beta_1}{4\epsilon }\mathbb{E}\big[|y_t^\epsilon|^{p}\big]+\frac{c}{\epsilon }\mathbb{E}\big[|x_t^\epsilon|^{p}+1\big]\cr
&\le - \frac{p\beta_1}{4\epsilon }\mathbb{E}\big[|y_t^\epsilon|^{p}\big]+\frac{c}{\epsilon }.
\end{aligned}
\end{equation}
where \(|\sigma (x,y)|\) is the Hilbert-Schmidt norm of the matrix \(\sigma(x,y)\).
Then by virtue of the Gronwall inequality, that is 
for any absolutely continuous function \(f:[0,T]\to \mathbb{R}\), for \(b,c\in \mathbb{R}\), it has
\(\frac{df_t}{dt} \le b f_t+c\),
then, \(f_t\le (f_0+\frac{c}{b})e^{bt}-\frac{c}{b}\).
So, take \(b=-p\beta_1/(4\epsilon)\) and \(c=c/\epsilon\),
\begin{equation}\label{y-4}
\mathbb{E}[{| \tilde {y}_{t}^{\epsilon} |^p}] 
\le |y_0|^p e^{-\frac{p\beta_1}{4\epsilon} t} + \frac{ c }{p\beta_1 }(1-e^{-\frac{p\beta_1}{4\epsilon} t})\le |y_0|^p+C_{\alpha, \beta_1,p}.
\end{equation}
Next, we intend to show (ii). We denote that \(\mathbb{E}^B\) and\(\mathbb{E}^W\) are the expectation with respect to \(B\) and \(W\) respectively, so it has that \(\mathbb{E}= \mathbb{E}^B\times \mathbb{E}^W\).
Then, we give the estimation for (ii) with respect to \(\mathbb{E}^W\).
Since we assume that $b(x,y):=-\Gamma y+\zeta(x,y)$ where $\Gamma$ is a matrix. We also assume that the function $\zeta$ is linear growth. Then, the fast equation is rewritten as:
\begin{equation}\label{y-8}
\dd y^{\epsilon}_t =-\frac{1}{\epsilon} \Gamma y^{\epsilon}_t \dd t+ \frac{1}{\epsilon} \zeta( x^{\epsilon}_t, y^{\epsilon}_t)\dd t + \frac{1}{{\sqrt \epsilon}}\sigma( x^{\epsilon}_t, y^{\epsilon}_t)\dd W_{t},\quad y^{\epsilon}_0=y\in \mathbb{R}^{m}.
\end{equation}
With the help of the Duhamel principle, we have that
\begin{equation}\label{y-8}
y^{\epsilon}_t =y e^{-\Gamma t/\epsilon}+\frac{1}{\epsilon} \int_0^te^{-\Gamma (t-s)/\epsilon}\zeta( x^{\epsilon}_s, y^{\epsilon}_s)\dd s +\frac{1}{\sqrt \epsilon} \int_0^te^{-\Gamma (t-s)/\epsilon}\sigma( x^{\epsilon}_s, y^{\epsilon}_s)\dd W_{s}.
\end{equation}
Then for \(0\le s<t\le T\), it has 
\begin{equation}\label{y-9}
\begin{aligned}
&y^{\epsilon}_t-y^{\epsilon}_s\cr
&=[I-e^{-\Gamma (t-s)/\epsilon}]y^{\epsilon}_s+\frac{1}{\epsilon} \int_s^te^{-\Gamma (t-r)/\epsilon}\zeta( x^{\epsilon}_r, y^{\epsilon}_r)\dd r +\frac{1}{\sqrt \epsilon} \int_s^te^{-\Gamma (t-r)/\epsilon}\sigma( x^{\epsilon}_r, y^{\epsilon}_r)\dd W_{r}\cr
&=:A_1(s,t)+A_2(s,t)+A_3(s,t).
\end{aligned}
\end{equation}
We firstly estimate \(A_2(s,t)\) and \(A_3(s,t)\). Recall that \(|e^{-\Gamma t y}|\le e^{-\gamma t |y|}\) with \(\gamma\) is the smallest eigenvalue for $\Gamma$. 
Next, for the first term in \eqref{y-9},
\begin{equation}\label{y-10}
\begin{aligned}
\mathbb{E}^W[|A_2(s,t)|]
&\le \frac{1}{\epsilon} \int_s^te^{-\gamma (t-r)/\epsilon}(1+\mathbb{E}^W[|x^{\epsilon}_r|]+\mathbb{E}^W[|y^{\epsilon}_r|])\dd r\\
&\le \frac{1}{\epsilon} \int_s^te^{-\gamma (t-r)/\epsilon}\dd r\sup_{r\in[s,t]} (1+\mathbb{E}^W[|x^{\epsilon}_r|]+\mathbb{E}^W[|y^{\epsilon}_r|])\\
&\le \frac{(t-r)^{1/2}}{\epsilon} \int_s^te^{-\gamma (t-r)/\epsilon}(t-r)^{-1/2}dr \big[ \sup_{r\in[s,t]}(1+\mathbb{E}^W[|x^{\epsilon}_r|]+\mathbb{E}^W[|y^{\epsilon}_r|])\big]\\
&\le \frac{(t-s)^{1/2}}{\sqrt \epsilon} \sup_{r\in[s,t]}(1+\mathbb{E}^W[|x^{\epsilon}_r|]+\mathbb{E}^W[|y^{\epsilon}_r|]).
\end{aligned}
\end{equation}
Note that we have \(\int_s^te^{-\gamma (t-r)/\epsilon}(t-r)^{-1/2}\dd r\le (\gamma/\epsilon)^{-\frac{1}{2}}\sup_{z\in\mathbb{R}} \int_0^z e^{-s}(z-s)\dd s\le C (\gamma/\epsilon)^{-\frac{1}{2}}\).

Then we intend to solve the second term in \eqref{y-9}. Denote that \(Y_\alpha(t)=\int_0^t(t-s)^{-\alpha} S(t-s) \Phi(s) d W_s\) for all \(t \geq 0\) with \(S^\epsilon(t-s) :=e^{-\gamma (t-s)/\epsilon}\) and \(\Phi(s):=\sigma(x^{\epsilon}_s, y^{\epsilon}_s)\). And 
\begin{equation}\label{y-11}
\begin{aligned}
& \frac{\sin a \pi}{\pi} \int_0^t(t-s)^{a-1} S^\epsilon(t-s) Y_\alpha(s) \dd s \\
& \quad=\frac{\sin a \pi}{\pi} \int_0^t(t-s)^{a-1} S^\epsilon(t-s)\left[\int_0^s(s-\sigma)^{-a} S^\epsilon(s-\sigma) \Phi(\sigma) \dd W_\sigma\right] \dd s \\
& \quad=\frac{\sin a \pi}{\pi} \int_0^t\left[\int_\sigma^t(t-s)^{a-1}(s-\sigma)^{-a} \dd s\right] S^\epsilon(t-\sigma) \Phi(\sigma) \dd W_\sigma
\end{aligned}
\end{equation}
By leveraging the result that 
\(\int_\sigma^t(t-s)^{a-1}(s-\sigma)^{-a} \dd s=\frac{\pi}{\sin a \pi}\) with\( 0 \leq \sigma \leq t\) and \( a \in(0,1)\),
then,
\begin{equation}\label{y-12}
\begin{aligned}
[A_3(s,t)]
\le\frac{1}{\sqrt \epsilon} \frac{\sin a \pi}{\pi} \int_s^t(t-r)^{\theta-1} e^{-\gamma (t-r)/\epsilon} \big[\int_0^r(r-\tau)^{-\alpha} e^{-\gamma (r-\tau)/\epsilon}\Phi(s) \dd W_\tau\big]\dd r.
\end{aligned}
\end{equation}
So that for \(\theta<1/2\), by the Burkholder-Davis-Gundy inequality,
\begin{equation}\label{y-13}
\mathbb{E}^W\big[\sup_{t\neq s}\frac{|A_3(s,t)|}{|t-s|^{\theta}}\big]^p
\le c \epsilon^{-\frac{p}{2}} \sup_{r\in[s,t]} [1+\mathbb{E}^W[x^{\epsilon}_r]^p+\mathbb{E}^W[y^{\epsilon}_r]^p].
\end{equation}
Next, we intend to estimate the term \(A_1(s,t)\). This term could be written as \(A_1(s,t)=[I-e^{-\Gamma (t-s)/\epsilon}](A_2(0,s)+A_3(0,s))\). Then, by applying \eqref{y-10}, 
\begin{equation}\label{y-14}
\begin{aligned}
&\mathbb{E}^W[|[I-e^{-\Gamma (t-s)/\epsilon}]A_2(0,s)|]\cr
&=\mathbb{E}^W[| [I-e^{-\Gamma (t-s)/\epsilon}]^{\theta}[I-e^{-\Gamma (t-s)/\epsilon}]^{1-\theta}A_2(0,s)|]\cr
&\le c(t-s)^\theta \epsilon^{-1-\theta}\int_0^s e^{-\gamma (s-r)/\epsilon}dr\sup_{r\in[0,t]}\mathbb{E}^W[1+|x^{\epsilon}_r|+|y^{\epsilon}_r|]\cr
&\le \epsilon^{-1-\theta}\epsilon^{\frac{1}{2}+\theta}\int_0^s (s-r)^{-\frac{1}{2}-\theta}dr\sup_{r\in[0,t]}\mathbb{E}^W[1+|x^{\epsilon}_r|+|y^{\epsilon}_r|]\cr
&\le \epsilon^{-\frac{1}{2}}(t-s)^{\frac{1}{2}}\sup_{r\in[0,t]}\mathbb{E}^W[1+|x^{\epsilon}_r|+|y^{\epsilon}_r|].
\end{aligned}
\end{equation}
For \(\theta<1/2\), with Burkholder-Davis-Gundy inequality,
\begin{equation}\label{y-13}
\mathbb{E}^W\big[\sup_{t\neq s}\frac{|[I-e^{-\Gamma (t-s)/\epsilon}]A_3(s,t)|}{|t-s|^{\theta}}\big]^p
\le c \epsilon^{-\frac{p}{2}} \sup_{r\in[s,t]} [1+\mathbb{E}^W[x^{\epsilon}_r]^p+\mathbb{E}^W[y^{\epsilon}_r]^p].
\end{equation}
Then by taking the expectation with respect to \(B\), the result in (ii) is arrived. 
This proof is completed.\qed

\bibliography{bibliography}

@article{Gehring-Li,
   title={Functional limit theorems for {V}olterra processes and applications to homogenization*},
   volume={35},
   ISSN={1361-6544},
   url={http://dx.doi.org/10.1088/1361-6544/ac4818},
   DOI={10.1088/1361-6544/ac4818},
   number={4},
   journal={Nonlinearity},
   publisher={IOP Publishing},
   author={Gehringer, Johann and Li, Xue-Mei  and Sieber, Julian },
   year={2022},
   month=mar, 
   pages={1521–1557} }

@article{rockner2021strong,
  title={Strong convergence order for slow--fast {M}cKean--{V}lasov stochastic differential equations},
  journal={Annales de l’Institut Henri Poincar{\'e}-Probabilit{\'e}s et Statistiques},
  volume={57},
  url={https://doi.org/10.1214/20-AIHP1087},
  DOI={10.1214/20-AIHP1087},
  number={4},
  number={1},
  author={R{\"o}ckner, Michael and Sun, Xiaobin and Xie, Yingchao},
  pages={547--576},
  year={2021},
  month=feb
}

@article{1999Decreusefond,
  title={Stochastic analysis of the fractional Brownian motion},
  author={Decreusefond, Laurent and {\"U}st{\"u}nel, Ali S},
  journal={Potential analysis},
  volume={10},
  number={2},
  pages={177--214},
  year={1999},
  publisher={Springer}
}

@article{mandelbrot1968fractional,
  title={Fractional Brownian motions, fractional noises and applications},
  author={Mandelbrot, Benoit B and Van Ness, John W},
  journal={SIAM review},
  volume={10},
  number={4},
  pages={422--437},
  year={1968},
  publisher={SIAM}
}

@article{Nualart03092008,
author = {Guerra, Joao and Nualart, David},
title = {Stochastic Differential Equations Driven by Fractional {B}rownian Motion and Standard {B}rownian Motion},
journal = {Stochastic Analysis and Applications},
volume = {26},
number = {5},
pages = {1053--1075},
year = {2008},
publisher = {Taylor \& Francis},
doi = {10.1080/07362990802286483},
URL = {https://doi.org/10.1080/07362990802286483},
eprint = { https://doi.org/10.1080/07362990802286483},
month = Aug
}

@article{Hairer-Li,
   title={Averaging dynamics driven by fractional {B}rownian motion},
   volume={48},
   ISSN={0091-1798},
   URL ={http://dx.doi.org/10.1214/19-AOP1408},
   DOI={10.1214/19-aop1408},
   number={4},
   journal={The Annals of Probability},
   publisher={Institute of Mathematical Statistics},
   author={Hairer, Martin and Li, Xue-Mei},
   year={2020},
   month=jul 
   }

@article{li2025fluctuations,
  title={Fluctuations from a random fractional averaging limit},
  author={Li, Xue-Mei and Piernot, Colin and Sobczak, Szymon and Ying, Kexing},
  journal={arXiv preprint arXiv:2512.08621},
  year={2025},
  month=dec 
}

@article{Hairer:22,
   title={Generating Diffusions with Fractional {B}rownian Motion},
   volume={396},
   ISSN={1432-0916},
   url={http://dx.doi.org/10.1007/s00220-022-04462-2},
   DOI={10.1007/s00220-022-04462-2},
   number={1},
   journal={Communications in Mathematical Physics},
   publisher={Springer Science and Business Media LLC},
   author={Hairer, Martin and Li, Xue-Mei},
   year={2022},
   month=aug, 
   pages={91--141} 
  }

@article{lamperti1962convergence,
  title={On convergence of stochastic processes},
  author={Lamperti, John},
  journal={Transactions of the American Mathematical Society},
  volume={104},
  number={3},
  pages={430--435},
  year={1962},
  month=sep, 
  publisher={JSTOR}
}

@article{hong2023central,
  title={Central limit type theorem and large deviation principle for multi-scale {M}cKean--{V}lasov SDEs},
  author={Hong, Wei and Li, Shihu and Liu, Wei and Sun, Xiaobin},
  journal={Probability Theory and Related Fields},
  volume={187},
  number={1-2},
  pages={133--201},
  year={2023},
  publisher={Springer},
  month = jun,
}

@article{Pardoux2001,
author = {Pardoux , {\`E}.  and Veretennikov, A. Yu. },
title = {{On the {P}oisson Equation and Diffusion Approximation. I}},
volume = {29},
journal = {The Annals of Probability},
number = {3},
publisher = {Institute of Mathematical Statistics},
pages = {1061 -- 1085},
year = {2001},
doi = {10.1214/aop/1015345596},
URL = {https://doi.org/10.1214/aop/1015345596},
month = jul,
}

@article{Pardoux2003,
author = {Pardoux , {\`E}.  and Veretennikov, A. Yu. },
title = {{On {P}oisson equation and diffusion approximation 2}},
volume = {31},
journal = {The Annals of Probability},
number = {3},
publisher = {Institute of Mathematical Statistics},
pages = {1166 -- 1192},
year = {2003},
doi = {10.1214/aop/1055425774},
URL = {https://doi.org/10.1214/aop/1055425774},
month = jul,
}

@article{Pardoux2005,
author = {Pardoux , {\`E}.  and Veretennikov, A. Yu. },
title = {{On the {P}oisson equation and diffusion approximation 3}},
volume = {33},
journal = {The Annals of Probability},
number = {3},
publisher = {Institute of Mathematical Statistics},
pages = {1111 -- 1133},
year = {2005},
doi = {10.1214/009117905000000062},
URL = {https://doi.org/10.1214/009117905000000062},
month = may,
}

@book{stroock2007multidimensional,
  title={Multidimensional diffusion processes},
  author={Stroock, Daniel W. and Varadhan, S. R. Srinivasa},
  year={2007},
  publisher={Springer}
}

@article{Khoa2020,
author = {Khoa, L{\^e}},
title = {{A stochastic sewing lemma and applications}},
volume = {25},
journal = {Electronic Journal of Probability},
number = {none},
publisher = {Institute of Mathematical Statistics and Bernoulli Society},
pages = {1--55},
year = {2020},
doi = {10.1214/20-EJP442},
URL = {https://doi.org/10.1214/20-EJP442},
month = may,
}

@article{Pei2023,
author = {Pei, Bin and Inahama, Yuzuru and Xu, Yong},
title = {{Averaging principles for mixed fast-slow systems driven by fractional {B}rownian motion}},
volume = {63},
journal = {Kyoto Journal of Mathematics},
number = {4},
publisher = {Duke University Press},
pages = {721 -- 748},
year = {2023},
doi = {10.1215/21562261-2023-0001},
URL = {https://doi.org/10.1215/21562261-2023-0001},
month = nov,
}

@article{Pei2024,
author = {Pei, Bin and Schmalfuss, Bjoern and Xu, Yong},
title = {Almost Sure Averaging for Evolution Equations Driven by Fractional {B}rownian Motions},
journal = {SIAM Journal on Applied Dynamical Systems},
volume = {23},
number = {4},
pages = {2807-2852},
year = {2024},
doi = {10.1137/23M1554448},
URL = {https://doi.org/10.1137/23M1554448},
month = dec,
}

@article{PEI2021202,
title = {Averaging principle for fast-slow system driven by mixed fractional {B}rownian rough path},
journal = {Journal of Differential Equations},
volume = {301},
pages = {202-235},
year = {2021},
issn = {0022-0396},
doi = {https://doi.org/10.1016/j.jde.2021.08.006},
url = {https://www.sciencedirect.com/science/article/pii/S0022039621005040},
author = {Pei, Bin and Inahama, Yuzuru and Xu, Yong},
month = nov
}

@article{Yuzuru2025,
author = {Inahama, Yuzuru},
title = {{Averaging principle for slow-fast systems of rough differential equations via controlled paths}},
volume = {77},
journal = {Tohoku Mathematical Journal},
number = {2},
publisher = {Tohoku University, Mathematical Institute},
pages = {189 -- 227},
year = {2025},
doi = {10.2748/tmj.20230525},
URL = {https://doi.org/10.2748/tmj.20230525},
month =jun
}

@article{inahama2025averaging,
  title={Averaging principle for rough slow-fast systems of level 3},
  author={Inahama, Yuzuru},
  journal={arXiv preprint arXiv:2504.03110},
  year={2025},
}

@article{LI2025104683,
title = {Averaging principle for semilinear slow–fast rough partial differential equations},
journal = {Stochastic Processes and their Applications},
volume = {188},
pages = {104683},
year = {2025},
issn = {0304-4149},
doi = {https://doi.org/10.1016/j.spa.2025.104683},
url = {https://www.sciencedirect.com/science/article/pii/S0304414925001243},
author = {Li, Miaomiao and Li, Yunzhang and Pei, Bin and Xu, Yong},
month =oct
}

@article{li2025averaging,
  title={Averaging principle for slow-fast systems of {PDE}s with rough drivers},
  author={Li, Miaomiao and Pei, Bin and Xu, Yong and Yue, Xiaole},
  journal={arXiv preprint arXiv:2510.22269},
  year={2025}
}

@article{gehringer2020diffusive,
  title={Diffusive and rough homogenisation in fractional noise field},
  author={Gehringer, Johann and Li, Xue-Mei},
  journal={arXiv preprint arXiv:2006.11544},
  year={2020}
}

@article{gehringer2022functional,
  title={Functional limit theorems for the fractional {O}rnstein--{U}hlenbeck process},
  author={Gehringer, Johann and Li, Xue-Mei},
  journal={Journal of Theoretical Probability},
  volume={35},
  number={1},
  pages={426--456},
  year={2022},
  publisher={Springer},
  month = mar
}

@article{Li_2008,
doi = {10.1088/0951-7715/21/4/008},
url = {https://doi.org/10.1088/0951-7715/21/4/008},
year = {2008},
month = {mar},
publisher = {},
volume = {21},
number = {4},
pages = {803},
author = {Li, Xue-Mei},
title = {An averaging principle for a completely integrable stochastic Hamiltonian system},
journal = {Nonlinearity},
}

@article{Bezemek-Spiliopoulos, 
author = {Bezemek, Zachary William and Spiliopoulos, Konstantinos},
title = {Rate of homogenization for fully-coupled {M}cKean–{V}lasov SDEs},
journal = {Stochastics and Dynamics},
volume = {23},
number = {02},
pages = {2350013},
year = {2023},
doi = {10.1142/S0219493723500132},
URL = { https://doi.org/10.1142/S0219493723500132},
eprint = { https://doi.org/10.1142/S0219493723500132}
}

@article{Li-Wu-Xie,
author = {Li, Yun and Wu, Fuke and Xie, Longjie},
title = {Poisson Equation on Wasserstein Space and Diffusion Approximations for Multiscale {M}cKean–{V}lasov Equation},
journal = {SIAM Journal on Mathematical Analysis},
volume = {56},
number = {2},
pages = {1495-1524},
year = {2024},
doi = {10.1137/22M1536856},
URL = { https://doi.org/10.1137/22M1536856},
eprint = { https://doi.org/10.1137/22M1536856},
month = mar
}

@article{hairer2004periodic,
  title={Periodic homogenization for hypoelliptic diffusions},
  author={Hairer, Martin and Pavliotis, GA},
  journal={Journal of Statistical Physics},
  volume={117},
  number={1},
  pages={261--279},
  year={2004},
  publisher={Springer},
  month=oct
}

@article{BAILLEUL20174894,
title = {Rough flows and homogenization in stochastic turbulence},
journal = {Journal of Differential Equations},
volume = {263},
number = {8},
pages = {4894-4928},
year = {2017},
issn = {0022-0396},
doi = {https://doi.org/10.1016/j.jde.2017.06.006},
url = {https://www.sciencedirect.com/science/article/pii/S0022039617303121},
author = {I. Bailleul and R. Catellier},
month = oct
}

@article{Bailleul2019,
  title={Rough flows},
  author={Isma\"el Bailleul and Sebastian Riedel},
  journal={Journal of the Mathematical Society of Japan},
  volume={71},
  number={3},
  pages={915-978},
  year={2019},
  doi={10.2969/jmsj/80108010}
}

@book{jacod2003limit,
  title={Limit theorems for stochastic processes},
  author={Jacod, Jean and Shiryaev, Albert},
  volume={288},
  year={2003},
  publisher={Springer Science \& Business Media}
}

@book{komorowski2012fluctuations,
  title={Fluctuations in {M}arkov processes: time symmetry and martingale approximation},
  author={Komorowski, Tomasz and Landim, Claudio and Olla, Stefano},
  volume={345},
  year={2012},
  publisher={Springer Science \& Business Media}
}

@article{kipnis1986central,
  title={Central limit theorem for additive functionals of reversible {M}arkov processes and applications to simple exclusions},
  author={Kipnis, Claude and Varadhan, S. R. Srinivasa},
  journal={Communications in Mathematical Physics},
  volume={104},
  number={1},
  pages={1--19},
  year={1986},
  publisher={Springer},
  month=mar
}

@article{HONG2025405,
title = {Diffusion approximation for multi-scale {M}cKean-{V}lasov SDEs through different methods},
journal = {Journal of Differential Equations},
volume = {414},
pages = {405-454},
year = {2025},
issn = {0022-0396},
doi = {https://doi.org/10.1016/j.jde.2024.09.012},
url = {https://www.sciencedirect.com/science/article/pii/S0022039624005928},
author = {Wei Hong and Shihu Li and Xiaobin Sun},
month =jan
}

@article{sun2025diffusion,
  title={Diffusion Approximation for Slow-Fast SDEs with State-Dependent Switching},
  author={Sun, Xiaobin and Wang, Jue and Xie, Yingchao},
  journal={arXiv preprint arXiv:2503.08047},
  year={2025}
}

@article{li2023functional,
  title={Functional law of large numbers and central limit theorem for slow-fast {M}cKean-{V}lasov equations.},
  author={Li, Yun and Xie, Longjie},
  journal={Discrete \& Continuous Dynamical Systems-Series S},
  volume={16},
  number={5},
  year={2023},
  month=jan 
}

@article{nourdin2016multivariate,
  title={Multivariate central limit theorems for averages of fractional {V}olterra processes and applications to parameter estimation},
  author={Nourdin, Ivan and Nualart, David and Zintout, Rola},
  journal={Statistical Inference for Stochastic Processes},
  volume={19},
  number={2},
  pages={219--234},
  year={2016},
  publisher={Springer}
}

@book{nourdin2012normal,
  title={Normal approximations with {M}alliavin calculus: from {S}tein's method to universality},
  author={Nourdin, Ivan and Peccati, Giovanni},
  volume={192},
  year={2012},
  publisher={Cambridge University Press}
}

@book{Nualart,
  title={The Malliavin calculus and related topics},
  author={Nualart, David},
  year={2006},
  publisher={Springer}
}

@book{Friz-Hairer,
  title={A Course on Rough Paths: With an Introduction to Regularity Structures},
  author={Friz, Peter K. and Hairer, Martin},
  year={2020},
  publisher={Springer Nature}
}

@book{Friz-Victoir,
  title={Multidimensional stochastic processes as rough paths: theory and applications},
  author={Friz, Peter K and Victoir, Nicolas B},
  volume={120},
  year={2010},
  publisher={Cambridge University Press}
}

@article{Bourguin-Spiliopoulos,
title = {Quantitative fluctuation analysis of multiscale diffusion systems via {M}alliavin calculus},
journal = {Stochastic Processes and their Applications},
volume = {180},
pages = {104524},
year = {2025},
issn = {0304-4149},
doi = {https://doi.org/10.1016/j.spa.2024.104524},
url = {https://www.sciencedirect.com/science/article/pii/S0304414924002321},
author = {Bourguin,  Solesne and Spiliopoulos, Konstantinos},
month = feb
}

@article{Kelly-Melbourne,
author = {Kelly, David and Melbourne, Ian},
title = {{Smooth approximation of stochastic differential equations}},
volume = {44},
journal = {The Annals of Probability},
number = {1},
publisher = {Institute of Mathematical Statistics},
pages = {479 -- 520},
year = {2016},
doi = {10.1214/14-AOP979},
URL = {https://doi.org/10.1214/14-AOP979},
month = jan
}

@article{bourguin-Gailus,
  title={Typical dynamics and fluctuation analysis of slow--fast systems driven by fractional {B}rownian motion},
  author={Bourguin, Solesne and Gailus, Siragan and Spiliopoulos, Konstantinos},
  journal={Stochastics and Dynamics},
  volume={21},
  number={07},
  pages={2150030},
  year={2021},
  publisher={World Scientific}
}

@article{Djurdjevac-Helena,
author = {Djurdjevac, Ana and Kremp, Helena and Perkowski, Nicolas},
title = {Rough homogenization for {L}angevin dynamics on fluctuating {H}elfrich surfaces},
journal = {Stochastic Analysis and Applications},
volume = {43},
number = {3},
pages = {423--445},
year = {2025},
publisher = {Taylor \& Francis},
doi = {10.1080/07362994.2025.2505736},
URL = { https://doi.org/10.1080/07362994.2025.2505736},
eprint = { https://doi.org/10.1080/07362994.2025.2505736},
month = jun
}

@article{GU20121069,
title = {Random homogenization and convergence to integrals with respect to the {R}osenblatt process},
journal = {Journal of Differential Equations},
volume = {253},
number = {4},
pages = {1069-1087},
year = {2012},
issn = {0022-0396},
doi = {https://doi.org/10.1016/j.jde.2012.05.007},
url = {https://www.sciencedirect.com/science/article/pii/S0022039612002100},
author = {Gu, Yu and Bal, Guillaume},
month = aug
}

@article{LI-Gao,
title = {Averaging principle for slow-fast {SPDE}s driven by mixed noises},
journal = {Journal of Differential Equations},
volume = {430},
pages = {113209},
year = {2025},
issn = {0022-0396},
doi = {https://doi.org/10.1016/j.jde.2025.02.080},
url = {https://www.sciencedirect.com/science/article/pii/S0022039625002086},
author = {Li, Haoyuan and Gao, Hongjun and Qu, Shiduo},
month = jun
}

@ARTICLE{Rockner-Xie,
       author = {{R{\"o}ckner}, Michael and {Xie}, Longjie},
        title = "{Averaging Principle and Normal Deviations for Multiscale Stochastic Systems}",
      journal = {Communications in Mathematical Physics},
         year = 2021,
        month = may,
       volume = {383},
       number = {3},
        pages = {1889-1937},
          doi = {10.1007/s00220-021-04069-z},
archivePrefix = {arXiv},
       eprint = {2008.04822},
 primaryClass = {math.PR},
       adsurl = {https://ui.adsabs.harvard.edu/abs/2021CMaPh.383.1889R}
}

@article{ROCKNER-Xie23,
title = {Asymptotic behavior of multiscale stochastic partial differential equations with {H}\"older coefficients},
journal = {Journal of Functional Analysis},
volume = {285},
number = {9},
pages = {110103},
year = {2023},
issn = {0022-1236},
doi = {https://doi.org/10.1016/j.jfa.2023.110103},
url = {https://www.sciencedirect.com/science/article/pii/S0022123623002604},
author = {R\"ockner, Michael and Xie, Longjie and Yang, Li},
month = nov
}

@article{YANG2022107897,
title = {The central limit theorem for slow–fast systems with {L}\'evy noise},
journal = {Applied Mathematics Letters},
volume = {128},
pages = {107897},
year = {2022},
issn = {0893-9659},
doi = {https://doi.org/10.1016/j.aml.2021.107897},
url = {https://www.sciencedirect.com/science/article/pii/S0893965921004882},
author = {Yang, Xiaoyu and Xu, Yong and Wang, Ruifang and Jiao, Zhe},
month = jun
}

@article{WANG2013822,
title = {Slow manifold and averaging for slow–fast stochastic differential system},
journal = {Journal of Mathematical Analysis and Applications},
volume = {398},
number = {2},
pages = {822-839},
year = {2013},
issn = {0022-247X},
doi = {https://doi.org/10.1016/j.jmaa.2012.09.029},
url = {https://www.sciencedirect.com/science/article/pii/S0022247X12007615},
author = {Wang, Wei and Roberts, Anthony John},
month =feb
}

@book{lyons-Qian,
  title={System control and rough paths},
  author={Lyons, Terry and Qian, Zhongmin},
  year={2002},
  publisher={Oxford University Press}
}

@article{engel2021homogenization,
  title={Homogenization of coupled fast-slow systems via intermediate stochastic regularization},
  author={Engel, Maximilian and Gkogkas, Marios Antonios and Kuehn, Christian},
  journal={Journal of Statistical Physics},
  volume={183},
  number={2},
  pages={25},
  year={2021},
  publisher={Springer},
  month=apr
}

@article{KELLY20174063,
title = {Deterministic homogenization for fast-slow systems with chaotic noise},
journal = {Journal of Functional Analysis},
volume = {272},
number = {10},
pages = {4063-4102},
year = {2017},
issn = {0022-1236},
doi = {https://doi.org/10.1016/j.jfa.2017.01.015},
url = {https://www.sciencedirect.com/science/article/pii/S0022123617300162},
author = {Kelly, David and Melbourne, Ian},
month = may
}

\end{document}